\documentclass[11pt,reqno]{amsart}
\usepackage{geometry}
\usepackage{mathrsfs,cite}
\usepackage{amsthm}
\usepackage{amssymb}
\usepackage{graphicx}
\usepackage{color}
\usepackage{amsfonts}
\usepackage{amsmath} %%% provide equation structures
\usepackage{bm}           %%% blod math symbols \bm{}
\usepackage{mathabx}
\usepackage{tikz}
\usepackage{xcolor}
\usepackage{enumerate}

\newcommand{\R}{\mathbb{R}}

\newcommand{\dd}{{\rm d\hspace{0.1mm}}}
\newcommand{\dr}{{\rm d}}
\newcommand{\D}{\mathbb{D}}

\newcommand{\CS}{\mathcal{S}}
 \newcommand{\xx}{\mathbf{x}}
\newcommand{\gb}{\beta}

\newcommand{\gD}{\Delta} 
\newcommand{\po}{\partial} 
 
\newcommand{\ve}{\varepsilon} 
 \newcommand{\gl}{\lambda}

\newcommand{\del}{{\partial}}
\renewcommand{\O}{{\mathbf O}}
\renewcommand{\d}{\delta} \renewcommand{\l}{\lambda}
 
\renewcommand{\b}{\beta} 
\newcommand{\s}{\sigma} \newcommand{\g}{\gamma}

\newcommand{\bR}{ {\mathbb{R}}} %%%% real
\newcommand{\bRnpr}{{\mathbb{R}^{d-1}}} %%%% real (n-1)-dim
\newcommand{\bTpr}{{\mathbb{T}^{d-1}}} %%%%  (n-1)-dim flat torus
\newcommand{\Sn}{{\mathbb{S}^{d-1}}} %%%% sphere S^{n-1}
\newcommand{\xpr}{{\x}'}  %%%% x' = (x_2, \cdots, x_n).
\newcommand{\dist}{ \mbox{dist}}
\newcommand{\divg}{ \mbox{div\,}}

\newcommand{\setK}{{\mathcal K}_{M}}

\newcommand{\mapJ}{{J}}

\newcommand{\defd}{:=}
\newcommand{\Le}{\tilde{\mathcal L}}
\newcommand{\Trho}{\tilde\rho}
\newcommand{\epsP}{\sigma}
\newcommand{\bnNot}{{\bn_0}}
\newcommand{\Sqr}{ {(0, 1)^{d-1}}}
\newcommand{\RDom}{{\Omega_{\rm e}}}

\newcommand{\ExCFn}{{\mathcal {P}}}

\newcommand{\bea}{{\begin{eqnarray}}}
\newcommand{\eea}{{\end{eqnarray}}}

\numberwithin{equation}{section}
\theoremstyle{plain}
\newtheorem{definition}{Definition}[section]
\newtheorem{theorem}{Theorem}[section]

\newtheorem{corollary}[definition]{Corollary}
\newtheorem{remark}[definition]{Remark}

\newtheorem{problem}{Problem}[section]

\newcommand{\nnu}{{\boldsymbol \nu}}

\newcommand{\ttau}{{\boldsymbol \tau}}

\newcommand{\eps}{\varepsilon}
\newcommand{\ee}{\mathbf{e}}

\newcommand{\sS}{\mathcal S}

\newcommand{\zz}{\mathbf{z}}
\newcommand{\ff}{\mathbf{f}}
\newcommand{\kk}{\mathbf{k}}

\newcommand{\bb}{\mathbf{b}}
\newcommand{\uu}{\mathbf{u}}

\newcommand{\ga}{\alpha}
\newcommand{\gc}{\gamma}
\newcommand{\vph}{\varphi}
\newcommand{\vae}{\varepsilon}

\newcommand{\gd}{\delta}
\renewcommand{\O}{{O}}
\renewcommand{\d}{\delta}
\renewcommand{\l}{\lambda}

\renewcommand{\b}{\beta}

\newcommand{\yy}{\mathbf{y}}

\newcommand{\hs}{\hat{\sigma}}
\newcommand{\pw}{\partial W}

\newcommand \der{\partial}
\newcommand \vphi{\varphi}

\newcommand \gam{\gamma}
\newcommand \alp{\alpha}

\newcommand \ol{\overline}

\newcommand{\x}{\mathbf{x}}
\newcommand{\y}{\mathbf{y}}

\newcommand \Sonic{\Gamma_{\rm sonic}}
\newcommand \Shock{\Gamma_{\rm shock}}
\newcommand \Wedge{\Gamma_{\rm wedge}}

\newcommand \sonic{\Gamma_{\rm sonic}}
\newcommand \shock{\Gamma_{\rm shock}}

\newcommand \ivphi{\varphi_0}
\newcommand \iu{u_{10}}

\newcommand \irho{\rho_0}

\newcommand{\pSi}{\varphi}
\newcommand{\PtIncW}{{P_0}}
\newcommand{\PtUpL}{{P_1}}
\newcommand{\PtLwL}{{P_2}}
\newcommand{\PtLwR}{{P_3}}
\newcommand{\PtUpR}{{P_4}}
\newcommand \bPhi{\bar{\Phi}}
\newcommand{\xxi}{{\boldsymbol\xi}}
\newcommand{\Gso}{\Gamma_{\text{\rm sonic}}}

\newcommand{\Gsh}{\Gamma_{\text{\rm shock}}}

\newcommand{\bt}{{\boldsymbol\tau}}
\newcommand{\bn}{{\boldsymbol\nu}}

\newcommand{\mR}{\mathbb{R}}

\newcommand{\cxi}{{\xi_1}}
\newcommand{\ceta}{{\xi_2}}

\renewcommand{\d}{\delta}
\renewcommand{\l}{\lambda}

\renewcommand{\r}{\rho}
\newcommand{\z}{\zeta}

\renewcommand{\bb}{{\boldsymbol{\beta}}}
\newcommand \Symm{\Gamma_{\rm sym}}

\numberwithin{equation}{section}
\numberwithin{figure}{section}
\begin{document}
\title[Transonic Shock Waves and Free Boundary Problems]{Multidimensional Transonic
Shock Waves \\and Free Boundary Problems}
\author{Gui-Qiang G. Chen}
\address{Gui-Qiang G. Chen:$\,$ OxPDE, Mathematical Institute, University of Oxford,
Oxford, OX2 6GG, UK\\
chengq@maths.ox.ac.uk}
\author{Mikhail Feldman}
\address
{Mikhail Feldman: $\,$ Department of Mathematics, University of Wisconsin-Madison, WI 53706-1388, USA\\
         feldman@math.wisc.edu}

\begin{abstract}
We are concerned with free boundary problems arising from the analysis of multidimensional transonic
shock waves for the Euler equations in compressible fluid dynamics.
In this expository paper, we survey some recent developments in the analysis of multidimensional
transonic shock waves and corresponding free boundary problems for the compressible Euler equations
and related nonlinear partial differential equations (PDEs) of mixed type.
The nonlinear PDEs under our analysis include the steady Euler equations for potential flow, the steady full Euler equations,
the unsteady Euler equations for potential flow, and related nonlinear PDEs of mixed elliptic-hyperbolic type.
The transonic shock problems include the problem of steady transonic flow past solid wedges,
the von Neumann problem for shock reflection-diffraction, and the Prandtl-Meyer problem for unsteady
supersonic flow onto solid wedges.
We first show how these longstanding multidimensional
transonic shock problems can be formulated as free boundary problems
for the compressible Euler equations and related nonlinear PDEs of mixed type.
Then we present an effective nonlinear method and related ideas and techniques to solve these free boundary problems.
The method, ideas, and techniques should be useful to analyze
other longstanding and newly emerging free boundary problems for nonlinear PDEs.
\end{abstract}


\maketitle
\tableofcontents

\section{Introduction}
We are concerned with free boundary problems arising from the analysis of multidimensional transonic
shock waves for the Euler equations in compressible fluid dynamics.
The purpose of this expository paper is to survey some recent developments
in the analysis of multidimensional (M-D) transonic shock waves and
corresponding free boundary problems
for the Euler equations and related
nonlinear partial differential equations (PDEs) of mixed type.
We show how several M-D transonic shock problems can be formulated as
free boundary problems for the compressible Euler equations and related nonlinear PDEs of mixed-type,
and then present an efficient nonlinear method and
related ideas and techniques to solve these free boundary problems.

Shock waves are steep wavefronts, which are fundamental in
high-speed fluid flows
({\it e.g.}, \cite{BD,Busemann,Chapman,CF-book2018,CF,Da,GlimmMajda,Majda,Mises,VD,Neumann0,Whitham}).
Such flows are governed by the compressible Euler equations in fluid dynamics.
The time-dependent compressible Euler equations are a
second-order nonlinear wave equation for potential flow,
or a first-order nonlinear system of hyperbolic conservation laws
for full Euler flow ({\it e.g.}, \cite{Chen,CF-book2018,CF,Da}).
One of the main features of such nonlinear PDEs is that, no mater how smooth
the given initial data start with,
the solution develops singularity in a finite time to form shock waves (shocks, for short) generically,
so that the classical notion of solutions has to be extended to the notion of entropy solutions
in order to accommodate such discontinuity waves for physical variables,
{\it that is}, the weak solutions satisfying
the entropy condition that is consistent with the second law of thermodynamics
({\it cf}. \cite{CF-book2018,CF,Da,Lax}).

General entropy solutions involving shocks for such PDEs have extremely complicated and rich structures.
On the other hand, many fundamental problems in physics and engineering
concern steady solutions ({\it i.e.}, time-independent solutions)
or self-similar solutions ({\it i.e.}, the solutions depend only on the self-similar variables
with form $\frac{\x}{t}$ for the space variables $\x$ and time-variable $t$);
see \cite{CF-book2018,CF,Da,GlimmMajda} and the references cited therein.
Such solutions are governed by the steady or self-similar compressible Euler equations
for potential flow, or full Euler flow.
These governing PDEs in the new forms are time-independent and often are of mixed elliptic-hyperbolic type.

Mathematically, {\it M-D transonic shocks} are codimension-one discontinuity fronts in the solutions
of the steady or self-similar Euler equations and related nonlinear PDEs of mixed elliptic-hyperbolic type,
which separate two phases:
one of them is supersonic phase ({\it i.e.}, the fluid speed is larger than the sonic speed) which is hyperbolic;
the other is subsonic phase ({\it i.e.}, the fluid speed is smaller than the sonic speed) which is elliptic for potential flow,
or elliptic-hyperbolic composite for full Euler
flow ({\it i.e.}, elliptic equations coupled with some hyperbolic transport equations).
They are formed in many physical situations,
for example, by smooth supersonic flows or supersonic shock waves impinging onto solid wedges/cones or passing through de Laval nozzles,
around supersonic or near-sonic flying bodies, or other physical processes.
The mathematical analysis of shocks at least dates back to Stokes \cite{Stokes} and Riemann \cite{Riemann},
starting from the one-dimensional (1-D) case.
The mathematical understanding of M-D transonic shocks has been one of the most challenging and longstanding
scientific research directions ({\it cf}. \cite{CF-book2018,CC,Sxchen-book2020,CF,Da,GlimmMajda,Guderley}).
Such transonic shocks can be formulated as {\it free boundary problems} (FBPs) in the mathematical theory of
nonlinear PDEs involving mixed elliptic-hyperbolic type.

Generally speaking, a {\it free boundary problem} is a boundary value problem for a PDE or system of PDEs which is defined in a domain,
a part of whose boundary is {\it a priori unknown}; this part is accordingly named as a {\it free boundary}.
The mathematical problem is then to determine both the location of the free boundary and the solution
of the PDE/system in the resulting domain, which requires to combine analysis and geometry in sophisticated ways.
The mathematical analysis of FBPs is one of the most important research directions in the analysis of PDEs,
with wide applications across the sciences and real-world problems.
On the other hand, it is widely regarded as a truly challenging field of mathematics.
See \cite{CJK,Caff-Salsa,CSV,DH,Friedman,KinderlehrerNirenberg} and the references cited therein.

Transonic shock problems for steady or self-similar solutions are typically formulated as
boundary value problems for a nonlinear PDE or system of mixed elliptic-hyperbolic type,
whose type at a point is determined by the solution, as well as its gradient for some cases.
For a system, the type is more complicated and may be either hyperbolic
or mixed-composite elliptic-hyperbolic (also called {\it mixed}, for the sake of brevity when no confusion arises).
General solutions of such nonlinear PDEs can be nonsmooth and of complicated structures
({\it e.g.} \cite{CCY2,CCY3,CH,Chen,CF-book2018,GlimmMajda,KTa,LaxLiu,LZY,SCG,Serre,ZZ,Zhe}),
so that even the uniqueness issue has not been settled in many cases.
However, in many problems, especially those motivated by physical phenomena,
the expected structures of solutions are known from experimental/numerical results and underlying physics.
The solutions are expected to be piecewise smooth, with some hyperbolic/elliptic regions
separated by shocks, or sonic curves/surfaces
of continuous type-transition ({\it i.e.}, the type of equations changes continuously in the physical variables
such as the velocity, density, {\it etc.}).
In this paper, we present the problems in which the hyperbolic part of the solution is {\it a priori} known,
or can be determined separately from the elliptic part, in some larger regions.
Then the problem is reduced to determining the region in which the underlying PDE
is elliptic, with the transonic shock
as a part of its boundary and the elliptic solution in that region.
In other words, we need to solve a free boundary problem for the elliptic phase of the solution,
with the transonic shock as a free boundary.
Since the type of equations depends on the solution itself, the ellipticity in the region is a part of
the results to be established.
We remark that, in some other problems involving shocks, FBPs also need to be solved in order to find the hyperbolic part
of the solution, which is beyond the scope of this paper.

For several problems under our discussion below,
the PDEs involved are single second-order quasilinear PDEs,
whose coefficients and types (elliptic, hyperbolic, or mixed) depend on the gradient of the solution.
In the other problems, the PDEs are first-order nonlinear systems, whose types
are hyperbolic or composite-mixed elliptic-hyperbolic, and are determined by the solution only.
In all the problems,
the PDEs (or  parts of the systems) are expected to be elliptic for our solutions in the regions
determined by the free boundary problems.
That is, we solve an expected elliptic free boundary problem.
However, the available methods and approaches of elliptic FBPs do not directly apply to our problems,
such as
the variational methods
of Alt-Caffarelli \cite{AC} and Alt-Caffarelli-Friedman \cite{ACF,AltCafFried_Compres-83,AltCafFried_Compres},
the Harnack inequality approach of Caffarelli \cite{Ca1,Ca2,Ca3},
and other methods and approaches in many further works.
The main reason is that the type of equations needs to be first controlled in order to apply these methods,
which requires some strong estimates {\it a priori}.
To overcome the difficulties, we exploit the global structure of the problems,
which allows us to derive certain properties of
the solution (such as the monotonicity, {\it etc.}) so that the type of equations
and the geometry of the problem can be controlled.
With this, we solve the free boundary problem by the iteration procedure.

\vspace{3pt}
Notice that the existence of multiple wild solutions for the Cauchy problem of the
compressible Euler equations has been shown; see \cite{CDK,KKMM} and the references cited therein
for both the isentropic and full Euler cases.
In this paper, we focus on the solutions of specific structures motivated by underlying physics;
for these solutions, the uniqueness can be shown for all the cases as we discuss below.
Since we are interested in the solutions of specific structures,
we construct the solution in a carefully chosen class of solutions,
called admissible solutions.
This class of solutions needs to be defined with two somewhat opposite features:
the conditions
need not only to be flexible enough so that this class contains all possible solutions
of the problem which are of the desired structure, but also to be rigid enough
to force the desired structure of the solutions
with
the sufficient analytic and geometric control such that the expected estimates
for these solutions can be derived,
so that eventually a solution  can be constructed in this class by the iteration procedure.
In order to define such a class, we start with the solutions near some background solutions:
\begin{enumerate}
\item[\rm (i)]
to make sure that the solutions obtained are still in the same desired structure
via careful estimates,
which is the structure of transonic shock solutions in our application;

\item[(ii)]
to gain the insight and motivation for the structure and properties of the solutions that
are not near the background solution but have the required configuration
to form the conditions on which the {\it a priori} estimates and fixed point argument are based.
\end{enumerate}
In several problems,  we consider only the solutions near the background solution, as in \S2--\S 3 below.
In the other problems, say in \S \ref{self-simSect}--\S 5,
we carry out both steps described above
and construct admissible solutions which are not close to any known background solution.

\smallskip
Furthermore,
we emphasize that the elliptic and hyperbolic regions may
be separated not only by shocks, which are discontinuity fronts
of physical variables such as the velocity and the density,
but also by sonic curves/surfaces
where the type of equations changes continuously
in the physical variables,
as pointed out earlier.
This means that the ellipticity and
 hyperbolicity degenerate near the sonic curves/surfaces.
This presents additional difficulties in
 the analysis of such solutions.
Moreover, the sonic curves/surfaces may intersect
 the transonic shocks (see {\it e.g.}  Fig. \ref{figure: free boundary problems-1}, point $P_1$)
 so that, near such points, the analysis of solutions is even more involved.

\smallskip
The organization of this paper is as follows:
In \S 2, we start with our presentation of M-D transonic shocks
and free boundary problems for the compressible Euler equations for potential flow
in a setup as simple as possible,
and show how a transonic shock problem can be formulated as a free boundary problem for the
corresponding nonlinear PDEs of mixed elliptic-hyperbolic type.
Then we describe an efficient nonlinear method and related ideas and techniques,
first developed in Chen-Feldman \cite{CF-JAMS2003},
with focus on the key points in solving such free boundary problems
through this simplest setup.
In \S 3, we describe how they
can be applied
to establishing the existence, stability, and asymptotic behavior of
2-D steady transonic flows with transonic shocks past curved wedges
for the full Euler equations, by reformulating the problems as free boundary problems
via two different approaches.
In \S4, we describe how the nonlinear method and related ideas and techniques presented in \S 2--3
can be extended to the case of self-similar shock reflection/diffraction for the compressible Euler equations
for potential flow, including the von Neumann problem for shock reflection-diffraction
and the Prandtl-Meyer problem for unsteady supersonic flow onto solid wedges,
where the solutions have
the sonic arcs in addition to the transonic shocks.
In \S 5, we discuss some recent developments in the analysis of geometric properties of transonic shocks
as free boundaries in the 2-D self-similar coordinates for compressible fluid flows with focus on
the convexity properties of the self-similar transonic shocks obtained in \S 4.

\section{Multidimensional Transonic Shocks and Free Boundary
Problems for the Steady Euler Equations for Potential Flow}

\smallskip
For clarity, we start with our presentation of M-D transonic shocks
and free boundary problems for the compressible Euler equations in a setup as simple as possible,
and show how a transonic shock problem can be formulated as a free boundary problem for the
corresponding nonlinear PDEs of mixed elliptic-hyperbolic type.
Then we describe a
method first developed in Chen-Feldman \cite{CF-JAMS2003},
with focus on the key points to solve such free boundary problems
through this simplest setup.

The steady Euler equations for potential flow, consisting of the conservation law of mass and the
Bernoulli law for the velocity, can be written as the following
second-order nonlinear PDE of mixed elliptic-hyperbolic type
for the velocity
potential $\varphi: \bR^d\rightarrow\bR$ ({\it i.e.}, $\mathbf{u}=D\varphi$ is the velocity):
\begin{equation}\label{PotenEulerCompres}
\divg(\rho(|D\varphi|^2)D\varphi)=0,
\end{equation}
by scaling so that the density function $\rho(q^2)$ has the form:
\begin{equation}\label{densFunc}
\rho(q^2)=\big(1-\frac{\gamma-1}{2} q^2\big)^{\frac{1}{\gamma-1}},
\end{equation}
where
$\gamma>1$ is the adiabatic exponent
and $D:=(\partial_{x_1}, \dots, \partial_{x_d})$ is the gradient with respect to $\x=(x_1,\dots, x_d)\in \mathbb{R}^d$.

Equation \eqref{PotenEulerCompres} can be written in the non-divergence form:
\begin{equation}\label{PotenEulerCompres-b}
\sum_{i,j=1}^d \big(\rho(|D\varphi|^2) \delta_{ij}+2\rho'(|D\varphi|^2) \varphi_{x_i}\varphi_{x_j}\big)\varphi_{x_ix_j}=0,
\end{equation}
where the coefficients of the second-order nonlinear PDE \eqref{PotenEulerCompres-b} depend on $D\varphi$,
the gradient of the unknown function $\varphi$.

\medskip
The nonlinear PDE \eqref{PotenEulerCompres}, or equivalently \eqref{PotenEulerCompres-b} for smooth solutions,
is strictly elliptic at $D\varphi$ with $|D\varphi|=q$ if
\begin{equation}\label{elliptic}
\rho(q^2)+ 2q^2\rho'(q^2)>0,
\end{equation}
and is strictly hyperbolic if
\begin{equation}\label{hyperbolic}
\rho(q^2)+ 2q^2\rho'(q^2)<0.
\end{equation}
In fluid dynamics, the elliptic regions of equation \eqref{PotenEulerCompres}
correspond to the {\it subsonic flow}, the hyperbolic regions of
 \eqref{PotenEulerCompres} to the {\it supersonic flow}, and the regions
with $\rho(q^2)+ 2q^2\rho'(q^2)=0$ for $q=|D\varphi|$ to the {\it sonic flow}.

\smallskip
\subsection{Steady Transonic Shocks and Free Boundary Problems}
Let $\Omega\subset\bR^d$ be a domain ({\it i.e.}, simply connected open subset).
A function $\varphi\in W^{1,\infty}(\Omega)$ is a {\it weak solution}
of (\ref{PotenEulerCompres}) in $\Omega$ if

\begin{enumerate}\renewcommand{\theenumi}{\roman{enumi}}

\item $|D\varphi(\xx)|\leq \sqrt{2/(\gamma-1)}$ \, {\it a.e.}
   $\xx\in\Omega$,  that is, the physical region so that $\rho(|D\varphi(\xx)|^2)$
   is well defined via \eqref{densFunc} for {\it a.e.} $\xx\in\Omega$;

\item for any test function $\zeta\in C^\infty_0(\Omega)$,
\begin{equation}\label{PotenEulerCompresWeak}
\int_\Omega \rho(|D\varphi|^2)D\varphi\cdot D\zeta\,\dd \xx=0.
\end{equation}
\end{enumerate}

\smallskip
We are interested in the weak
solutions with shocks ({\it i.e.}, the surfaces of jump discontinuity
of $D\varphi$ of the solution $\varphi$ with codimension one) satisfying
the physical entropy condition that is consistent with the Second Law of Thermodynamics
in Continuum Physics.
More precisely, let $\Omega^+$ and $\Omega^-$ be open nonempty subsets of $\Omega$
such that
$$
\Omega^+\cap \Omega^-=\emptyset, \qquad
\overline{\Omega^+}\cup \overline{\Omega^-}=\overline\Omega,
$$
and $\sS:=\partial\Omega^+\setminus\partial\Omega$.
Let $\varphi\in W^{1,\infty}(\Omega)$ be a weak solution of
(\ref{PotenEulerCompres}) so that
$\varphi\in C^2({\Omega^\pm})\cap C^1(\overline{\Omega^\pm})$
and $D\varphi$ has a jump across $\sS$.

We now derive the necessary conditions on $\sS$ that is a $C^1$--surface
of codimension one.
First, the requirement that $\varphi$ is in $W^{1,\infty}(\Omega)$
yields ${\rm curl}(D\varphi)=0$ in the sense of distributions,
which implies
\begin{equation}\label{tangDerivEqual}
\varphi^+_{\boldsymbol{\tau}}=\varphi^-_{\boldsymbol{\tau}}\;\;\;\;\;\;\;\;
\mbox{on $\sS$},
\end{equation}
where
$$
\varphi^\pm_{{\boldsymbol\tau}}
:=D\varphi^\pm -(D\varphi^\pm\cdot{\boldsymbol{\nu}}){\boldsymbol{\nu}}
$$
are the trace values of the tangential gradients of $\varphi$ on $\sS$
in the tangential space with $(d-1)$-dimension
on the $\Omega^\pm$ sides, respectively, and
$\nnu$ is the unit normal to $\sS$ from $\Omega^-$ to $\Omega^+$.
Then we simply write $\varphi_{\boldsymbol\tau}:=\varphi^\pm_{\boldsymbol \tau}$ on $\sS$
and choose
\begin{equation}\label{continuity}
\varphi^+=\varphi^-\; \qquad \mbox{on $\sS$}
\end{equation}
to be consistent with the $W^{1,\infty}$--requirement of $\varphi$.

\smallskip
Now, for $\zeta\in C^\infty_0(\Omega)$,
we use (\ref{PotenEulerCompresWeak})
to compute
\begin{align*}
0&=\left(\int_{\Omega^+}+\int_{\Omega^-}\right)\rho(|D\varphi|^2)
    D\varphi\cdot D\zeta\,\dd \x \\
&=-\int_{\partial\Omega^+}\rho(|D\varphi|^2)
    D\varphi\cdot \nnu\, \zeta\,\dd{\mathcal H}^{d-1}
    +\int_{\partial\Omega^-}\rho(|D\varphi|^2)
    D\varphi\cdot \nnu \, \zeta\,\dd{\mathcal H}^{d-1}\\
&=\int_\sS\left(-\rho(|D\varphi^+|^2)
     D\varphi^+\cdot \nnu +\rho(|D\varphi^-|^2)
   D\varphi^-\cdot \nnu\right)\zeta\,\dd{\mathcal H}^{d-1},
\end{align*}
where ${\mathcal H}^{d-1}$ is the $(d-1)$-D Hausdorff measure,
{\it i.e.}, the surface area measure.
Thus, the other condition on $\sS$, which measures the trace jump
of the normal derivative of $\varphi$ across $\sS$, is
\begin{equation}\label{FBCondition-1}
\rho(|D\varphi^+|^2) \varphi_\nnu^+
 =\rho(|D\varphi^-|^2)
  \varphi_\nnu^-\;\;\;\;\;\;\;\; \mbox{on $\sS$},
\end{equation}
where
$\varphi_{\bn}^\pm=D\varphi^\pm\cdot \bn$ are the trace values of the
normal derivative of $\varphi$ along $\sS$ on the $\Omega^\pm$ sides, and
$$
\rho(|D\varphi^\pm|^2)
=\Big(1-\frac{\gamma-1}{2}\big(|\varphi^\pm_\bt|^2
+|\varphi^\pm_\bn|^2\big)\Big)^{\frac{1}{\gamma-1}},
$$
respectively.

Conditions \eqref{continuity}--\eqref{FBCondition-1} are called the Rankine-Hugoniot conditions for potential flow
in fluid dynamics.
On the other hand, it can also be shown that any $\varphi\in C^2({\Omega^\pm})\cap C^1(\overline{\Omega^\pm})$ that is
a $C^2$--solution of (\ref{PotenEulerCompres}) in $\Omega^\pm$ respectively, such that $D\varphi$ has a jump across $\sS$ satisfying the Rankine-Hugoniot conditions  \eqref{continuity}--\eqref{FBCondition-1},
must be a weak solution of
(\ref{PotenEulerCompres}) in the whole domain $\Omega$.
Therefore, the necessary and sufficient conditions for
$\varphi\in C^2({\Omega^\pm})\cap C^1(\overline{\Omega^\pm})$ that is a solution of (\ref{PotenEulerCompres}) in $\Omega^\pm$ respectively
 to be a weak solution of
(\ref{PotenEulerCompres}) in the whole domain $\Omega$
are the Rankine-Hugoniot conditions \eqref{continuity}--\eqref{FBCondition-1}.

\bigskip
For given $K>0$, consider the function:
\begin{equation}\label{Phi-K}
\displaystyle\Phi_K(p)\defd\big(K-\frac{\gamma-1}{2} p^2\big)^{\frac{1}{\gamma-1}}p
\qquad\,\,\mbox{ for $\displaystyle p\in [0, \sqrt{2K/(\gamma-1)}]$}.
\end{equation}
Then $\Phi_K\in C([0, \sqrt{2K/(\gamma-1)}])$ and
\begin{align}
&\mbox{$\Phi_K(p)> 0$ for $p\in (0, \sqrt{2K/(\gamma-1)})$, $\quad\,\,$
  $\Phi_K(0)=\Phi_K(\sqrt{2K/(\gamma-1)})=0$}, \label{FBCOndLemma_2}\\[1mm]
&\mbox{$0<\Phi_K'(p)\le K^{\frac{1}{\gamma-1}}\,\, $ for $p\in (0, p^K_{\rm sonic})$},
     \label{FBCOndLemma_3}\\[1mm]
&\mbox{$\Phi_K'(p) <0\,\,$  for $p\in (p^K_{\rm sonic},\sqrt{2K/(\gamma-1)})$}, \label{FBCOndLemma_3b}\\[1mm]
&\mbox{$\Phi_K^{''}(p)<0\,$ for $p\in (0, p^K_{\rm sonic}]$,}\label{FBCOndLemma_4}
\end{align}
where
\begin{equation}\label{Pcritical}
p^K_{\rm sonic}:=\sqrt{2K/(\gamma+1)}.
\end{equation}

By direct calculation, condition \eqref{elliptic} is equivalent to
$\Phi_1'(q)>0$, and condition  \eqref{hyperbolic} is equivalent to $\Phi_1'(q)<0$. Thus, using \eqref{FBCOndLemma_3},
we obtain that PDE \eqref{PotenEulerCompres}
is strictly elliptic at $D\varphi$ if
$|D\varphi|<p^1_{\rm sonic}$ and is strictly hyperbolic if $|D\varphi|>p^1_{\rm sonic}$,
where we have used notation \eqref{Pcritical}.

\medskip
Suppose that $\varphi(x)$ is a solution satisfying
\begin{equation}\label{TnasonicCond}
|D\varphi|<p^1_{\rm sonic}=\sqrt{2/(\gamma+1)}
\;\; \mbox{in $\Omega^+$,}\;
\qquad \;\; |D\varphi|>p^1_{\rm sonic} \;\;\mbox{in $\Omega^-$,}
\end{equation}
and
\begin{equation} \label{2.9a}
 D\varphi^\pm\cdot\bn>0 \qquad \text{on $\sS$},
\end{equation}
besides \eqref{continuity} and \eqref{FBCondition-1}.
Then $\varphi(x)$ is a {\it transonic shock solution}
with {\it transonic shock} $\sS$ that divides
the {\em subsonic region} $\Omega^+$
from the {\em supersonic region} $\Omega^-$.
In addition, $\varphi(\x)$ satisfies the {\it physical entropy
condition} (see Courant-Friedrichs \cite{CF}; also see \cite{Da,Lax}):
\begin{equation} \label{2.8a}
\rho(|D\varphi^-|^2)<\rho(|D\varphi^+|^2),
\end{equation}
which implies, by (\ref{2.9a}), that {\it the density $\rho$ increases in the
flow direction}; that is, the transonic shock solution is an entropy solution.
Note that equation (\ref{PotenEulerCompres})
is elliptic in the subsonic region $\Omega^+$ and
hyperbolic in the supersonic region $\Omega^-$.

For clarity of presentation of the nonlinear method,
first developed in Chen-Feldman \cite{CF-JAMS2003},
we focus first on the free boundary problem in the simplest setup,
while the method and related ideas and techniques
have been applied to more general free boundary problems
involving transonic shocks for the nozzle problems and other important problems,
some of which will be discussed in \S 3--\S5.

\medskip
Let $(\xpr,x_d)$ be the coordinates of $\bR^d$ with
$\xpr=(x_1,\cdots,x_{d-1})\in\bRnpr$ and $x_d\in \bR$.
{}From now on, in this section, we focus on $\Omega:=\Sqr\times(-1, 1)$ for simplicity,
without loss of our main objectives.

Let $\displaystyle q^-\in (p^1_{\rm sonic}, \sqrt{2/(\gamma-1)})$
and $\varphi_0^-(x)\defd q^-x_d$. Then $\varphi_0^-$
is a supersonic solution in $\Omega$.
From \eqref{FBCOndLemma_2}--\eqref{FBCOndLemma_4},
there exists a unique $q^+\in (0, p^1_{\rm sonic})$ such that
\begin{equation}\label{pPlusMinCondit}
\big(1-\frac{\gamma-1}{2}(q^+)^2\big)^{\frac{1}{\gamma-1}}q^+=
\big(1-\frac{\gamma-1}{2}(q^-)^2\big)^{\frac{1}{\gamma-1}}q^-.
\end{equation}
In particular, $q^+<q^-$. Define $\varphi_0^+(\x):=q^+x_d$ in $\Omega$.
Then the function:
\begin{equation} \label{linearTransonicSol}
\varphi_0(\x)=\min(\varphi^+_0(\x),\varphi_0^-(\x))
\end{equation}
is a transonic shock solution in $\Omega$, in which
$\Omega^-_0=\{x_d\le 0\}\cap\Omega$
and $\Omega^+_0=\{x_d\ge 0\}\cap\Omega$
are the supersonic and subsonic regions of $\varphi_0(\x)$, respectively.
Also, the boundary condition: $(\varphi_0)_{\bn}=0$ holds on
$\partial \Sqr \times [-1, 1]$.

\medskip
We start with perturbations of the background solution $\varphi_0(\x)$ defined in (\ref{linearTransonicSol}).
We use the following H\"{o}lder norms:
For $\alpha\in (0, 1)$ and any non-negative integer $k$,
\begin{eqnarray}
&&[u]_{k,0,\Omega}
=\sum_{|\bb|=k}\sup_{x\in\Omega}\,|D^\bb u(\x)|,
\qquad [u]_{k,\alpha,\Omega}
=\sum_{|\bb|=k}\sup_{\x,\y\in\Omega, \x\ne \y}
  {\frac{|D^\bb u(\x)-D^\bb u(\y)|}{|\x-\y|^\alpha}},
\label{holdNrms_usual}\\
&&\|u\|_{k,0,\Omega}
=\sum_{j=0}^k [u]_{j,0,\Omega},
\qquad \|u\|_{k,\alpha,\Omega}=\|u\|_{k,0,\Omega} +[u]_{k,\alpha,\Omega},
\nonumber
\end{eqnarray}
where $\bb=(\beta_1, \cdots,\beta_d), \,\beta_l\ge 0$ integers,
$D^\bb=\partial_{x_1}^{\beta_1}\cdots\partial_{x_d}^{\beta_d}$,
and $|\bb|=\beta_1+\cdots + \beta_d$.

\medskip
Then the transonic shock problem can be formulated as the following problem:

\medskip
\begin{problem}\label{problem-a}
Given a supersonic solution $\varphi^-$ of \eqref{PotenEulerCompres}
in $\Omega$, which is a $C^{2,\alpha}$--perturbation of $\varphi_0^-${\rm :}
\begin{equation}\label{smallPert}
\|\varphi^- - \varphi_0^-\|_{2,\alpha, \Omega}\le \epsP
\end{equation}
for some $\alpha\in (0,1)$ with small $\sigma>0$ and satisfies
\begin{equation}\label{smallPertBdryCond}
\varphi^-_\bn=0\; \;
\qquad  \mbox{on $\,\,\partial\Sqr \times [-1, 1]$,}
\end{equation}
find a transonic shock solution $\varphi$ in $\Omega$
such that
$$
\varphi=\varphi^-\qquad \mbox{ in $\;\Omega^-:=\Omega\setminus \overline{\Omega^+}$},
$$
where
$\Omega^+:=\{\x\in\Omega\;:\;\;|D\varphi(\x)|<p^1_{\rm sonic}\}$
is the subsonic region of $\varphi$, which is the complementary
set of the supersonic region of $\varphi$ in $\Omega$,
and
\begin{equation} \label{bdryConditions}
\begin{cases} \displaystyle
(\varphi, \varphi_{x_n})=(\varphi^-, \varphi^-_{x_n})\; \qquad &\mbox{on $\,\,\Sqr \times \{-1\}$}, \\
\varphi=\varphi_0^+\; \qquad &\mbox{on $\,\,\Sqr \times \{1\}$}, \\
\varphi_\bn=0\; \qquad &
\mbox{on $\,\,\partial \Sqr \times [-1, 1]$}.
\end{cases}
\end{equation}
\end{problem}

\medskip
Since $\varphi=\varphi^-$ in $\Omega^-$,
$|D\varphi| < p^1_{\rm sonic} < |D\varphi^-|$ in $\Omega^+$,
$|D\varphi^-|\thicksim\partial_{x_d}\varphi^- > p^1_{\rm sonic}$ in $\Omega$,
and it is expected that  $\Omega^+=\{x_d>f(\x')\}\cap \Omega$ and
$|D\varphi|\thicksim\partial_{x_d}\varphi< p^1_{\rm sonic} $ in $\Omega^+$ with
\eqref{continuity} across the transonic shock $\sS=\{x_d=f(\x')\}\cap \Omega$,
then $\varphi$ should satisfy
\begin{equation} \label{transonicInequalities-1}
\varphi(\x)\leq \varphi^-(\x)\;\;\; \qquad \hbox{for $\, \x\in\Omega$.}
\end{equation}
This motivates the following reformulation of Problem \ref{problem-a}
as a free boundary problem for the subsonic (elliptic)
part of the solution:

\medskip
\begin{problem}[Free Boundary Problem]\label{problem-b}
Find $\varphi\in C(\overline\Omega)$ such that
\begin{enumerate}\renewcommand{\theenumi}{\roman{enumi}}
\item[\rm (i)]\label{transonicFBP-item1}
  $\varphi$ satisfies \eqref{transonicInequalities-1} in $\Omega$
  and \eqref{bdryConditions}
  on $\partial \Omega${\rm ;}
\item[\rm (ii)] \label{transonicFBP_item3}
 $\varphi\in C^{2,\alpha}(\overline{\Omega^+})$ is a solution
 of \eqref{PotenEulerCompres} in
 $\Omega^+= \{\x\in\Omega\;:\;\varphi(\x)<\varphi^-(\x)\}$,
  the non-coincidence set{\rm ;}
\item[\rm (iii)] \label{transonicFBP-item6}
 the free boundary $\sS=\partial\Omega^+\cap\Omega$ is
 given by $x_d=f(\xpr)$ for $\xpr \in \Sqr$
 so that $\Omega^+=\{x_d>f(\x')\;:\; x'\in \Sqr\}$ with
 $f\in C^{2,\alpha}([0,a]^{d-1})${\rm ;}
\item[\rm (iv)] \label{transonicFBP-item7}
 the free boundary condition
 \eqref{FBCondition-1} holds on $\sS$.
\end{enumerate}
\end{problem}

\medskip
\noindent
In the free boundary problem (Problem \ref{problem-b}) above,
phase $\varphi^-$ is not required to be a solution of
\eqref{PotenEulerCompres} and $\varphi$ is not necessary to be subsonic
in $\Omega^+$, although we require the subsonicity in Problem \ref{problem-a}
so that the free boundary is a transonic shock.

It is proved in Chen-Feldman \cite{CF-JAMS2003} that,
if perturbation $\varphi^--\varphi^-_0$ is
small enough in $C^{2,\alpha}$, then the free boundary
problem (Problem \ref{problem-b}) has a solution that is subsonic on $\Omega^+$,
so that Problem \ref{problem-a}
has a transonic shock solution.
Furthermore, the transonic shock is stable under any small
$C^{2,\alpha}$--perturbation of $\varphi^-$.

\begin{theorem}[Chen-Feldman \cite{CF-JAMS2003}]\label{mainFBPTh}
Let $\displaystyle q^+\in (0, p^1_{\rm sonic})$ and
$\displaystyle q^- \in (p^1_{\rm sonic}, \sqrt{2/(\gamma-1)} )$
satisfy {\rm (\ref{pPlusMinCondit})}.
Then there exist positive constants $\epsP_0$, $C_1$, and $C_2$ depending
only on  $(q^+, d, \gamma)$ and $\Omega$ such that,
for every $\epsP\leq\epsP_0$ and any function $\varphi^-$ satisfying
{\rm (\ref{smallPert})}--{\rm (\ref{smallPertBdryCond})},
there exists a unique solution $\varphi$ of the free boundary problem,
{\rm Problem \ref{problem-b}}, satisfying
$$
\|\varphi - \varphi_0^+\|_{2,\alpha, \Omega^+}\le C_1\epsP
$$
and $|D\varphi|< p^1_{\rm sonic}$ in $\Omega^+$.
Moreover, $\Omega^+=\{x_d> f(\xpr)\}\cap\Omega$ with
$f: \bR^{d-1} \rightarrow \bR$ satisfying
\begin{eqnarray*}
\|f\|_{2, \alpha, \bRnpr}\leq C_2\epsP,\,\qquad\,\,
D_{\x'}f(\x')=0\;\;\mbox{on $\partial\Sqr$},
\end{eqnarray*}
that is, the free boundary
$\sS=\{(\x',x_d)\;:\,\; x_d=f(\x'), \x'\in \bRnpr\}\cap\Omega $
is in $C^{2, \alpha}$ and orthogonal to
$\partial \Omega$ at their intersection points.
\end{theorem}

In particular, we obtain
\begin{corollary}\label{mainFBPTh-cor}
Let $q^\pm$ be as in {\rm Theorem \ref{mainFBPTh}}, and let $\sigma_0$ be the constant
defined in {\rm Theorem  \ref{mainFBPTh}}.
If $\varphi^-(\x)$ is a supersonic solution  of \eqref{PotenEulerCompres}
satisfying {\rm (\ref{smallPert})}--{\rm (\ref{smallPertBdryCond})} with $\sigma\le\sigma_0$,
then there exists  a transonic shock solution $\varphi$ of {\rm Problem \ref{problem-a}}
with shock $\sS=\{(\x',x_d)\;:\,\; x_d=f(\x'), \x'\in \bRnpr\}\cap\Omega$
such that $\varphi$ and $f$ satisfy the properties stated in {\rm Theorem \ref{mainFBPTh}}.
\end{corollary}

Indeed, under the conditions of Corollary \ref{mainFBPTh-cor},
 solution $\varphi$ of {\rm Problem \ref{problem-b}} obtained in
 Theorem \ref{mainFBPTh}, along with the free boundary
$\sS=\{(\x',x_d)\;:\,x_d=f(\x'), \x'\in \bRnpr\}\cap\Omega$,
forms a transonic shock solution of {\rm Problem \ref{problem-a}}.

\medskip
The following features of equation (\ref{PotenEulerCompres})
and the free boundary condition (\ref{FBCondition-1}) are employed
in the proof of Theorem \ref{mainFBPTh}.
\begin{enumerate}
\item[(i)] The nonlinear equation (\ref{PotenEulerCompres}) is uniformly
elliptic only if $|D\varphi|<p^1_{\rm sonic}-\varepsilon$ in $\Omega^+$ for
some $\varepsilon>0$;
\item[(ii)] $|D\varphi^+|=\left(|\varphi^+_\bn|^2+|\varphi_\bt|^2\right)^{1/2}$
on $\sS$ is subsonic only if $\varphi_\bt$ is sufficiently small;

\item[(iii)] The free boundary condition (\ref{FBCondition-1})
is uniformly non-degenerate ({\it i.e.}, $\varphi^-_\bn-\varphi^+_\bn$
is bounded from below by a positive constant on $\sS$) only if
$\varphi^-_\bn > p^K_{\rm sonic}+\varepsilon$ on $\sS$ for some $\epsilon>0$
with $K=1-\frac{\gamma-1}{2}|\varphi_\bt|^2$.
\end{enumerate}
By (\ref{smallPert}), these conditions hold if, for any $\x\in \sS$,
the unit normal $\bn(\x)$  to $\sS$
is sufficiently close to being orthogonal to $\{x_d=0\}$.

\subsection{A Nonlinear Method for Solving the Free Boundary Problems for Nonlinear PDEs of Mixed Elliptic-Hyperbolic Type}
\label{existenceSection}

We now describe a nonlinear method
and related ideas and techniques, developed first in Chen-Feldman \cite{CF-JAMS2003},
for the construction of solutions of the free boundary problems for nonlinear PDEs of mixed elliptic-hyperbolic type,
through Problem \ref{problem-b} as the simplest setup.
We present the version of the method that is restricted to this setup.
The key ingredient is an iteration scheme, based on the non-degeneracy of the free boundary  condition: the jump of
the normal derivative of solutions across the free boundary has a strict lower bound.
Since the PDE is of mixed type, we make a cutoff (truncation) of the nonlinearity near the value related to
the background solution in order to fix the type of equation (to make it elliptic everywhere) and,
at the fixed point of the iteration, we remove the cutoff eventually by a required estimate.
The iteration set consists of the functions close to the background solution -- in the $C^{2,\alpha}$--norm in the present case.
Then, for each function from the iteration set, the nondegeneracy allows of using one of the Rankine-Hugoniot conditions,
equality \eqref{continuity}, to define the iteration free boundary, which is a smooth graph.
In domain $\Omega^+$ determined by the iteration free boundary, we solve a boundary value problem with the truncated PDE,
the condition on the shock derived from the other Rankine-Hugoniot condition \eqref {FBCondition-1}
by a truncation (similar to the truncation of the PDE)
and other appropriate modifications to achieve the uniform obliqueness,
and the same boundary conditions as in the original problem for the iteration problem
on the other parts of the boundary of the iteration domain.
The solution of this iteration problem defines
the iteration map.
We exploit the estimates for the iteration
problem to prove the existence of a fixed point of the iteration map, and then we show that a fixed point is
a solution of the original problem.

In some further problems, we look for the solutions that are not close to a known background solution.
Some of these problems, as well as the corresponding versions of the nonlinear method
described above, are discussed in \S \ref{self-simSect}. A related method for the construction of perturbations of transonic shocks for the steady transonic small disturbance model was proposed in \cite{CanicKeyfitz}, in which the type of equation depends on the solution only (but not on its gradient) so that the ellipticity can be controlled by the maximum principle; also see \cite{CKK1b}.

\medskip
\noindent
{\bf 2.2.1. Subsonic Truncations -- Shiffmanization}.
In order to solve the free boundary problem,
we first reformulate Problem \ref{problem-b} as a truncated one-phase free boundary problem,
motivated by the argument introduced originally
in Shiffman \cite{Sh}, now so called the {\it shiffmanization} ({\it cf}. Lax \cite{Lax2003}); also see
\cite[pp. 87--90]{AltCafFried_Compres}.
This is achieved by modifying both the nonlinear equation (\ref{PotenEulerCompres})
and the free boundary condition (\ref{FBCondition-1})
to make the equation uniformly elliptic
and the free boundary
condition non-degenerate. Then we solve the truncated one-phase free boundary problem with the modified equation in the downstream
region, the modified free boundary condition, and the given hyperbolic phase in the upstream region.
By a careful gradient estimate later on, we prove that the solution
in fact solves the original problem.
We note that, for the steady potential flow equation \eqref{PotenEulerCompres}, the coefficients
of its non-divergent form \eqref{PotenEulerCompres-b} depend on $D\varphi$,
so the type of equation depends on $D\varphi$.

We first recall that the ellipticity condition for (\ref{PotenEulerCompres})
at $|D\varphi|=q$ is \eqref{elliptic}, which is equivalent to
\begin{equation}\label{ellipticityCondPhi}
\Phi_1'(q)>0,
\end{equation}
where $\Phi_K(p)$ is the function defined in \eqref{Phi-K}.
By \eqref{FBCOndLemma_3},
inequality (\ref{ellipticityCondPhi}) holds for $q\in (0,p^1_{\rm sonic})$.

The truncation
is done by
modifying $\Phi_1(q)$ so that the new function
$\tilde\Phi_1(q)$ satisfies (\ref{ellipticityCondPhi}) uniformly for
all $q>0$ and, around $q^+$, $\tilde\Phi_1(q)=\Phi_1(q)$.
More precisely, the procedure consists of the following steps:

\medskip
{1.} Denote
$\varepsilon:= \frac{p^1_{\rm sonic} - q^+}{2}$.
Let $y=c_0q+c_1$ be the tangent line of the graph
of $y=\Phi_1(q)$ at $q=p^1_{\rm sonic}-\varepsilon$. Then,
using \eqref{FBCOndLemma_3}, we obtain
$
c_0=\Phi_1'(p^1_{\rm sonic}-\varepsilon)>0.
$
Define $\tilde\Phi_1: [0, \infty)\rightarrow \bR$ as
\begin{equation}\label{defTildePhi}
\tilde\Phi_1(q)=
\begin{cases}
\Phi_1(q)\; & \mbox{if $0\leq  q< p^1_{\rm sonic}-\varepsilon$},\\
c_0q+c_1\;  & \mbox{if $q>p^1_{\rm sonic}-\varepsilon$},
\end{cases}
\end{equation}
which satisfies $\tilde\Phi_1\in C^{1,1}([0, \infty))$.

\medskip
{2.} Define
\begin{equation}\label{defTrho}
\displaystyle \Trho(s)=\frac{\tilde\Phi_1(\sqrt{s})}{\sqrt{s}}\; \qquad\;
\mbox{for $s\in[0, \infty)$}.
\end{equation}
Then $\Trho\in C^{1,1}([0, \infty))$ and
\begin{equation}\label{goodTruncatEquat}
\Trho(q^2)=\rho(q^2)\;\;\qquad \mbox{if}\;\;
                  0\leq q<p^1_{\rm sonic}-\varepsilon.
\end{equation}
By \eqref{FBCOndLemma_3}--\eqref{FBCOndLemma_4}
and the definition of $\tilde\Phi_1$ in \eqref{defTildePhi},
$$
0< c_0=\Phi_1'(p^1_{\rm sonic}-\varepsilon) \leq \tilde\Phi_1'(q)=\Trho(q^2)+2q^2\Trho'(q^2)
\leq C \qquad\;\hbox{for $\, q\in(0, \infty)$}
$$
for some constant $C>0$.
Then the equation:
\begin{equation}\label{PotenEulerCompresTrunc}
\Le\varphi\defd \divg(\Trho(|D\varphi|^2)D\varphi)=0
\end{equation}
is uniformly elliptic, with ellipticity constants
depending only on $q^+$ and $\gamma$.

\medskip
{3.} We also do the corresponding truncation of the
free boundary condition \eqref{FBCondition-1}:
\begin{equation}\label{FBConditionConormTr}
\Trho(|D\varphi|^2) \varphi_\bn
 =\rho(|D\varphi^-|^2)
  D\varphi^-\cdot\bn\;\;\;\;\;\;\;\; \mbox{on $\sS$}.
\end{equation}
On the right-hand side of (\ref{FBConditionConormTr}),
we use the non-truncated function $\rho$ since
$\rho\ne\Trho$ on the range of $|D\varphi^-|^2$.
Note that
(\ref{FBConditionConormTr}), with the right-hand side considered as
a known function, is the conormal boundary
condition for the uniformly elliptic equation
(\ref{PotenEulerCompresTrunc}).

\medskip
{4.} Introduce the function:
$$
u:= \varphi^--\varphi.
$$
Then,
by \eqref{transonicInequalities-1},
the problem is to find
$u \in C(\overline\Omega)$ with $u\ge 0$ such that
\begin{enumerate}\renewcommand{\theenumi}{\roman{enumi}}
\item[\rm (i)]\label{C_transonicFBP_item3}
 $u\in C^{2,\alpha}(\overline{\Omega^+})$ is a solution
of
\begin{align}
&\divg A(Du,\x)=F(\x)  \quad && \mbox{ in $\Omega^+\defd \{u>0\}\cap\Omega$ (the non-coincidence set)},
                \label{generalFBP:1}  \\
&A(Du,\x)\cdot \bn = G(\bn,\x)\quad && \mbox{ on $\sS \defd \partial\Omega^+\setminus \partial \Omega$},
\label{generalFBP:2}
\end{align}
and the boundary condition on $\partial\Omega$  determined by \eqref{bdryConditions} and $\varphi^-(\x)$:
\begin{equation} \label{bdryConditions-b}
\begin{cases} \displaystyle
u=0 \qquad &\mbox{on $\,\Sqr \times \{-1\}$}, \\
u=\varphi^--\varphi_0^+\; \qquad &\mbox{on $\,\Sqr \times \{1\}$}, \\
u_\bn=0 \qquad &
\mbox{on $\,\partial \Sqr \times [-1, 1]$},
\end{cases}
\end{equation}
where $\bn$ is the unit normal to $\sS$ towards the unknown phase
and
\begin{eqnarray*}
&&A(P,\x)=\Trho(|D\varphi^-(\x)-P|^2) (D\varphi^-(\x)-P)-
  \Trho(|D\varphi^-(\x)|^2) D\varphi^-(\x) \qquad
  \mbox{for $ P \in \bR^d$}, \\
&&F(\x) =-{\rm \divg}(\Trho(|D\varphi^-(\x)|^2)D \varphi^-(\x)), \\
&&G(\bn,\x)=\big(\rho(|D\varphi^-(\x)|^2)- \Trho(|D\varphi^-(\x)|^2)\big)
   D\varphi^-(\x)\cdot\bn.
\end{eqnarray*}
Note that condition \eqref{smallPertBdryCond} has been used to determine the third condition
in
\eqref{bdryConditions-b}.

\smallskip
\item[\rm (ii)] \label{C_transonicFBP-item6}
 the free boundary $\sS:=\partial\Omega^+\cap\Omega=\{x_d=f(\xpr)\,:\,\xpr \in \Sqr\}$
 so that $\Omega^+=\{x_d>f(\x')\}\cap \Omega$ with $f\in C^{2,\alpha}([0,a]^{d-1})$ and
 $D_{\x'}f=0$ on $\partial(\Sqr\times [-1,1])$.
\end{enumerate}

\bigskip
\noindent
{\bf 2.2.2. Domain Extension}.
We then extend domain $\Omega$ of the truncated free boundary problem in \S 2.2.1 above
to domain $\Omega_{\rm e}$,
so that the whole free boundary lies in the interior of the extended domain.
This is possible owing to the simple geometry of the domain, as considered in this section.

\smallskip
Notice that, for a function $\phi\in C^{2,\alpha}(\overline{\Omega})$
with $\Omega:=\Sqr\times(-1, 1)$ satisfying
\begin{equation}\label{normDerivCond}
\phi_\bn=0\; \qquad \mbox{on} \,\,\, \partial\Sqr\times[-1, 1],
\end{equation}
we can extend $\phi$ to $\bRnpr\times [-1, 1]$
so that the extension (still denoted) $\phi$ satisfies
$$
\phi\in C^{2,\alpha}(\bRnpr \times [-1, 1]),
$$
and, for every $m=1,\cdots, n-1$, and $k=0, \pm1, \pm2, \cdots$,
\begin{equation}\label{symmetricFunc}
\phi(x_1, \cdots, x_{m-1}, k-z, x_{m+1},\cdots, x_d)
=\phi(x_1, \cdots, x_{m-1}, k+z, x_{m+1},\cdots, x_d),
\end{equation}
that is, $\phi$ is symmetric with respect
to every hyperplane $\{x_m=k\}$.
Indeed, for ${\bf k}=(k_1, \cdots, k_{d-1},0)$
with integers $k_j, j=1,\cdots d-1$,  we define
$$
\phi(\x+{\bf k})= \phi(\eta({x_1,k_1}),\cdots,\eta(x_{d-1},k_{d-1}),x_d)\;
\;\;\;\;\;\;\mbox{for $\; \x \in\Sqr\times [-1,1]$}
$$
with
$$
\eta(t, k)=
\begin{cases}
t\;\;\;\;&\mbox{if $k$ is even},\\
1-t\; &\mbox{if $k$ is odd}.
\end{cases}
$$
It follows from (\ref{symmetricFunc}) that $\phi(\x',x_d)$
is $2$-periodic in each variable of $(x_1, \cdots, x_{d-1})$:
$$
\phi(\x+2 \ee_m)=\phi(\x) \qquad \hbox{for}\,\,\,
  \x\in\bRnpr\times [-1, 1], \,\,m=1, \cdots, d-1,
$$
where $\ee_m$ is the unit vector in the direction of $x_m$.

Thus, with respect to the $2$-periodicity,
we can consider $\phi$ as a function on
$\RDom\defd\bTpr\times [-1, 1]$, where $\bTpr$
is a  flat torus in $d-1$ dimensions with its coordinates
given by cube $(0, 2)^{d-1}$.
Note that (\ref{symmetricFunc}) represents an
extra symmetry condition, in addition to
$\phi\in C^{2,\alpha}(\bTpr\times[-1, 1])$,
and (\ref{symmetricFunc}) implies (\ref{normDerivCond}).

\medskip
Then we can extend $\varphi^-$ in the same way by \eqref{smallPertBdryCond},
that is,
$\varphi^-\in C^{2,\alpha}(\RDom)$ satisfies (\ref{symmetricFunc}).
Notice that $\varphi^\pm_0$ can also be considered as the functions
in $\RDom$ satisfying (\ref{symmetricFunc}), since
$\varphi^\pm_0(\x)=q^\pm x_d$ in $\bRnpr\times [-1, 1]$
which are independent of $\xpr$.

\medskip
Therefore, we have reduced the transonic shock problem, Problem \ref{problem-b},
into the following free boundary problem:

\begin{problem}\label{problem-c}
Find $u\in C(\overline\RDom)$ with $u\ge 0$ such that
\begin{enumerate}[{\rm (i)}]
\item
 $u\in C^{2,\alpha}(\overline{\Omega_{\rm e}^+})$ is a solution
of \eqref{generalFBP:1}
in $\Omega_{\rm e}^+:=\{u(\x)>0\}\cap \RDom$,
 the non-coincidence set{\rm ;}
\item the first two conditions in \eqref{bdryConditions-b} hold on $\partial \RDom$,
i.e., $u=0$ on $\partial \RDom\cap\{x_n=-1\}$ and
\begin{equation}\label{bdryCondNozleExitU}
 u=\varphi^--\varphi_0^+\qquad\mbox{on $\partial \RDom\cap\{x_n=1\}$};
\end{equation}
\item
 the free boundary $\sS=\partial\Omega^+\cap\RDom$ is
 given by  $x_d=f(\xpr)$ for $\xpr \in \bTpr$
 so that $\Omega^+=\{x_d>f(\x')\;:\;\xpr\in \bTpr\}$ with $f\in C^{2,\alpha}(\bTpr)${\rm ;}

\smallskip
\item
 the free boundary condition \eqref{generalFBP:2}
 holds on $\sS$.
\end{enumerate}
\end{problem}

\medskip
As indicated in \S 1,
one of the main difficulties for solving the modified
free boundary problem, Problem \ref{problem-c},
is that the methods presented in the previous works for elliptic free boundary problems do not
directly apply.
Indeed, equation \eqref{generalFBP:1} is quasilinear, uniformly elliptic,
but does not have a clear variational structure, while
 $G(\bn,\x)$  in the free boundary condition
\eqref{generalFBP:2} depends on $\bn$. Because of these features,
 the variational methods
in \cite{AC,ACF} do not directly apply to Problem \ref{problem-c}.
Moreover, the nonlinearity in our problem makes it difficult
to apply the Harnack inequality approach of Caffarelli in \cite{Ca1,Ca2,Ca3}.
In particular, a boundary comparison principle for positive solutions
of elliptic equations in Lipschitz domains is unavailable in our
case that the nonlinear PDEs are not homogeneous
with respect to $(D^2u, Du, u)$ here.
Therefore, a different method is required to overcome
these difficulties
for solving Problem \ref{problem-c}.

\bigskip
\noindent
{\bf 2.2.3. Iteration Scheme for Solving Free Boundary Problems.}
The iteration scheme, developed in Chen-Feldman \cite{CF-JAMS2003},
is based on the non-degeneracy of the free boundary condition:
the jump of the normal derivative of a solution across the free boundary
has a strictly positive lower bound.

Denote $u_0:=\varphi^--\varphi_0^+$. Note that $u_0$ satisfies the nondegeneracy condition:
$\partial_{x_d}u_0=q^--q^+>0$ in $\RDom$. Assume that \eqref{smallPert}
holds with $\sigma\le \frac{q^--q^+}{10}$.
 Let a function $w$ on $\RDom$ be given such that
$\|w-(\varphi^--\varphi_0^+)\|_{C^{2,\alpha}(\overline{\RDom})} \le \frac{q^--q^+}{10}$,
which implies that
$w$ satisfies the  nondegeneracy condition:
$\partial_{x_d}w\ge \frac{q^--q^+}2>0$ in $\RDom$.
Define domain $\Omega^+(w):=\{w>0\}\subset \RDom$. Then
$$
\Omega^+(w)=\{x_d>f(\xpr)\;:\; \xpr\in \bTpr\}, \qquad
\sS(w)\defd\partial\Omega^+(w)\setminus
\partial\RDom=\{x_d=f(\xpr)\;:\; \xpr\in \bTpr\}
$$
with
$f\in C^{2, \alpha}(\bTpr)$.
We solve the oblique derivative problem  (\ref{generalFBP:1})--\eqref{generalFBP:2} and \eqref{bdryCondNozleExitU}
in $\Omega^+(w)$
to obtain a solution $u\in C^{2, \alpha}(\overline{\Omega^+(w)})$.
However, $u$ is not identically zero on $\sS(w)$ in general, so that $u$ is not a solution of the free boundary problem.
Next, the estimates for problem (\ref{generalFBP:1})--\eqref{generalFBP:2} and \eqref{bdryCondNozleExitU} in $\Omega^+(w)$  show that
$\| u-(\varphi^--\varphi_0^+)\|_{C^{2,\alpha}(\overline{\Omega^+(w)})}$ is small.
Then we extend $u$ to the whole domain $\RDom$ so that
$\|u-(\varphi^--\varphi_0^+)\|_{C^{2,\alpha}(\overline{\RDom})}$ is small.
This defines the iteration map: $w \mapsto u$.
The fixed point $u=w$ of this process determines a solution of
the free boundary problem, since $u$ is a solution
of (\ref{generalFBP:1})--\eqref{generalFBP:2} and \eqref{bdryCondNozleExitU}
in $\Omega^+(u)$, and $u$
satisfies $u=w>0$ on $\Omega^+(u)=\Omega^+(w):=\{w>0\}$ and $u=w=0$ on
$\sS\defd\partial\Omega^+(w)\setminus\partial\RDom$.
Then it remains to show the existence of a fixed point.
Since the right-hand side of the free boundary condition (\ref{generalFBP:2}) depends on $\nnu$,
we need to exploit the structure of our problem, in addition to the elliptic estimates,
to obtain the better estimates for the iteration and prove the existence of a fixed point.
More precisely, the nonlinear method can be described in the following
five steps:

\medskip
{1.}
 {\bf Iteration set}. Let $M \geq 1$. Set
\begin{equation}\label{defSetK}
\begin{array}{ll} \displaystyle
\setK \defd
\big\{
w \in C^{2,\alpha}(\overline{\RDom})
  \;\; : \;\;
\mbox{$w$ satisfies \eqref{symmetricFunc} and
$\|w-(\varphi^--\varphi^+_0)\|_{2,\alpha,\RDom} \leq M\epsP$}
\big\},
\end{array}
\end{equation}
where $\varphi^+_0(\x)= q^+x_d$.
Then $\setK$ is convex and
compact in $C^{2,\beta}(\RDom)$ for $0<\beta<\alpha$.

Let $w\in\setK$.
Since $q^->q^+$, it follows that, if
\begin{equation}\label{smallnessTheta_1}
\epsP\leq \frac{q^-- q^+}{10(M+1)},
\end{equation}
then combining (\ref{smallPert}) and (\ref{defSetK}) with \eqref{smallnessTheta_1}
implies
\begin{equation}\label{nondegeneracy}
w_{x_d}(\x)\geq \frac{q^-- q^+}{2}>0.
\end{equation}
By the implicit function theorem,
$\Omega^+(w)\defd \{w(\x)>0\}\cap \RDom$
has the form:
\begin{equation}\label{OmegaPL-f}
\Omega^+(w)
=\{x_d>f(\xpr)\;:\; \xpr\in \bTpr\},
\qquad
\|f\|_{2,\alpha,\bTpr} \leq CM\epsP<1,
\end{equation}
where $C$ depends on $q^--q^+$, and the last inequality is obtained by choosing small $\sigma$.
The corresponding unit normal on $\sS(w) \defd \{x_d=f(\xpr)\}$ is
$$
\bn(\x')=\frac{(-D_{\x'}f(\x'),1)}{\sqrt{1+|D_{\x'}f(\x')|^2}}
\in C^{1,\alpha}(\bTpr;\Sn)
$$
with
\begin{equation}\label{nuOnSetK_Gamma}
\begin{array}{l}
\displaystyle\|\bn-\bnNot\|_{1, \alpha, \bRnpr}
\leq CM\epsP,
\end{array}
\end{equation}
where $\bnNot$ is defined by
\begin{equation}\label{defNuNot}
\bnNot\defd
{\frac{D(\varphi^-_0-\varphi^+_0)}{|D(\varphi^-_0-\varphi^+_0)|}}= (0, \cdots, 0, 1 )^\top.
\end{equation}
Also,  $\bn(\cdot)$ can be considered as a function
on $\sS(w)$.
Since $\Omega^+(w)=\{w(x)>0\}\cap \RDom$,
from the definition of $f(\x')$ in \eqref{OmegaPL-f},
it follows that, for $\x\in \sS(w)$,
\begin{equation}\label{nuOnOmega}
\bn(\x) =\frac{Dw(\x)}{|Dw(\x)|}.
\end{equation}
By the definition of $\setK$ and
(\ref{smallnessTheta_1}) with (\ref{smallPert}),
$\bn(\x)$ can be extended to $\RDom$ via formula (\ref{nuOnOmega})
and
\begin{equation}\label{nuOnSetK}
\|\bn-\bnNot\|_{1, \alpha, \RDom} \leq CM\epsP
\end{equation}
with $C=C(q^+, q^-)$.
Motivated by the free boundary condition (\ref{FBConditionConormTr}),
we define a function $G_w$ on $\RDom$:
\begin{equation}\label{FBCond_u_2_k}
G_w(\x):=
\big(\rho(|D\varphi^-(\x)|^2)- \Trho(|D\varphi^-(\x)|^2)\big)D\varphi^-(\x)\cdot\bn(\x),
\end{equation}
where $\bn(\cdot)$ is defined by (\ref{nuOnOmega}).

We now solve the following fixed boundary value problem for $u$ in domain
$\Omega^+(w)$:
\begin{align}
&\divg A(Du,\x)=F(\x)&&\mbox{ in $\Omega^+\defd \{w>0\}$},
              \label{mixedProblemEq} \\
&A(Du,\x)\cdot \bn = G_w(\x)
&&\mbox{ on $\sS(w)\defd \partial\Omega^+(w)\setminus \partial \RDom$},
\label{mixedProblemBdryCond_1}\\
&u= \varphi^- - q^+
&&\mbox{ on $ \{x_d=1\}=\partial \Omega^+(w) \setminus \sS(w) $},
\label{mixedProblemBdryCond_2}
\end{align}
and show that its unique solution $u$ can be extended to the whole
domain $\RDom$ so that $u\in\setK$.

\bigskip
2. {\bf Existence and uniqueness of the solution for the
fixed boundary value
problem \eqref{mixedProblemEq}--\eqref{mixedProblemBdryCond_2}}.
We establish the existence and uniqueness of solution $u$
for problem (\ref{mixedProblemEq})--(\ref{mixedProblemBdryCond_2})
and show that $u$ is close
in $C^{2, \alpha}(\overline{\Omega^+(w)})$
to the unperturbed subsonic solution $\varphi^--\varphi^+_0$:
For $M\geq 1$, there is $\epsP_0>0$,
depending only on $(M, q^+, d, \gamma, \Omega)$,
such that, if $\epsP\in (0, \epsP_0)$,
$\varphi^-$ satisfies \eqref{smallPert}, and  $w\in \setK$,
there exists a unique solution
$u\in C^{2, \alpha}(\overline{\Omega^+(w)})$
of problem \eqref{mixedProblemEq}--\eqref{mixedProblemBdryCond_2}
satisfying \eqref{symmetricFunc} and
\begin{equation}\label{UcloseToPhi_0_Est}
\|u-(\varphi^--\varphi_0^+)\|_{2, \alpha,\Omega^+(w)}
\leq C\epsP,
\end{equation}
where $C$ depends only on  $(q^+, d, \gamma,\Omega)$
and is independent of $M$,
$w\in\setK$, and $\epsP\in (0, \epsP_0)$.

To achieve this, it requires to combine the existence arguments
with careful Schauder estimates for nonlinear oblique boundary value problems
for nonlinear elliptic PDEs, based on the results in
\cite{GilbargTrudinger,Lieberman86,LiebermanTrudinger,Trudinger85}
and the references cited therein.
Moreover, the independence of $C$ from $M$ is achieved by employing a cancellation based
on the structure in \eqref{mixedProblemBdryCond_1} with the explicit expressions
of $A(Du, \x)$ and $G_w(\x)$ and on the Rankine-Hugoniot condition
for the background solution.

\bigskip
3. {\bf Construction and continuity of the iteration map}.
We now construct the iteration map by an extension of the
unique solution of (\ref{mixedProblemEq})--(\ref{mixedProblemBdryCond_2}), which
satisfies (\ref{UcloseToPhi_0_Est}), and show the continuity
of the iteration map:
Let $w\in \setK$, and let $u(\x)$ be a solution of problem
\eqref{mixedProblemEq}--\eqref{mixedProblemBdryCond_2} in domain
$\Omega^+(w)$ established in Step 2 above. Then $u(\x)$ can
be extended to the whole domain $\RDom$ in such a way that this extension,
denoted as $\ExCFn_w u(\x)$, satisfies the following two
properties{\rm :}
\begin{enumerate}[{\rm (i)}]
\item \label{Uextended_EstK}
There exists $C_0>0$, which depends only on  $(q^+, d, \gamma, \Omega)$
and is independent
of $(M, \epsP)$ and $w(\x)$, such that
\begin{equation}\label{UcloseToPhi_0_Est_ext}
\| \ExCFn_w u -(\varphi^--\varphi_0^+)\|_{2, \alpha, \RDom}
\leq C_0\epsP.
\end{equation}
\item \label{Uextended_Cont}
Let $\beta \in (0, \alpha)$. Let $w_j \in \setK$ converge in
$C^{2, \beta}(\overline{\RDom})$
to $w \in  \setK$. Let
$u_j\in C^{2, \alpha}(\overline{\Omega^+(w_j)})$ and
$u \in C^{2, \alpha}(\overline{\Omega^+(w)})$
be the solutions of
problems \eqref{mixedProblemEq}--\eqref{mixedProblemBdryCond_2}
for $w_j(x)$ and $w(x)$, respectively.
Then $\ExCFn_{w_j} u_j \rightarrow \ExCFn_w u$ in
$C^{2, \beta}(\overline{\RDom})$.
\end{enumerate}

\smallskip
Define
the iteration map $\mapJ: \setK \rightarrow C^{2, \alpha}(\overline{\RDom})$
by
\begin{equation}\label{defMapJ}
\mapJ w\defd {\mathcal P}_w u,
\end{equation}
where $u(\x)$ is the unique solution of
problem (\ref{mixedProblemEq})--(\ref{mixedProblemBdryCond_2})
for $w(\x)$. By (ii),
$\mapJ$ is continuous  in the $C^{2, \beta}(\overline{\RDom})$--norm
for any positive $\beta<\alpha$.

Now we denote by $u(\x)$ both the function $u(\x)$
in $\Omega^+(w)$ and its extension $\ExCFn_w u(\x)$.
Choose $M$ to be the constant $C_0$ from (\ref{UcloseToPhi_0_Est_ext}).
Then, for $w\in\setK$, we see that $u\defd \mapJ w \in\setK$
if $\epsP>0$ is sufficiently
small, depending only on  $(q^+, d, \gamma, \Omega)$,
since $M$ is now fixed.
Thus, (\ref{defMapJ}) defines the iteration map
$\mapJ: \setK \rightarrow \setK$ and,
from \eqref{UcloseToPhi_0_Est_ext},
$\mapJ$ is continuous on $\setK$ in the $C^{2, \beta}(\overline{\RDom})$--norm
for any positive $\beta<\alpha$.

\medskip
4. {\bf Existence of a fixed point of the iteration map}.
We then prove the existence of solutions
of the free boundary problem, Problem \ref{problem-b}.

First, in order to solve Problem \ref{problem-c},
we seek a fixed point of map $\mapJ$.
 We use the Schauder fixed point
theorem ({\it cf}. Gilbarg-Trudinger \cite[Theorem 11.1]{GilbargTrudinger})
in the following setting:

Let $\epsP>0$ satisfy the conditions in Step 2 above.
Let $\beta \in (0, \alpha)$. Since $\RDom$ is a compact manifold with boundary
and $\setK$ is a  bounded convex subset of $C^{2,\alpha}(\overline{\RDom})$, it follows that
$\setK$ is a compact convex subset of $C^{2,\beta}(\overline{\RDom})$.
We have shown that $\mapJ(\setK)\subset\setK$, and $J$ is continuous
in the $C^{2,\beta}(\overline{\RDom})$--norm.
Then, by the Schauder fixed point theorem,
$\mapJ$ has a fixed point $\varphi \in \setK$.

If $u(\x)$ is such a fixed point, then
$$
\tilde{u}(\x)\defd\max(0, u(\x))
$$
is a classical solution
of Problem \ref{problem-c}, and $\sS(u)$ is its free boundary.

It follows that $\varphi:=\varphi^--\tilde u$ is a solution
of Problem \ref{problem-b},
provided that $\epsP$ is small enough so that (\ref{UcloseToPhi_0_Est})
implies that $|D\varphi|=|D(\varphi^--u)|< p^1_{\rm sonic}-\varepsilon$ on $\Omega^+(u)$,
where $\varepsilon=\frac{p_{\rm sonic}-q^+}{2}$ defined in \S 2.2.1.
Indeed, then (\ref{goodTruncatEquat})
implies that $\varphi(\x)$
lies in the non-truncated region for equation
(\ref{PotenEulerCompresTrunc}). Note also that the boundary condition
$\varphi_\bn=0$ on $\,\partial \Sqr \times [-1, 1]$ is satisfied
because $u$ and $\varphi^-$ satisfy
\eqref{symmetricFunc} on $\bTpr\times [-1, 1]$.

For such values of $\epsP$, if $\varphi^-(\x)$ is a supersonic solution  of \eqref{PotenEulerCompres}
satisfying the conditions stated in Problem {\rm \ref{problem-a}},
the defined function $\varphi(x)$ is a solution of Problem \ref{problem-a}.
Indeed,
$|D\varphi|=|D(\varphi^--\tilde{u})|< p^1_{\rm sonic}-\varepsilon$
on $\Omega^+(\varphi):=\{\varphi<\varphi^-\}=\{\tilde{u}(\x)>0\}$
since $\tilde{u}=u$ on $\Omega^+(\tilde{u})$ and $|D\varphi|=|D\varphi^-|> p^1_{\rm sonic}$
on $\Omega\setminus \Omega^+(\varphi)$,  equation \eqref{PotenEulerCompres} is
satisfied in both $\Omega^+(\varphi)$ and $\Omega\setminus \Omega^+(\varphi)$,
and the Rankine-Hugoniot conditions \eqref{continuity}--\eqref{FBCondition-1} are satisfied
on $\sS=\partial \Omega^+(\varphi)\setminus \partial \Omega$.

\smallskip
This completes the construction of the global solution.
The uniqueness and stability of the solution of the free boundary
problem are obtained by using the regularity
and nondegeneracy of solutions.

\vspace{5pt}
\begin{remark}\label{rmk-1}
For clarity, in this section, we focus on the simplest setup of the domain as $\Omega=(0,1)^{d-1}\times (-1, 1)$,
which can be extended directly to $\Omega_R=\Pi_{j=1}^{d-1}(0,a_j)\times (-1, R)$ for any $R>0$, then to
$\Omega_\infty=\Pi_{j=1}^{d-1}(0,a_j)\times (-1, \infty)$ by analyzing the asymptotic behavior of the solution when $R\to \infty$,
as well as to $\Omega=\R^{d-1}\times (-1, \infty)${\rm ;} see Chen-Feldman {\rm \cite{CF-JAMS2003,ChenFeldman2,ChenFeldman4}}.
See also {\rm Chen \cite{Sxchen3}} for the extension to the isentropic Euler case.

If the hyperbolic phase is $C^\infty$, then the solution and
its corresponding free boundary in {\rm Theorem \ref{mainFBPTh}}
are also $C^\infty$.
Furthermore, our results can be extended to the problem
with a steady $C^{1,\alpha}$--perturbation of the upstream supersonic flow
and/or  general Dirichlet data $h(\x'), \x'\in \R^{d-1},$
at $x_d=1$ satisfying
$$
\|h-\varphi^+_0\|_{1,\alpha, \R^{d-1}}\le C\sigma
\qquad \mbox{for $\alpha\in (0,1)$}.
$$
Also, the Dirichlet data in Problem {\rm \ref{problem-b}} may be replaced by
the corresponding Neumann data satisfying the global solvability
condition.

The global uniqueness of piecewise constant transonic shocks in straight ducts
modulo translations was analyzed in {\rm \cite{CY,Fang-Liu-Yuan}}.
\end{remark}

\begin{remark}\label{rmk-2}
The domains in the setup of {\rm Problems 2.1--2.2}
have also been extended to M-D infinite nozzles
of arbitrary cross-section in Chen-Feldman {\rm \cite{ChenFeldman3Arch}}{\rm ;}
also see Xin-Yin {\rm \cite{XY}}, Yuan {\rm \cite{Yuan1}}, and the references cited therein
for the $2$-D case with the downstream pressure exit.
For the analysis of geometric effects of the nozzles on the
uniqueness and stability of steady transonic shocks, see
{\rm \cite{BaeFeldman,CY2,Liu-Yuan,Liu-Xu-Yuan,LXY}}
and the references cited therein.
\end{remark}

\begin{remark}\label{rmk-3}
The iteration scheme can also be reformulated in a way such that the free boundary normal $\nnu$ is unknown
in the iteration by replacing the known function $w$ in \eqref{nuOnOmega} by the unknown $u$, that is,
by replacing $\bn(\x)$ in \eqref{nuOnOmega} via
\begin{equation}\label{nuOnOmegaU}
\bn(\x)
=\frac{Du(\x)}{|Du(\x)|}.
\end{equation}
Note that \eqref{nuOnOmegaU} coincides with \eqref{nuOnOmega} at the fixed point $u=w$,
i.e., defines the normal to $\sS$.
Using expression \eqref{nuOnOmegaU} for $\bn$ in the iteration boundary condition,
we improve the regularity and structure of the boundary condition;
in particular, it is made
independent of the regularity and constants in the iteration set.
This is useful in many cases, see e.g. {\rm \cite{ChenFeldman}}. Moreover,
this allows us to obtain the compactness of the iteration map,
which has been used in {\rm \cite{CF-book2018}}.
\end{remark}

This nonlinear method and related ideas and techniques described above for free boundary problems
have played a key role in many recent developments in the analysis
of M-D transonic shock problems, as shown in {\rm \S 3}--{\rm \S 5} below.

\smallskip
\section{Two-Dimensional Transonic Shocks and Free Boundary Problems for the Steady Full Euler Equations}

\newcommand{\LU}{\mathcal{U}}
We now describe how the nonlinear method
and related ideas and techniques
presented in \S 2 can be applied
to establish the existence, stability, and asymptotic behavior of
2-D steady transonic flows with transonic shocks past curved wedges
for the full Euler equations, by reformulating the problems as free boundary problems,
via two different approaches.

The 2-D steady Euler
equations for polytropic gases are of the form
({\it cf.} \cite{CF-book2018,CF}):
\begin{equation}\label{Euler1}
\begin{cases}
{\rm div}(\rho \uu)=0, \\[0.5mm]
{\rm div}(\rho{\uu\otimes\uu})+\nabla p=0,\\[0.5mm]
{\rm div}\big(\rho\uu(E+ \frac{p}{\rho})\big)=0,
\end{cases}
\end{equation}
where
$\uu=(u_1, u_2)$ is the velocity, $\rho$ the density, $p$ the pressure, and
$
E=\frac{1}{2}|\uu|^2+e
$
the total energy with internal energy $e$.

Choose pressure $p$ and density $\rho$ as the independent thermodynamical variables.
Then the constitutive relations can be written as
$$
(e,T, S)=(e(p, \rho), T(p,\rho), S(p,\rho))
$$
 governed by
$$
T dS=de-\frac{p}{\rho^2}d\rho,
$$
where $T$ and $S$ represent the temperature and the entropy, respectively.
For a polytropic gas,
\begin{equation}\label{gas-2}
e=e(p,\rho)=\frac{p}{(\gamma-1)\rho},\quad
T=T(p,\rho)=\frac{p}{(\gamma-1)c_v\rho}, \quad
S=S(p,\rho)=c_v\ln (\frac{p}{\kappa\rho}),
\end{equation}
where
$\gamma>1$ is the adiabatic exponent,
$c_v>0$ the specific heat at constant volume,
and $\kappa>0$  any constant under scaling.

System \eqref{Euler1} can be written as a first-order system of conservation laws:
\begin{equation}\label{FullSystInUvar}
\partial_{x_1} F(U) +\partial_{x_2} G(U)=0, \,\quad\qquad U=(\uu, p, \rho)\in \R^4.
\end{equation}
Solving ${\rm det}(\lambda \nabla_UF(U)-\nabla_UG(U))=0$ for $\lambda$, we obtain four eigenvalues:
$$
\lambda_1=\lambda_2=\frac{u_2}{u_1}, \qquad\,\,
\lambda_j=\frac{u_1u_2+(-1)^jc\sqrt{|\uu|^2-c^2}}{u^2_1-c^2}\,\,\,\, \mbox{for $j=3,4$},
$$
where
\begin{equation}\label{1.4aa}
c=\sqrt{\frac{\gamma p}{\rho}}
\end{equation}
is the sonic speed of the flow for a polytropic gas.

The repeated eigenvalues $\lambda_1$ and $\lambda_2$ are real
and correspond to the two linear degenerate
characteristic families which generate vortex sheets and entropy waves, respectively.
The eigenvalues $\lambda_3$ and $\lambda_4$ are real when the flow is supersonic ({\it i.e.}, $|\uu|> c$),
and complex when the flow is subsonic ({\it i.e.}, $|\uu|<c$) in which case the corresponding two
equations are elliptic.

For a transonic flow in which both the supersonic and subsonic phases occur in the flow,
system \eqref{Euler1}
is of mixed-composite hyperbolic-elliptic type, which consists
of two equations of mixed elliptic-hyperbolic type
and two equations of hyperbolic type ({\it i.e.}, two transport-type equations).

\smallskip
In the regimes with $\rho |\uu|>0$, from the first equation in \eqref{Euler1},
in any domain containing the origin,
there exists a unique
stream function $\psi$
such that
\begin{equation}\label{psi-potential}
D\psi= (-\rho u_2, \rho u_1)
\qquad\,\, \mbox{with $\psi(\mathbf{0})=0$}.
\end{equation}
We use the following Lagrangian coordinate
transformation:
\begin{equation}\label{def-coord}
(x_1,x_2)\, \to\, (y_1,y_2)= (x_1, \psi(x_1, x_2)),
\end{equation}
under which the original curved streamlines become straight.
In the new coordinates $\yy=(y_1, y_2)$, we still denote the unknown
variables $U(\xx(\yy))$ by $U(\yy)$ for simplicity of notation.
Then the original Euler equations in \eqref{Euler1} become the
following equations in divergence form:
\begin{eqnarray}
  &&\big(\frac{1}{\rho u_1}\big)_{y_1}
    -\big(\frac{u_2}{u_1} \big)_{y_2}= 0,\label{eqn-euler1}\\[0.5mm]
  &&\big( u_1 + \frac{p}{\rho u_1}\big)_{y_1}
   - \big( \frac{p u_2}{u_1}\big)_{y_2}= 0, \label{eqn-euler2}\\[0.5mm]
  &&(u_2)_{y_1} + p_{y_2}= 0,\label{eqn-euler3} \\[0.5mm]
  &&\big( \frac{1}{2}|\uu|^2 + \frac{\gc p}{(\gc-1)\rho} \big)_{y_1}=
  0.\label{eqn-euler4}
\end{eqnarray}
One of the advantages of the Lagrangian coordinates is to straighten the streamlines
so that the streamline may be employed as one of the coordinates to simplify
the formulations, since the Bernoulli variable is conserved
along the streamlines.
Note that the entropy is also conserved along the streamlines in the continuous part of the flow.

\subsection{Steady Supersonic Flow onto Solid Wedges and Free Boundary Problems}
\begin{figure}
 	\centering
 	\includegraphics[height=35mm]{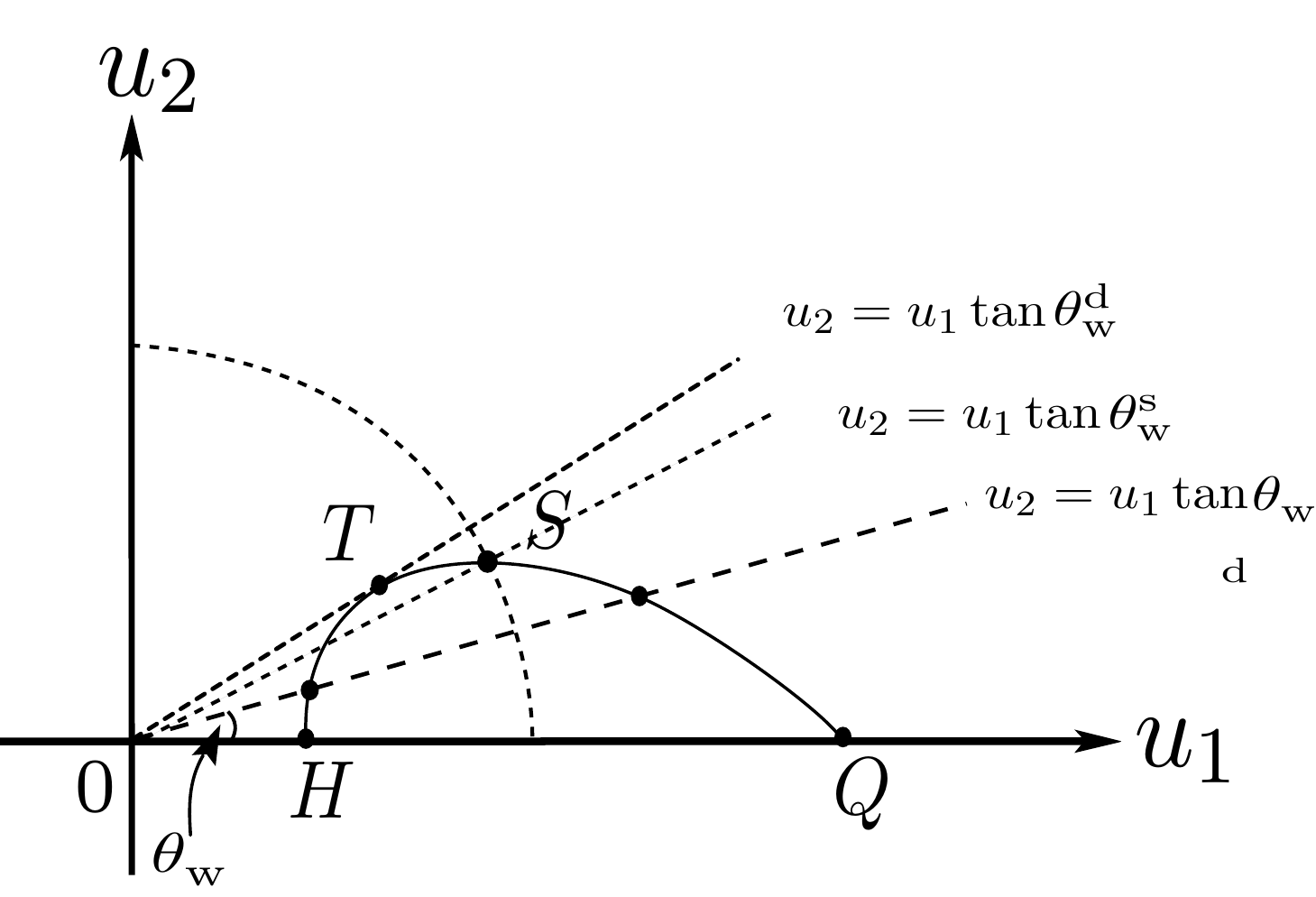}
\hspace{12mm}
 		\includegraphics[height=35mm]{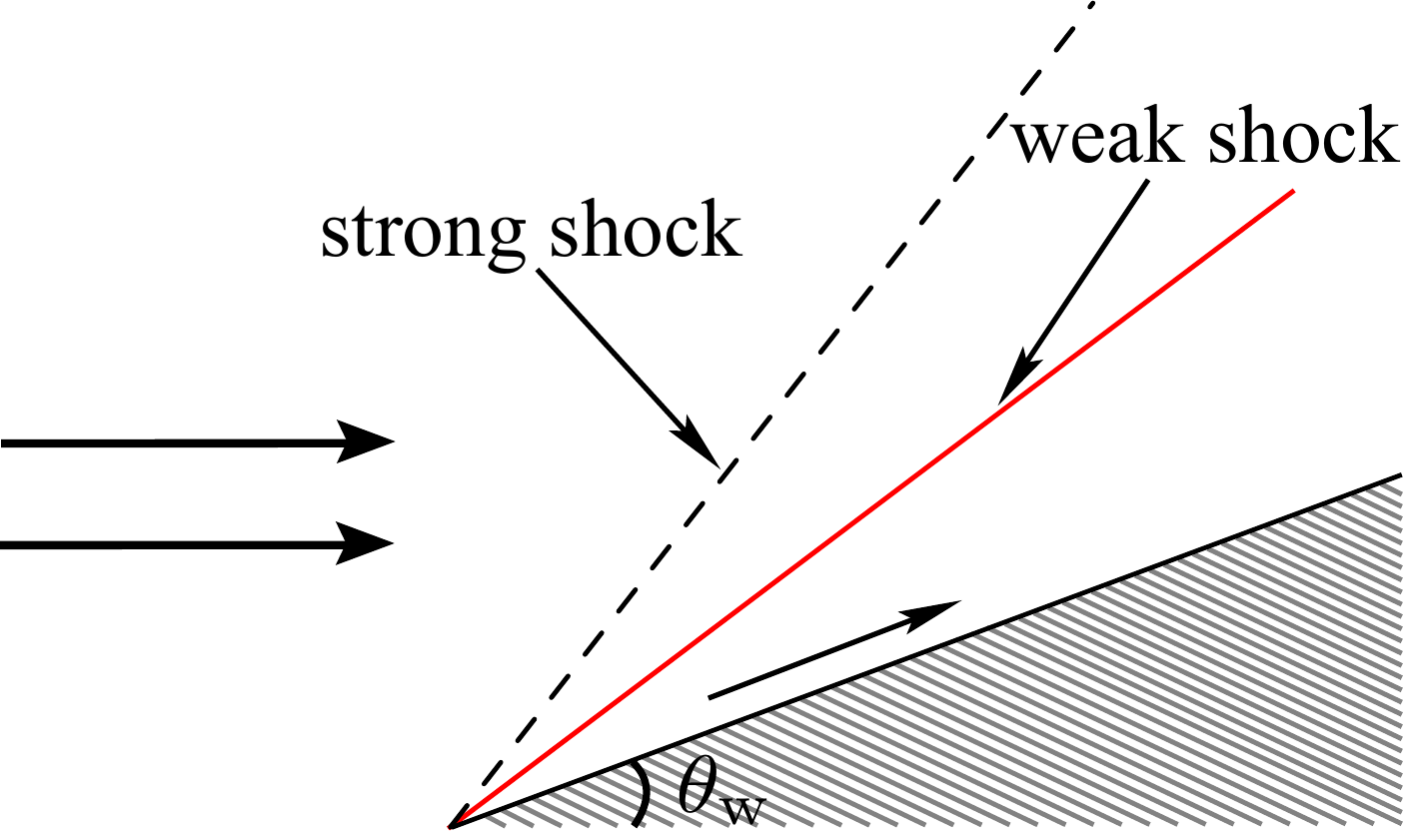}
		\caption{The shock polar in the $\uu$-plane and uniform steady (weak/strong) shock flows (see \cite{Chen2})}
		\label{Figure2}
 \end{figure}

For an upstream steady uniform supersonic
flow past a symmetric straight-sided wedge (see Fig. \ref{Figure2}):
\begin{equation}\label{wedge-1}
W:=\{\xx=(x_1,x_2)\in \R^2\,:\, |x_2|<x_1\tan \theta_{\rm w}, x_1>0\}
\end{equation}
whose angle $\theta_{\rm w}$
is less than
the detachment angle $\theta^{\rm d}_{\rm w}$,
there exists an oblique shock emanating from the
wedge vertex.
Since the upper and lower subsonic regions do not interact with each other,
it suffices to study the upper part.
More precisely, if the upstream steady flow is a uniform supersonic state,
we can find the corresponding
constant downstream flow along the straight-sided upper wedge boundary,
together with a straight shock separating
the two states.
The downstream flow is determined by the shock polar
whose states in the phase space are governed by
the Rankine-Hugoniot conditions
and the entropy condition; see Fig.~\ref{Figure2}.
Indeed, Prandtl in \cite{Prandtl} first employed the shock polar analysis
to show that there are two possible
steady oblique shock configurations when the wedge angle $\theta_{\rm w}$ is less than
the detachment angle $\theta_{\rm w}^{\rm d}$ --- The steady weak shock with supersonic or subsonic downstream flow
(determined by the wedge angle that is less or larger than
the sonic angle $\theta^{\rm s}_{\rm w}$)
and the steady strong shock
with subsonic downstream flow, both of which satisfy the entropy
condition, provided that no additional conditions are assigned at downstream.
See also \cite{Busemann,Chen2,CF,Meyer,Prandtl} and the references cited therein.

The fundamental issue -- whether one or both of the steady weak and strong shocks are physically
admissible -- has been vigorously debated over the past seven decades
({\it cf.} \cite{Chen2,CF,Liu,Neumann,Neumann10, Serre2}).
Experimental and numerical results have strongly indicated
that the steady weak shock
solution would be physically admissible, as Prandtl conjectured in \cite{Prandtl}.
One natural approach to single out the physically admissible steady
shock solutions is via the
stability analysis: the stable ones are physical.
See Courant-Friedrichs \cite{CF} and von Neumann \cite{Neumann,Neumann10}; see also \cite{Liu,Serre2}.

A piecewise smooth solution $U=(\uu, p, \rho)\in \R^4$
separated by a front $\mathcal{S}:=\{\xx\,:\,  x_2=\sigma(x_1), x_1\ge 0\}$
becomes a weak solution of the Euler equations \eqref{Euler1} as in \S 2.1 if and only if the
Rankine-Hugoniot conditions are satisfied along
$\mathcal{S}$:
\begin{equation} \label{con-RH}
\begin{cases}
\sigma'(x_1)[\,\rho u_1\,]=[\,\rho u_2\,],\\[1mm]
\sigma'(x_1)[\,\rho u_1^2 + p\,]= [\,\rho u_1 u_2\,],\\[1mm]
\sigma'(x_1)[\,\rho u_1 u_2\,] =  [\,\rho u_2^2 + p\,],\\[1mm]
\sigma'(x_1)[\,\rho u_1(E+\frac{p}{\rho})\,] =[\,\rho u_2(E+\frac{p}{\rho})\,],
\end{cases}
\end{equation}
where $[\, \cdot\, ]$ denotes
the jump between the quantities of two states across front $\mathcal{S}$ as before.

Such a front $\mathcal{S}$ is called a shock if the entropy condition holds along $\mathcal{S}$:
{\it The density increases in the fluid direction across $\mathcal{S}$.}

\smallskip
For given state $U^-$, all states $U$ that can be connected
with $U^-$ through the relations in
\eqref{con-RH}
form a curve in the state space $\R^4$;
the part of the curve whose states satisfy the entropy condition is called the {\it shock polar}.
The projection of the shock polar onto the $\uu$--plane is shown in Fig.~\ref{Figure2}.
In particular,  for an upstream uniform horizontal flow
$U_0^-=(u_{10}^-,0, p_0^-, \rho_0^-)$ past the upper part of a straight-sided wedge whose angle
is $\theta_{\rm w}$, the downstream constant flow can be
determined by the shock polar (see Fig.~\ref{Figure2}).
Note that the downstream flow must be parallel to the wedge, and the upstream flow
is parallel to the axis of wedge, so the angle between the upstream and downstream flow is equal to the (half) wedge angle.
According to the shock polar, the two flow angles (or, equivalently, wedge angles) are
particularly important:

One is the detachment angle $\theta_{\rm w}^{\rm d}$ such that
line $u_2=u_1 \tan\theta_{\rm w}^{\rm d}$
is tangential to the shock polar at point $T$ and there is no intersection between
line $u_2=u_1 \tan \theta_{\rm w}$ and
the shock polar when $\theta_{\rm w}>\theta_{\rm w}^{\rm d}$.
For any wedge angle $\theta_{\rm w}\in (0, \theta_{\rm w}^{\rm d})$,
there are two intersection points
of line  $u_2=u_1 \tan \theta_{\rm w}$ and
the shock polar: one intersection point is on arc $\wideparen{TH}$ which determines velocity
$\uu^{\rm sg}=(u_1^{\rm sg}, u_2^{\rm sg})$ of the downstream flow corresponding
to the strong shock, and the other intersection point is on arc $\wideparen{TQ}$  which determines velocity
$\uu^{\rm wk}=(u_1^{\rm wk}, u_2^{\rm wk})$ of the downstream flow corresponding
to the weak shock.
Thus, for any wedge angle $\theta_{\rm w}\in (0, \theta_{\rm w}^{\rm d})$,
the shock polar ensures the existence of two attached shocks at the wedge: strong versus weak.

Since each point on the shock polar defines a downstream flow that is a constant state,
we can use \eqref{1.4aa} to compute its sonic speed $c_0$ and then determine whether this downstream
state is subsonic or supersonic.
It can be shown that there exists the unique point $S$ on the shock polar
so that all downstream states are subsonic for the points on $\wideparen{HS} \setminus \{S\}$,
supersonic for the points of $\wideparen{SQ} \setminus \{S\}$,
and sonic for the state at $S$.
Moreover, $S$ lies in the interior of arc $TQ$.
Then, denoting by $\theta_{\rm w}^{\rm s}$ the angle corresponding to point $S$,
we see that $\theta_{\rm w}^{\rm s} < \theta_{\rm w}^{\rm d}$.
The wedge angle
$\theta_{\rm w}^{\rm s}$ is called the sonic angle.
Point $T$  divides arc $\wideparen{HS}$, which corresponds to the transonic shocks,
into the two open arcs
$\wideparen{TS}$ and $\wideparen{TH}$; see Fig.~\ref{Figure2}.
The nature of these two cases, as well as the case for arc $\wideparen{SQ}$,
is very different.
When the wedge angle $\theta_{\rm w}$ is between
$\theta_{\rm w}^{\rm s}$ and $\theta_{\rm w}^{\rm d}$,
there are two subsonic solutions (corresponding to the strong and weak shocks);
while, for the wedge angle
$\theta_{\rm w}$ is smaller than $\theta_{\rm w}^{\rm s}$,
there are one subsonic
solution (for the strong shock) and one supersonic solution  (for the weak shock).
Such an oblique shock $\mathcal{S}_0$ is  straight,
described by $x_2= s_0 x_1$ with $s_0$ as its slope.
The question is whether the steady oblique shock
solution is stable under a perturbation of both the upstream supersonic
flow and the wedge boundary.

\begin{figure}
 \centering
\includegraphics[height=47mm]{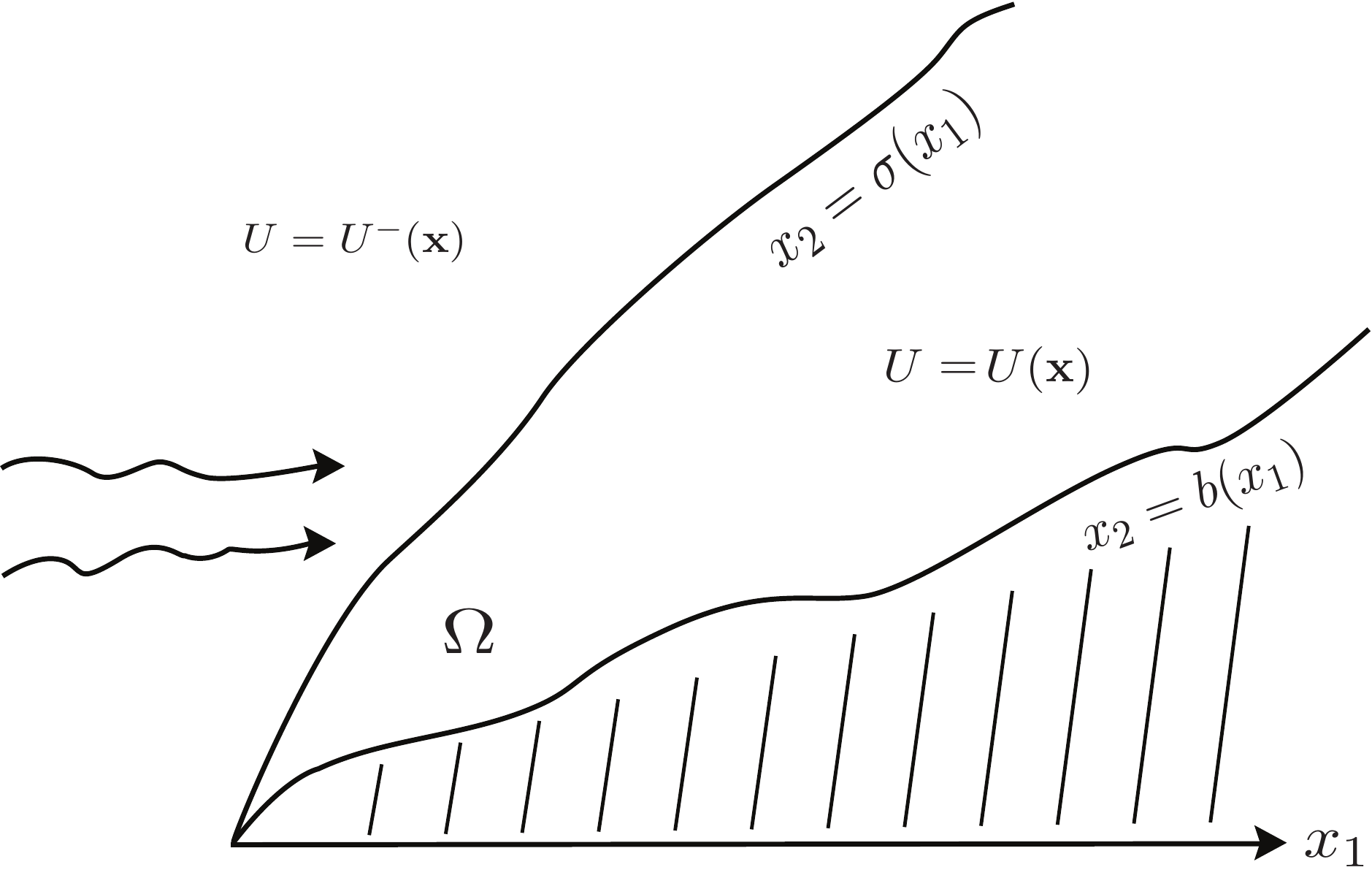}
\caption[]{The leading steady shock $x_2=\sigma(x_1)$ as a free boundary under the perturbation (see \cite{Chen2})}
\label{Figure2a}
\end{figure}

Since we are interested in determining the downstream flow, we can restrict the domain to the
first quadrant; see Fig. \ref{Figure2a}.
Fix a constant upstream flow $U^-_0$,
a wedge angle $\theta_{\rm w}\in (0, \theta_{\rm w}^{\rm d})$, and a constant downstream state $U_0^+$
which is one of the downstream states (weak or strong) determined by
the shock polar for the chosen upstream flow and wedge angle.
States $U^-_0$ and $U^+_0$ determine the oblique shock $x_2= s_0 x_1$,
and the transonic shock solution $U_0$ in $\{\xx\; : \; x_1>0,\,x_2>0\}\setminus W$ such that $U_0=U^-_0$
in $\Omega_0^-=\{\xx\in \R^2\,:\,  x_2>s_0 x_1, x_1>0\}$ and
$U_0=U^+_0$
in $\Omega_0^+=\{\xx\in \R^2\,:\, x_1\tan\theta_{\rm w}< x_2< s_0x_1, x_1>0\}$;
see Fig. \ref{Figure2}.
We refer to this solution as a constant transonic solution $(U^-_0,U^+_0)$.

Assume that the perturbed upstream flow $U^-_I$ is close to $U_0^-$
(so that $U^-_I$ is supersonic and almost horizontal)
and that the perturbed wedge is close
to a straight-sided wedge.  Then, for any suitable wedge angle (smaller
than the detachment angle), it is expected that there should be a shock
attached to the wedge vertex; see  Fig. \ref{Figure2a}.
We now use a function $b(x_1)\ge 0$ to describe the upper perturbed wedge boundary:
\begin{equation}\label{wall}
\partial{W}=\{\xx\in \R^2\,:\,  x_2= b(x_1),\ x_1>0\} \qquad\mbox{with} \; b(0)=0.
\end{equation}
Then the wedge problem can be formulated as the following problem:

\begin{problem}[Initial-Boundary Value Problem]
\label{pr-3.1}
Find a global solution of system \eqref{Euler1}
in $\Omega:=\{x_2>b(x_1), x_1>0\}$
such that the following conditions hold{\rm :}
\begin{enumerate}
\item[\rm (i)] Cauchy condition at $x_1=0${\rm :}
\begin{equation}
U|_{x_1=0}=U^-_I(x_2){\rm ;}
\end{equation}

\item[\rm (ii)] Boundary condition on $\partial{W}$ as the slip boundary{\rm :}
\begin{equation}\label{slipcon}
\uu\cdot \nnu_{\rm w}|_{\partial{W}}=0,
\end{equation}
where $\nnu_{\rm w}$ is the outer unit normal vector to $\partial{W}$.
\end{enumerate}
\end{problem}

Note that the background shock is the straight line
given by $x_2= \s_0(x_1)$ with $\s_0(x_1):=s_0 x_1$.
When the upstream steady supersonic
perturbation $U^-_I(x_2)$ at $x_1=0$
is suitably regular and small under some natural norm,
the upstream steady supersonic smooth solution $U^-(\xx)$ exists
in region $\Omega^- = \{\xx\,:\, x_2>\frac{s_0}{2} x_1, x_1\ge 0\}$,
beyond the background shock, and $U^-$ in $\Omega^-$
is still close to $U_0^-$.

Assume that the shock-front
$\mathcal{S}$ to be determined is
\begin{equation}\label{shock:1}
\sS=\{\xx \,: \, x_2 = \s(x_1),  \,x_1 \ge 0\} \qquad\mbox{with $\s(0)=0\, $ and $\s(x_1)>0$ for $x_1>0$}.
\end{equation}
The domain for the downstream flow behind $\mathcal{S}$ is denoted by
\begin{equation}\label{domain:1}
\Omega=\{\xx\in \R^2\,:\,  b(x_1)<x_2<\s(x_1), x_1>0\}.
\end{equation}
Then Problem \ref{pr-3.1} can be reformulated into
the following free boundary problem with $\sS$ as a free boundary:

\begin{problem}[Free Boundary Problem; see Fig.~\ref{Figure2a}]\label{pr-3.2}
Let $(U^-_0, U^+_0)$ be a constant transonic solution
for the wedge angle $\theta_{\rm w}\in (0, \theta_{\rm w}^{\rm d})$
with transonic shock
$\mathcal{S}_0:=\{x_2=\sigma_0(x_1)\,:\, x_1>0\}$ for $\sigma_0(x_1):=s_0x_1$.
For any upstream flow $U^-$ for system \eqref{Euler1} in domain $\Omega^-$
which is a small perturbation of $U^-_0$,
and any wedge boundary function $b(x_1)$ that is a  small perturbation of
$b_0(x_1)=x_1\tan\theta_{\rm w}$,
find a shock $\sS$  as a free boundary $x_2= \s(x_1)$
and a solution $U$ in $\Omega$,
which are small perturbations of $\mathcal{S}_0$ and $U^+_0$ respectively, such that
\begin{enumerate}[\rm (i)]
\item $U$ satisfies \eqref{Euler1} in domain $\Omega$,
\item  The slip condition \eqref{slipcon} holds  along the wedge boundary,
\item The Rankine-Hugoniot conditions in \eqref{con-RH} as free boundary conditions
hold along the transonic shock-front  $\CS$.
\end{enumerate}

There are three subcases based on $U^+_0${\rm :}
For a weak supersonic shock $\mathcal{S}_0$ given by $U^+_0$ corresponding to a supersonic state on arc $\wideparen{SQ}$,
we denote the problem by {\rm Problem \ref{pr-3.2}(WS)}{\rm ;}
for a weak transonic shock $\mathcal{S}_0$ given by $U^+_0$ corresponding to a subsonic state on arc $\wideparen{TS}$,
we denote the problem by {\rm Problem \ref{pr-3.2}(WT)}{\rm ;}
finally, for a strong transonic shock $\mathcal{S}_0$ given by $U^+_0$ corresponding to a subsonic state
on arc $\wideparen{TH}$, we denote the problem by {\rm Problem \ref{pr-3.2}(ST)}.
\end{problem}

In general, the uniqueness for the initial-boundary value  problem (Problem \ref{pr-3.1})
is not known (as it is a problem for a nonlinear system of a composite elliptic-hyperbolic type),
so it may not yet be excluded that {Problem \ref{pr-3.1}} has solutions which are not of steady oblique shock structure,
{\it i.e.}, are not solutions of {Problem \ref{pr-3.2}}.
On the other hand, the global solution of the free boundary problem (Problem \ref{pr-3.2})
provides the global structural stability of
the steady oblique shock, as well as more detailed structure of the solution.

Supersonic ({\it i.e.}, supersonic-supersonic) shocks correspond to
arc $\wideparen{SQ}$ which is a weaker shock (see Fig.~\ref{Figure2}).
The local stability of such shocks was first established in \cite{Gu,Li,Schaeffer}.
The global stability of the supersonic shocks
for potential flow past piecewise smooth perturbed curved wedges was established
in Zhang \cite{Zh2}; also see \cite{ChS3,CXY,ChenFang}
and the references therein.
The global stability and uniqueness of the supersonic shocks
for the full Euler equations, Problem \ref{pr-3.2}(WS),
were solved for more general perturbations of both the initial data and wedge boundary
even in $BV$ in Chen-Zhang-Zhu  \cite{ChenZhangZhu} and Chen-Li \cite{ChenLi2008}.

\smallskip
For transonic ({\it i.e.}, supersonic-subsonic) shocks,
the strong shock case corresponding to arc $\wideparen{TH}$
was first studied in Chen-Fang \cite{ChenFang} for the potential flow (see Fig.~\ref{Figure2}).
In Fang \cite{Fang},  the full Euler equations were studied
with a uniform Bernoulli constant
for both weak and strong transonic shocks.
Because the framework is a weighted Sobolev space, the asymptotic behavior of the shock
slope or subsonic solution was not derived.
In Yin-Zhou \cite{YinZhou}, the H\"older norms were used for
the estimates of solutions of the full Euler equations
with the assumption on the sharpness of the wedge angle,
which means that the subsonic state
is near point $H$ in the shock polar,
by Approach I introduced first in \cite{CCF} which is described in \S \ref{Sect3.2} below.
In Chen-Chen-Feldman \cite{Chen-Chen-Feldman},
the weaker transonic shock, which corresponds to
arc $\wideparen{TS}$, was first investigated by Approach I as described in \S \ref{Sect3.2} below.
Then, in \cite{CCF3},  the weak and strong transonic shocks,
which correspond to arcs $\wideparen{TS}$ and $\wideparen{TH}$, respectively,
were solved, by Approach II which is described in \S \ref{Sect3.3} below,
so that the existence, uniqueness, stability, and asymptotic
behavior of subsonic solutions of both Problem \ref{pr-3.2}(WT) and Problem \ref{pr-3.2}(ST)
in a weighted H\"{o}lder space were obtained.

We now describe two approaches for the wedge problem, based on the nonlinear method
and related ideas and techniques presented in \S 2.
First, we need to introduce the weighed H\"older norms
in the subsonic domain $\Omega$,
where $\Omega$ is either a truncated triangular domain or an unbounded domain with the vertex
at origin $\O$ and one side as the wedge boundary.
There are two weights: One is the distance function to origin $\O$ and the other is to
the wedge boundary $\partial W$.
For any $\xx, \xx'\in \Omega$,  define
\begin{align*}
&\gd^{\rm o}_\xx := \min(|\xx|,1),\quad \gd^{\rm o}_{\xx,\xx'} := \min (\gd^{\rm o}_\xx, \gd^{\rm o}_{\xx'}),
\quad \gd^{\rm w}_\xx := \min(\textrm{dist}(\xx,\partial W),1), \quad \gd^{\rm w}_{\xx,\xx'}
   := \min (\gd_\xx^{\rm w},\gd_{\xx'}^{\rm w}), \\[1mm]
&\gD_\xx := |\xx|+1, \quad \gD_{\xx, \xx'}:=\min(\gD_\xx, \gD_{\xx'}),
\quad \widetilde{\gD}_\xx := \textrm{dist}(\xx,\partial W)+1,
\quad \widetilde{\gD}_{\xx, \xx'}:=
\min(\widetilde{\gD}_\xx, \widetilde{\gD}_{\xx'}).
\end{align*}
Let $\ga \in (0,1)$  and $l_1,l_2, \gamma_1, \g_2 \in \R$ with $\g_1 \ge \g_2$,
and let $k\ge 0$ be an integer.
Let $\kk = (k_1, k_2)$ be an integer-valued vector, where $k_1, k_2 \ge 0$,
$|\kk|=k_1 +k_2$, and $D^{\kk}= \po_{x_1}^{k_1}\po_{x_2}^{k_2}$.
We define
\begin{align}
&[ f ]_{k,0;(l_1,l_2); \Omega}^{(\g_1;\O)(\g_2;\partial W)}
=\sup_{\begin{subarray}{c}
\xx\in \Omega\\
 |\kk|=k
\end{subarray}}
\begin{array}{l}
\big\{(\gd^{\rm o}_{\xx})^{\hat{\gamma}_0}
(\gd^{\rm w}_{\xx})^{\max\{k+\g_2,0\}}
\,\gD_\xx^{l_1} \widetilde{\gD}_\xx^{l_2+k} |D^\kk f(\xx)|\big\},
\end{array}\label{def-normk0}\\[2mm]
&{[ f ]}_{k,\ga;(l_1,l_2); \Omega}^{(\g_1;\O)(\g_2;\partial W)}
=\sup_{
\begin{subarray}{c}
   \xx, \xx'\in \Omega\\
 \xx \ne \xx', |\kk|=k
 \end{subarray}}
\Big\{
  (\gd^{\rm o}_{\xx,\xx'})^{\hat{\gamma}_\alpha}
  (\gd^{\rm w}_{\xx,\xx'})^{\max\{k+\ga+\g_2,0\}}
  \gD_{\xx,\xx'}^{l_1} \widetilde{\gD}_{\xx,\xx'}^{l_2+k+\ga}\frac{|D^\kk f(\xx)-D^\kk
f(\xx')|}{|\xx-\xx'|^\ga}
\Big\},
\label{def-normkga}
\\
& \|f\|_{k,\ga;(l_1,l_2); \Omega}^{(\g_1;\O)(\g_2;\partial W)}
= \sum_{i=0}^k {[f]}_{i,0; (l_1,l_2);\Omega}^{(\g_1;\O)(\g_2;\partial W)}
+ {[f]}_{k,\ga;(l_1,l_2);\Omega}^{(\g_1;\O)(\g_2;\partial W)}, \label{def-norm}
\end{align}
where $\hat{\gamma}_\beta=\max\{\g_1 + \min\{k+\beta,-\g_2\},0\}$ for $\beta\in [0, 1)$.
Similarly, the H\"{o}lder norms for a function of one variable
on $\mathbb{R}^+:=(0,\infty)$ with the weight near $\{0\}$ and the decay at infinity
are denoted by $\|f\|_{k,\ga;(l);\mathbb{R}^+}^{(\gamma_2; 0)}$.

For a vector-valued function $\ff=(f_1, f_2, \cdots,f_n )$, we
define
\[
 \|\ff\|_{k,\ga;(l_1,l_2); \Omega}^{(\g_1;\O)(\g_2;\partial W)} =
\sum_{i=1}^n  \|f_i\|_{k,\ga;(l_1,l_2); \Omega}^{(\g_1;\O)(\g_2;\partial W)}.
\]
 Let
\begin{equation}\label{def-C}
C^{k,\ga;(l_1,l_2)}_{(\g_1;\O)(\g_2;\partial W)}(\Omega)
= \left\{ \ff\,:\,\|\ff\|_{k,\ga;(l_1,l_2); \Omega}^{(\g_1;\O)(\g_2;\partial W)} <\infty \right\}.
\end{equation}

The requirement $\g_1 \ge \g_2$ in the definition above means that the regularity up
to the wedge boundary is no worse than the regularity up to the wedge vertex.
When $\g_1 = \g_2$, the $\d^{\rm o}$--terms disappear so that $(\g_1,\O)$ can be dropped
in the superscript.
If there is no weight $(\g_2,\pw)$ in the superscript, the $\d$--terms for the weights
should be understood as $(\d_{\xx}^{\rm o})^{\max\{k+\g_1,0\}}$
and $(\d_{\xx}^{\rm o})^{\max\{k + \ga +\g_1,0\}}$
in \eqref{def-normk0} and \eqref{def-normkga}, respectively.
Moreover, when no weight appears in the superscripts of the seminorms
in \eqref{def-normk0} and \eqref{def-normkga},
it means that neither $\delta^{\rm o}$ nor $\delta^{\rm w}$ is present.
For a function of one variable defined on $(0,\infty)$,
the weighted norm $\|f\|^{(\g_1;0)}_{k,\ga;(l);\R^+}$
is understood in the same as the definition
above with the weight to $\{0\}$ and the decay
at infinity.

In the study of Problem \ref{pr-3.2} for a transonic solution
$(U^-_0, U^+_0)$ with wedge angle $\theta_{\rm w}$,
the variables in $U$
are expected to have different levels of regularity.
Thus, we distinguish these variables by defining
\begin{equation}\label{var-u1-w-p-rho}
{U_1}=(\uu\cdot \ttau_{\rm w}^0,\rho), \qquad U_2=(w, p)\,\,\,\,\mbox{with $w=\frac{\uu\cdot \nnu_{\rm w}^0}{\uu\cdot \ttau_{\rm w}^0}$},
\end{equation}
where $\nnu_{\rm w}^0=(-\sin\theta_{\rm w}, \cos \theta_{\rm w})$ and
$\ttau_{\rm w}^0=(\cos\theta_{\rm w}, \sin \theta_{\rm w})$.
We note that, for the solutions under our consideration,
the denominator in the definition of $w$ is strictly positive, since
it is a positive constant for the background solution.

Note that $U_{10}^+=(|\uu_0^+|, \rho_0^+)$ and $U_{20}^+=(0,p_0^+)$ are
the corresponding  quantities for the background subsonic state.
Moreover, $\nnu_{\rm w}^0$ is the interior (for $\Omega_0$) unit normal to $\del W_0$,
and $\ttau_{\rm w}^0$ is the tangential  unit vector to $\del W_0$, where $\del W_0$ and
$\Omega_0$ are defined by \eqref{wall} and  \eqref{domain:1} for the background solution
$(U^-_0, U^+_0)$, {\it i.e.},
$\uu\cdot \ttau_{\rm w}^0$ and $\uu\cdot \nnu_{\rm w}^0$ are the components $u_1$ and $u_2$ of $\uu$
in the coordinates rotated clockwise
by angle $\theta_{\rm w}$, so that the background downstream flow becomes horizontal.

\begin{theorem}[Chen-Chen-Feldman \cite{CCF3}] \label{thm-main} \quad
Let $(U^-_0, U^+_0)$ be a constant transonic solution
 for the wedge angle $\theta_{\rm w}\in (0, \theta_{\rm w}^{\rm d})$.
There are positive constants $\ga , \b,  C_0$, and $\varepsilon$, depending only on the background
states $(U^-_0, U^+_0)$,  such that

\begin{enumerate}[{\rm (i)}]
\item \label{thm-main-i1}
If $(U^-_0, U^+_0)$ corresponds to the state on arc $\wideparen{TS}$,
and
\begin{equation}\label{cond-U-small-pert}
\| U^- -U^-_0 \|_{2,\ga;(1+\b,0);\Omega^-} +\| b'-\tan\theta_{\rm w}\|^{(-\ga;0)}_{1,\ga; (1+\b);\R^+} <\ve,
\end{equation}
then
there exist a solution $(U,\s)$ of {\rm {Problem \ref{pr-3.2}(WT)}}
and a function $U^\infty(\xx)=(u^\infty_1, 0, p_0^+, \rho^\infty)(\xx)=Z^\infty( -x_1\sin\theta_{\rm w}+x_2\cos\theta_{\rm w})$
with $Z^\infty:[0, \infty)\to \R^4$ of form  $Z^\infty=(z_1, 0, p_0^+, z_4)$
such that $U_1$ and $U_2$ defined by \eqref{var-u1-w-p-rho} satisfy
\begin{align}
&\qquad
\|U_1 - U_1^\infty\|^{(-\ga;\del  W)}_{2,\ga;(\gb,1);\Omega}
+\|U_2 - U_{20}^+\|^{(-\ga;  \O)(-1-\ga;\del  W)}_{2,\ga;(1+\gb,0);\Omega}
+ \| \s' - s_0\|^{(-\ga;0)}_{2,\ga;(1+\gb);\R^+}
+  \| Z_1^\infty - U^+_{10} \|^{(-\ga;0)}_{2,\ga;(1+\gb);[0, \infty)}
\nonumber\\[2mm]
&\qquad \le{}
  C_0 \left(\| U^- -U^-_0 \|_{2,\ga;(1+\b,0);\Omega^-}
     +\|b'-\tan\theta_{\rm w}\|^{(-\ga;0)}_{1,\ga; (1+\b);\R^+} \right),
\label{est-U-small-pert}
\end{align}
where we have denoted $U^\infty_1:=(\uu^\infty\cdot\ttau_{\rm w}^0, \rho^\infty)=
(u^\infty_1\cos\theta_{\rm w}, \rho^\infty)$
and $Z^\infty_1:=(z_1\cos\theta_{\rm w}, z_4)${\rm ;}

\smallskip
\item
\label{thm-main-i2}
If $(U^-_0, U^+_0)$ corresponds to the state on arc $\wideparen{TH}$
and
\begin{equation}\label{cond-U-small-pert2}
\|U^- -U^-_0\|_{2,\ga;(\b,0);\Omega^-} + \|b'-\tan\theta_{\rm w}\|^{(-\ga -1;0)}_{2,\ga; (\b);\R^+}<\ve,
\end{equation}
then
there exists a solution $(U,\s)$ of
 {\rm {Problem \ref{pr-3.2}(ST)}} such that $U_1$ and $U_2$ defined by \eqref{var-u1-w-p-rho} satisfy
\begin{equation}\label{est-U-small-pert2}
\begin{split}
&\|U_1 - U_{10}^+\|^{(-1-\ga;\pw)}_{2,\ga;(0,\gb);\Omega}
	+\|U_2 - U_{20}^+\|^{(-1-\ga; \O)}_{2,\ga;(\gb,0);\Omega}
+\|\s' - s_0\|^{(-1-\ga;0)}_{2,\ga;(\gb);\R^+}
\\[2mm]
&\le{} C_0  \left(\| U^- -U^-_0\|_{2,\ga;(\b);\Omega^-}
  +\|b'-\tan\theta_{\rm w}\|^{(-1-\ga;0)}_{2,\ga; (\b);\R^+}  \right).
\end{split}
\end{equation}
\end{enumerate}
The solution $(U,\s)$ is unique within the class of solutions
for each of {\rm {Problem \ref{pr-3.2}(WT)}} and {\rm {Problem \ref{pr-3.2}(ST)}}
when the left-hand sides of
\eqref{est-U-small-pert} for {\rm {Problem \ref{pr-3.2}(WT)}}
and \eqref{est-U-small-pert2}
for {\rm {Problem \ref{pr-3.2}(ST)}} are less than $C_0 \vae$ correspondingly.
\end{theorem}

The dependence of constants $\ga , \b,  C_0$, and $\varepsilon$ in Theorem \ref{thm-main}
is as follows:  $\ga$ and $\b$ depend on $(U^-_0,U^+_0)$
but are independent of $(C_0,\varepsilon)$,
$C_0$ depends on $(U^-_0, U^+_0,\ga, \beta)$ but is independent of $\ve$,
and $\ve$ depends on all $(U^-_0, U^+_0, \ga , \b, C_0)$.

The difference in the results of the two problems is that the solution of
{\rm {Problem \ref{pr-3.2}(WT)}} has less regularity at corner ${O}$
and decays faster with respect to $|\xx|$ (or the distance from the wedge boundary)
than the solution of {\rm {Problem \ref{pr-3.2}(ST)}}.

Notice that part \eqref{thm-main-i1} of Theorem \ref{thm-main} gives the asymptotics
of solution $U$ as $|\xx|\to \infty$ within
$\Omega$, and $U^\infty$ is an asymptotic profile.
Moreover, the convergence of $U_2$
to $U_2^\infty=U_{20}^+$ as $|\xx|\to \infty$
is of polynomial rate $|\xx|^{-(\beta+1)}$ that
is faster than the convergence rate of $U_1$, which is $|\xx|^{-\beta}$.
However, as $x_2\to \infty$, both $U_1$ and $U_2$
decay to $U_{10}^+$ and $U_{20}^+$, respectively,
with the decay rate $x_2^{-(\beta+1)}$,
which can be seen by combining the estimates of
the first and last terms on the right-hand side of
\eqref{est-U-small-pert} for $U_1$.
Part \eqref{thm-main-i2} of Theorem \ref{thm-main} does not give the asymptotic limit
of $U_1$ as $|\xx|\to \infty$, while $U_2$ converges
to $U_{20}^+$ with the decay rate $|\xx|^{-\beta}$. Also,
as $x_2\to\infty$, both $U_1$ and $U_2$
decay to $U_{10}^+$ and $U_{20}^+$, respectively,
with the decay rate $x_2^{-\beta}$ for part \eqref{thm-main-i2}.

Furthermore, for both parts \eqref{thm-main-i1} and \eqref{thm-main-i2} of Theorem \ref{thm-main},
the asymptotic profile in the Lagrangian coordinates is given in Theorem \ref{LagrCoordAppr2-thm}.

\subsection{Approach I for {Problem \ref{pr-3.2}(WT)}}
\label{Sect3.2}
We now describe Approach I for solving {\rm {Problem \ref{pr-3.2}(WT)}}.
We work in the Lagrangian coordinates introduced in \eqref{def-coord}.
From the slip condition \eqref{slipcon} on the wedge boundary $\del {W}$, it follows that
$\del {W}$  is a streamline so that
$\del {W}$ becomes the half-line
$\mathcal{L}_1=\{(y_1, y_2)\,:\;y_1>0, y_2=0\}$
in the Lagrangian coordinates.
Let $\sS=\{y_2= \hat{\sigma}(y_1)\}$ be a shock-front. Then,
from equations \eqref{eqn-euler1}--\eqref{eqn-euler4}, we can derive the Rankine-Hugoniot
conditions along $\sS$:
\begin{eqnarray}\label{con-RH1}
 &&\hat{\sigma}'(y_1)\big[\frac{1}{\rho u_1}\big]=-\big[\frac{u_2}{u_1} \big],\\[0.5mm]
 &&\hat{\sigma}'(y_1) \big[ u_1 + \frac{p}{\rho u_1}\big]=-\big[\frac{p u_2}{u_1}\big],
  \label{con-RH2}\\[0.5mm]
&&\hat{\sigma}'(y_1)[\,u_2 \,]= [\,p\,],
\label{con-RH3}\\[0.5mm]
&&\big[\frac{1}{2}|\uu|^2 + \frac{\gc p}{(\gc-1)\rho}\big]= 0.
\label{con-RH4}
\end{eqnarray}
The background shock-front in the Lagrangian coordinates is $\sS_0=\{y_2 = s_1 y_1\}$ with
$s_1 = \rho_0^+ u_{10}^+(s_0- \tan \theta_0)>0$.

Without loss of generality, we assume that, in the Lagrangian coordinates, the supersonic solution
$U^-$ exists in domain $\D^-$ defined by
\begin{equation}\label{def-super-D}
    \D^- = \left\{\yy\,:\, y_2>\frac{s_1}{2}y_1,\, y_1>0 \right\}.
\end{equation}
For a given shock function $\hat{\sigma}(y_1)$, let
\begin{eqnarray}
\D^-_{\hat{\sigma}}
&=&\big\{\yy: y_2>\hat{\sigma}(y_1),\, y_1>0\big\},
\label{3.11a}\\[1mm]
\D_{\hat{\sigma}}
&=& \big\{\yy: 0 <y_2<\hat{\sigma}(y_1),\, y_1>0\big\}. \label{3.12a}
\end{eqnarray}
Then {\bf Approach I} consists of three steps:

\medskip
1. {\bf Potential function $\phi(\yy)$}.
We first use a potential function to reduce the full Euler equations
\eqref{eqn-euler1}--\eqref{eqn-euler4}
to a scalar second-order nonlinear elliptic PDE in the subsonic region.
This method was first proposed in \cite{CCF} in which
the advantage of the conservation properties of the Euler system is taken for the reduction.

More precisely, since $\rho u_1 \ne 0$ in either the supersonic or subsonic region,
it follows from \eqref{eqn-euler1} that
there exists a potential function of the vector field
$(\frac{u_2}{u_1}, \frac{1}{\rho u_1})$ such that
\begin{equation}\label{eqn-phi-deri}
    D\phi= (\frac{u_2}{u_1}, \frac{1}{\rho u_1}) \qquad\,\, \mbox{with $\phi(\mathbf{0})=0$}.
\end{equation}
Equation \eqref{eqn-euler4} implies the Bernoulli law:
\begin{equation}\label{eqn-bernoulli}
\frac{1}{2}q^2 + \frac{\gc p}{(\gc-1)\rho} = B(y_2),
\end{equation}
where $q=|\uu|=\sqrt{u_1^2+u_2^2}$, and
$B=B(y_2)$ is completely determined by the given incoming flow $U^-$
at the initial position $\mathcal{I}$ because of the
Rankine-Hugoniot condition \eqref{con-RH4}.

From equations \eqref{eqn-euler1}--\eqref{eqn-euler4}, we find
\begin{equation}\label{3.30a}
\Big(\frac{p}{\rho^\gamma}\Big)_{y_1}=0,
\end{equation}
which implies
\begin{equation}\label{eqn-rho-p}
    p= A(y_2) \rho^\gc \qquad\,\,\mbox{in the subsonic region $\D_{\hat{\sigma}}$}.
\end{equation}

\smallskip
With equations \eqref{eqn-phi-deri} and \eqref{eqn-rho-p}, we can
rewrite the Bernoulli law \eqref{eqn-bernoulli} as
\begin{equation}\label{eqn-rho}
\frac{\phi_{y_1}^2 + 1 }{2 \phi_{y_2}^2} + \frac{\gc}{\gc
-1}A\rho^{\gc+1} = B \rho^2.
\end{equation}
In the subsonic region, $q=|\uu|< c=\sqrt{\frac{\gc p}{\rho}}$. Therefore,
the Bernoulli law \eqref{eqn-bernoulli} implies
\begin{equation}\label{con-subsonic}
    \rho^{\gc-1} > \frac{2(\gamma-1)B}{\gamma(\gamma+1)A}.
\end{equation}
Condition \eqref{con-subsonic} guarantees that  $\rho$ can be
solved from \eqref{eqn-rho} as a smooth function of $(D\phi,A,B)$.

\medskip
Assume that $A=A(y_2)$ has been determined. Then $(\uu, p, \rho)$ can be
expressed as functions of $D\phi$:
\begin{equation}\label{eqn-express-U}
\rho = \rho(D\phi,A, B),\qquad
\uu=(\frac{1}{\rho \phi_{y_2}}, \frac{\phi_{y_1}}{\rho \phi_{y_2}}),\qquad  p=A\rho^\gc,
\end{equation}
since $B=B(y_2)$ is given by the incoming flow.

Similarly, in the supersonic region $\D^-$, we employ the corresponding
variables $(\phi^-, A^-, B)$ to replace $U^-$, where $B$ is the
same as in the subsonic region because of the Rankine-Hugoniot
condition \eqref{con-RH4}.

\smallskip
We now choose \eqref{eqn-euler3} to derive a {\it second-order
nonlinear elliptic equation for $\phi$} so that the full Euler
system \eqref{eqn-euler1}--\eqref{eqn-euler4}
is reduced to the following nonlinear PDE in the subsonic region:
\begin{equation}\label{eqn-N}
  \big( N^1(D\phi,A,B)\big)_{y_1} + \big(N^2(D\phi,A,B)\big)_{y_2}=0,
\end{equation}
where $(N^1,N^2)(D\phi,A,B)=(u_2, p)(D\phi,A,B)$
are given by
\begin{equation}\label{3.36a}
N^1(D\phi,A,B)=
\frac{\phi_{y_1}}{\phi_{y_2}\rho(D\phi,A(y_2),B(y_2))},\quad
N^2(D\phi,A,B)=A(y_2)\big(\rho(D\phi,A(y_2),B(y_2))\big)^\gamma.
\end{equation}
Then a careful calculation shows
that
\begin{equation}
N^1_{\phi_{y_1}} N^2_{\phi_{y_2}} -
N^1_{\phi_{y_2}}N^2_{\phi_{y_1}}= \frac{c^2\rho^2 u_1^2}{c^2 - q^2}
>0
\end{equation}
in the subsonic region with $\rho u_1\ne 0$.
Therefore, when $\phi$ is sufficiently close to $\phi_0^+$ (determined
by the subsonic background state $U_0^+$) in the $C^1$--norm,
equation \eqref{eqn-N} is uniformly elliptic, and the Euler system
\eqref{eqn-euler1}--\eqref{eqn-euler4} is reduced to the elliptic
equation \eqref{eqn-N} in domain $\D_{\hat{\sigma}}$,
where $\hat{\sigma}$ is the function for the free boundary (transonic shock).

The boundary condition for $\phi$ on the wedge boundary $\{y_2 =0\}$ is
derived from the fact that $\phi(y_1, y_2)=x_2(y_1, y_2)$ by
\eqref{psi-potential}--\eqref{def-coord} and \eqref{eqn-phi-deri}.
Then, recalling that $\del {W}=\{\xx\,:\;x_2=b(x_1), x_1>0\}$ in the $\xx$--coordinates
which is
$\{ \yy\,:\;y_2=0, y_1>0\}$  in the $\yy$--coordinates
and using $y_1=x_1$ by \eqref{def-coord}, we obtain
\begin{equation}\label{Con-bd-phi-wall}
\phi (y_1, 0) = b(y_1).
\end{equation}
The condition on $\sS$ is derived from the Rankine-Hugoniot
conditions \eqref{con-RH1}--\eqref{con-RH3}. Condition
\eqref{con-RH1} is equivalent to the continuity of $\phi$ across
$\sS$:
\begin{equation}\label{continuity-b}
[\phi]|_{\sS}=0,
\end{equation}
which, by \eqref{eqn-phi-deri}, gives
\begin{equation}\label{for-tau-prime}
    \hat{\sigma}'(y_1)=-\frac{[\phi_{y_1}]}{[\phi_{y_2}]}(y_1,\hat{\sigma}(y_1)).
\end{equation}
Replacing $\hat{\sigma}'(y_1)$ in \eqref{con-RH2}--\eqref{con-RH3}
by \eqref{for-tau-prime} and using \eqref{eqn-rho-p} give rise to the
conditions on $\sS$:
\begin{eqnarray}
&&G(D\phi, A, U^-)\equiv
  [\phi_{y_1}]\big[\frac{1}{\rho \phi_{y_2}}+A\rho^\gc \phi_{y_2}\big]
   -[\phi_{y_2}][A\rho^\gc \phi_{y_1}] =0,  \label{eqn-G}\\[1mm]
&&H(D\phi, A, U^-)\equiv
[\phi_{y_1}][N^1]+[\phi_{y_2}][N^2]=0.\label{eqn-H}
\end{eqnarray}
We now combine the above two conditions into the boundary
condition for \eqref{eqn-N} by eliminating $A$.
Taking the partial derivative of $G$ and $H$ with respect
to $A$ respectively and making careful calculation, we have
\begin{align*}
G_A &= [\phi_{y_1}]
  \Big(\frac{N^1_A}{\phi_{y_1}}
   + \phi_{y_2} N^2_A \Big)-[\phi_{y_2}]\phi_{y_1} N^2_A \\
  &= \frac{u_2 \rho^\gc(q^2 + \frac{c^2}{\gc-1})}{u_1(c^2 -q^2)}
   \left[\frac{1}{\rho u_1}\right]
   - \frac{ \rho^{\gc -1}}{u_1(c^2-q^2)}
   \left(u_2^2 + \frac{c^2 -u_1^2}{\gc -1}\right)\left[\frac{u_2}{u_1} \right]<0,
\end{align*}
and
\begin{align*}
  H_A = N^1_A [\phi_{y_1}] + N^2_A[\phi_{y_2}]
   = \frac{\gc}{\gc -1}\frac{\rho^{\gc-1}u_2}{c^2 -q^2}
   \left[\frac{u_2}{u_1}\right]
   - \frac{\rho^\gc(q^2 + \frac{c^2}{\gc-1})}
   {c^2 -q^2}\left[\frac{1}{\rho u_1}\right]
   > 0,
\end{align*}
since $[\frac{1}{\rho u_1}]<0$ and $u_2^-$ is close to $0$.
Therefore, both equations \eqref{eqn-G} and \eqref{eqn-H} can be solved
for $A$ to obtain
$A= g_1(D\phi, U^-)$ and $A= g_2(D\phi,U^-)$, respectively.
With these, we obtain our desired condition
on the free boundary ({\it i.e.}, the shock-front):
\begin{equation}\label{con-phi-shock}
\bar{g}(D\phi,U^-):= (g_2 -g_1)(D\phi,U^-) =0.
\end{equation}
Then the original free boundary problem, Problem \ref{pr-3.2}, is reduced to
the following free boundary problem for the elliptic equation \eqref{eqn-N}:

\begin{problem}[Free Boundary Problem]\label{problem-3.3}
Seek $(\hat\sigma, \phi, A)$  such that $\phi$ is a  solution of the elliptic
equation \eqref{eqn-N} in the region with the fixed boundary condition
\eqref{Con-bd-phi-wall} and
conditions \eqref{continuity-b} and  \eqref{eqn-G}--\eqref{con-phi-shock}
on $\sS$.
\end{problem}

\medskip
2. {\bf Hodograph transformation and fixed boundary value problem}.
In order to solve the free boundary problem, we employ the
hodograph transformation so that the shock-front becomes a fixed boundary
by using the free boundary conditions \eqref{continuity-b} and \eqref{con-phi-shock}.
This allows us to find $\phi$ for each $A$ from an appropriately chosen set.
After that, we only need to perform an iteration for the unknown function $A$
to satisfy \eqref{eqn-G}--\eqref{eqn-H}.

\smallskip
Note that the expected solutions in Theorem \ref{thm-main} satisfy that
$\|U-U^+_0\|_{L^\infty(\Omega)}\le C_0\varepsilon$. Then, denoting by
$\phi_0^+$ the potential function \eqref{eqn-phi-deri} for the subsonic background  state $U^+_0$,
we obtain that $\phi$ is close to $\phi_0^+$ in $C^1$ on the closure of the subsonic region.
On the iteration, we consider (and eventually obtain) solutions $U$ for which the same property holds.
Thus, we assume that $\phi$ is close to $\phi_0^+$ in $C^1(\overline{\D_{\hat{\sigma}}})$ below;
see \eqref{3.12a}.

We now extend the domain of $\phi^-$ from $\D^-$ to the first quadrant $\D^- \cup \D_{\hat{\sigma}}$.
Let
$$
\phi^-_0 = \frac{1}{\rho^-_0 u_{20}^-} y_2,
$$
which is the potential function \eqref{eqn-phi-deri} for the supersonic background state $U^-_0$.
Then $\phi^-$ is close to $\phi^-_0$ in $C^1(\overline{\D^-})$ since $U^-$ is close to $U^-_0$ in
$L^\infty$ (and in the stronger norm; see Theorem \ref{thm-main}).  We can extend $\phi^-$ into
$\D^- \cup \D_{\hat{\sigma}}$
so that it remains close to $\phi^-_0$ in $C^1$ on the closure of $\D^- \cup \D_{\hat{\sigma}}$.
We then use the following partial hodograph transformation:
\begin{equation}\label{def-hodo}
(y_1,y_2) \to (z_1,z_2)=(\phi - \phi^-, y_2).
\end{equation}
Note that
$\partial_{y_1}(\phi_0^+-\phi_0^-)=\frac{u^+_{20}}{u^+_{10}}>0$ by using \eqref{eqn-phi-deri},
where $(u^+_{10}, u^+_{20})$ is the velocity of the background subsonic state $U_0^+$
and the fact that $u^-_{20}=0$ has been used.
Since $\phi$ and $\phi^-$ are  close  to $\phi_0^+$ and  $\phi_0^-$ in the $C^1$--norm respectively,
transformation \eqref{def-hodo} is invertible, so that
$y_1$ is a function of $\zz:=(z_1, z_2)$, denoted as $y_1 = \varphi(\zz)$.

Let
\begin{eqnarray*}
  M^1 (D\phi, A, U^-)= N^1(D\phi,A,B)
  + N^2(D\phi, A, B)\frac{[\phi_{y_2}]}{[\phi_{y_1}]}, \quad
  M^2 (D\phi, A, U^-) =\frac{N^2(D\phi, A,B)}{[\phi_{y_1}]},
\end{eqnarray*}
and
\begin{eqnarray*}\label{def-Mbar}
\overline{M}^i ( D\varphi,\varphi, A, \zz)
=-M^i( \partial_{y_1}{\phi^-}(\varphi, z_2)
     +\frac{1}{\varphi_{z_1}},
\,\partial_{y_2}{\phi^-}(\varphi, z_2)
     -\frac{\varphi_{z_2}}{\varphi_{z_1}}, A, U^-(\varphi, z_2)), \qquad i=1,2.
\end{eqnarray*}
Then equation \eqref{eqn-N} becomes
\begin{equation}\label{eqn-Mbar}
   \big(\overline{M}^1 (D\varphi, \varphi, A, \zz)\big)_{z_1}
   +\big(\overline{M}^2 (D\varphi, \varphi, A, \zz)\big)_{z_2}=0.
\end{equation}
Notice that
\begin{eqnarray*}
\overline{M}^1_{\vph_{z_1}} \overline{M}^2_{\vph_{z_2}}  -
\frac{1}{4}\big(\overline{M}^1_{\vph_{z_2}}+ \overline{M}^2_{\vph_{z_1}}\big)^2=
[\phi_{y_1}]^2 \big(N^1_{\phi_{y_1}} N^2_{\phi_{y_2}} -
(N^1_{\phi_{y_2}})^2\big) >0,
\end{eqnarray*}
which implies that equation \eqref{eqn-Mbar} is uniformly elliptic, for any solution
$\varphi$ that is close to $\varphi_0^+$ (determined by \eqref{def-hodo} with $\phi=\phi_0^+$)
in the $C^1$--norm.

Under transform \eqref{def-hodo}, the unknown shock-front $\sS$ becomes a fixed boundary, which
is the $z_2$-axis,
where we have used that $\phi$ is close in $C^1$ to $\phi_0^+$
in $\overline{\D_{\hat{\sigma}}}$
and to $\phi_0^-$ in $\overline{\D_{\hat{\sigma}}^-}$ in order
to conclude that
$\phi$ is Lipschitz across $\sS$ from \eqref{eqn-phi-deri} and then that $\phi=\phi^-$
on $\sS$ but $\phi\ne\phi^-$ in $\overline{\D_{\hat{\sigma}}}\setminus\sS$.
Along the $z_2$-axis, condition
\eqref{con-phi-shock} is now
\begin{align}\label{con-phi-gt}
\tilde{g} (D\varphi, \varphi, \zz)&:=
     \bar{g}(\partial_{y_1}{\phi^-}(\varphi, z_2)
     +\frac{1}{\varphi_{z_1}}, \partial_{y_2}{\phi^-}(\varphi, z_2)
     -\frac{\varphi_{z_2}}{\varphi_{z_1}}, U^-(\varphi,z_2))\nonumber\\
     &= 0      \qquad\quad\mbox{on } \{z_1=0,\,z_2>0\}.
 \end{align}

We also convert condition \eqref{eqn-H} into the $\zz$--coordinates:
Along the $z_2$-axis,
\begin{eqnarray}\label{con-Ht}
\widetilde{H}(D\varphi,\varphi, A, \zz)
:= {H}(\partial_{y_1}{\phi^-}(\varphi, z_2)
     +\frac{1}{\varphi_{z_1}},\,\partial_{y_2}{\phi^-}(\varphi, z_2)
     -\frac{\varphi_{z_2}}{\varphi_{z_1}}, A,  U^-(\varphi,z_2))=0.
\end{eqnarray}

\medskip
The condition on the $z_1$-axis can be derived from
\eqref{Con-bd-phi-wall} as follows: Restricted on $z_2 =0$, the
coordinate transformation \eqref{def-hodo} becomes
\begin{equation*}
    z_1 = b(y_1) - \phi_-(y_1, 0).
\end{equation*}
Then $y_1$ can be expressed in terms of $z_1$ as $y_1  = \widetilde{b}(z_1)$ so that
$\varphi(z_1, 0)=y_1$ becomes
\begin{equation}\label{con-phi-z1}
\varphi(z_1, 0) = \widetilde{b}(z_1)  \qquad\mbox{ on $\mathcal{L}_1:=\{z_2=0,\,z_1>0\}$}.
\end{equation}

Therefore, the original wedge problem has now been reduced to
the following problem on the first quadrant $\mathbb{Q}$.

\begin{problem}[Fixed Boundary Value Problem]\label{problem-3.4}
Seek  $(\varphi,A)$ such that $\varphi$ is a solution of the second-order nonlinear elliptic equation \eqref{eqn-Mbar}
in the unbounded domain $\mathbb{Q}$ with the boundary conditions \eqref{con-phi-gt}
and \eqref{con-phi-z1}, and such that \eqref{con-Ht} holds.
\end{problem}

\medskip
3. {\bf Solution to the fixed boundary value problem -- Problem \ref{problem-3.4}}.
Through the shock polar,
we can determine the values of $U$ at the origin so that $A(0)$ is fixed, depending
on the values of $U^-(\mathbf{0})$ and $b'(0)$.
Then we solve
\eqref{con-Ht} to obtain a unique solution
$\tilde A=h(\zz,\phi, D\phi)$ that defines the iteration map.

This is achieved by the following fixed point argument.
Consider a Banach space:
$$
X= \{\mathcal{A}\,: \, \mathcal{A}(0) =0,\, \|\mathcal{A}\|_{1,\alpha; (1+\beta); \mathbb{R}^+}^{(-\alpha); \{0\}}<\infty\}.
$$
Then we define our iteration map $\mathcal{J}:  X\longrightarrow X$
through the following:

\smallskip
First, we define a smooth cutoff
function $\chi$ on $[0,\infty)$ such that
\[
\chi(s)=
\begin{cases}
    1 &\,\, \mbox{for $0\le s <1$}, \\[0.5mm]
    0 &\,\, \mbox{for $s>2$}. \\
\end{cases}
\]
 Set
\begin{equation}\label{def-fix-a0}
A(0) := t(\omega(0), b'(0))
\qquad \mbox{for $\omega= U^- - U^-_0$},
\end{equation}
where $t$ is a function
determined by the Rankine-Hugoniot conditions
\eqref{eqn-G}--\eqref{eqn-H}.
Then we define $w_t(z_2)$ as
\begin{equation}\label{6.2-a}
w_t(z_2):= A^+_0 + \big( t(\omega(0),
b'(0))-A^+_0\big)\chi(z_2),
\end{equation}
where $A^+_0 = \frac{p_0^+}{(\rho_0^+)^\gc}$.

\smallskip
Consider any $A=A(z_2)$ so that $A-w_t\in X$ satisfying
\begin{equation}\label{4.11-a}
\|A-A_0^+\|_{1,\alpha; (1+\beta); \mathbb{R}^+}^{(-\alpha); \{0\}}\le C_0\varepsilon
\end{equation}
for some fixed constant $C_0>0$.
With this $A$, we solve equation \eqref{eqn-Mbar} for $\varphi=\varphi_A$ in
the unbounded domain $\mathbb{Q}$ with the boundary conditions \eqref{con-phi-gt}
and \eqref{con-phi-z1}, and with the asymptotic condition $\varphi\to \varphi^\infty$ as $\xx\to\infty$,
where the limit is understood in the appropriate sense, $\varphi^\infty$ is the solution of
\begin{equation}\label{eqn-varphi-inf}
z_1=(\phi^\infty-\phi^-)(\varphi^\infty,z_2),
\end{equation}
with $\phi^\infty=\frac{u_{20}^+}{u_{10}^+}y_1+l(y_2)$,
and
$l(y_2)$ is determined by the Bernoulli law \eqref{eqn-rho},
via replacing
$\phi$ and $\rho$ by their asymptotic values
$\phi^\infty$ and $\rho^\infty(y_2)=\big(\frac{p_0^+}{A(y_2)}\big)^{\frac 1\gamma}$
and noting that $B=B(y_2)$ is determined by the upstream state $U^-$.
More specifically, we show the existence of a solution  $\varphi$
of \eqref{eqn-Mbar}--\eqref{con-phi-gt} and \eqref{con-phi-z1}
in the set:
\begin{equation*}
\Sigma_\delta=\big\{\varphi :\, \|\varphi-\varphi^\infty\|_{2,\alpha;(\beta,0);\mathbb{Q}}^{(-1-\alpha);\mathcal{L}_1}\le \delta\big\}
\qquad\mbox{for sufficiently small $\delta>0$},
\end{equation*}
which is a compact and convex subset of the Banach space:
\begin{equation*}
\mathcal{B} = \big\{\varphi \,:\, \|\varphi-\varphi^\infty\|_{2,\ga';(\gb',0);\mathbb{Q}}^{(-1-\ga');\mathcal{L}_1} < \infty\big\}
\qquad\,\, \mbox{with $0< \ga'<\ga$ and $0<\beta'<\beta$}.
\end{equation*}
For $\varphi  \in \Sigma_\gd$, equation \eqref{eqn-Mbar}
is uniformly elliptic if $\delta>0$ is small.
This allows us to solve the problem for $\varphi=\varphi_A \in \Sigma_\gd$
by the Schauder fixed point theorem
if the perturbation is small, {\it i.e.},
if $\varepsilon$ is small in \eqref{4.11-a} and the conditions of Theorem \ref{thm-main}.
Then, with this $\varphi=\varphi_A$, we solve \eqref{con-Ht} to obtain
a unique $\tilde{A}$ that defines the iteration map $J$ by $\mathcal{J}(A-w_t):=\tilde{A}-w_t$.

\smallskip
Finally, by the implicit function theorem,
we prove that $\mathcal{J}$ has a fixed point $A-w_t$, for which $A$ satisfies
\eqref{4.11-a}.

\smallskip
For more details for this approach, see Chen-Chen-Feldman \cite{CCF,Chen-Chen-Feldman}.
This approach can also be applied to Problem \ref{pr-3.2}(ST);
see \cite{YinZhou} for the case when the wedge angle is sufficiently small.

\smallskip
\subsection{Approach II for Problem \ref{pr-3.2}(ST) and Problem \ref{pr-3.2}(WT)}
\label{Sect3.3}

\newcommand{\brt}{b_{\rm rot}}
We now describe the second approach, Approach II.
It allows us to handle both cases in
Theorem \ref{thm-main}: {\rm {Problem \ref{pr-3.2}(WT)}} and {\rm {Problem \ref{pr-3.2}(ST)}}.
In particular, for {\rm {Problem \ref{pr-3.2}(WT)}},
this approach yields a better asymptotic decay rate, as stated in \eqref{est-U-small-pert}.

It is convenient to rotate the $\xx$--coordinates clockwise
by the wedge angle $\theta_{\rm w}$, so that the background downstream flow becomes horizontal, as
discussed in the paragraph before Theorem \ref{thm-main}.
We still use the same notations in the rotated coordinates when no confusion arises;
in particular, we write $\xx=(x_1, x_2)$ and $\uu=(u_1, u_2)$ in the rotated basis.
Then, in the new coordinates,
\begin{equation}\label{rotatedCoor-Appr2}
\frac{u_{20}^-}{u_{10}^-}=-\tan\theta_{\rm w}, \quad
U_0^-=(u_{10}^-,\, -u_{10}^-\tan\theta_{\rm w},\,p_0^-, \,\rho_0^-),\quad
U_0^+=(u_{10}^+,\, 0,\,p_0^+, \,\rho_0^+).
\end{equation}
Since the velocity components $(u_1, u_2)$ are now in the basis
$(\ttau_{\rm w}^0, \nnu_{\rm w}^0)$,
{\it i.e.}, $u_1=\uu\cdot \ttau_{\rm w}^0$ and $u_2=\uu\cdot \nnu_{\rm w}^0$,
we see that, by \eqref{var-u1-w-p-rho},
\begin{equation}\label{var-u1-w-p-rho-rot}
{U_1}=(u_1,\rho), \qquad U_2=(w, p)\,\,\,\mbox{with $w=\frac{u_2}{u_1}$}
\end{equation}
in the new coordinates.
Furthermore, we obtain from \eqref{wall}
and \eqref{cond-U-small-pert} or \eqref{cond-U-small-pert2} with small $\varepsilon$
that, in the rotated coordinates,
\begin{equation}\label{wall-rot}
\partial{W}=\{\xx\in \R^2\,:\,  x_2= \brt(x_1),\  \brt(0)=0\},
\end{equation}
and function $\brt(x_1)$ satisfies the estimates in \eqref{cond-U-small-pert-rot}
or \eqref{cond-U-small-pert2-rot} below, respectively, with $C\varepsilon$
instead of $\varepsilon$
when $\varepsilon$ is small, where $C$ depends only on $b(\cdot)$.
For the background solution,
$b_{{\rm rot}, 0}=0$, {\it i.e.}, $\partial{W}_0$ is the positive $x_1$-axis.

We now construct a solution with a shock-front $\sS$ expressed as
\eqref{shock:1} in the rotated coordinates with a
function
$\tilde\sigma(x_1)$. The background shock
is now expressed as $\mathcal{S}_0:=\{x_2=\tilde\sigma_0(x_1)\,:\, x_1>0\}$ for
$\tilde{\sigma}_0(x_1):=\tilde s_0x_1$, where $\tilde s_0=\tan(\arctan s_0-\theta_{\rm w})$.
Then the subsonic region of the solution has the form:
\begin{equation}\label{domain:1-rot}
\Omega=\{\xx\in \R^2\,:\,  b_{\rm rot}(x_1)<x_2<\tilde \s(x_1), x_1>0\}.
\end{equation}
We can assume that the upstream steady supersonic smooth solution $U^-(\xx)$ exists
in region $\Omega^- = \{\xx\,:\, \frac{\tilde s_0}{2} x_1<x_2<2\tilde s_0 x_1, x_1\ge 0\}$,
beyond the background shock, but is still close to $U_0^-$.
Moreover, in part \eqref{thm-main-i1} of Theorem \ref{thm-main},
$U^\infty$ is independent of $x_1$ and $U^\infty=Z^\infty$ in the rotated
coordinates.

More specifically, we establish the following theorem in the rotated coordinates:
\begin{theorem}[Chen-Chen-Feldman \cite{CCF3}] \label{thm-main-rot} \quad
 Let $(U^-_0, U^+_0)$, given by \eqref{rotatedCoor-Appr2},
be a constant transonic solution
for the wedge angle $\theta_{\rm w}\in (0, \theta_{\rm w}^{\rm d})$.
There are positive constants $\ga , \b,  C_0$, and $\varepsilon$ depending only on the background
states $(U^-_0, U^+_0)$ such that

\begin{enumerate}[{\rm (i)}]
\item \label{thm-main-i1-b}
If $(U^-_0, U^+_0)$ corresponds to the state on arc $\wideparen{TS}$
and
\begin{equation}\label{cond-U-small-pert-rot}
\| U^- -U^-_0 \|_{2,\ga;(1+\b,0);\Omega^-} +\| \brt' \|^{(-\ga;0)}_{1,\ga; (1+\b);\R^+} <\ve,
\end{equation}
then there exist a solution $(U,\tilde \s)$ of {\rm {Problem \ref{pr-3.2}(WT)}}
and a function
$$
U^\infty(x_2)=(u^\infty_1(x_2), 0, p_0^+, \rho^\infty(x_2))
$$
so that $U_1$ and $U_2$ defined by \eqref{var-u1-w-p-rho} satisfy
\begin{align}
&\quad\,\,\,\, \|U_1 - U_1^\infty\|^{(-\ga;\del  W)}_{2,\ga;(\gb,1);\Omega}
+\|U_2 - U_{20}^+\|^{(-\ga;  \O)(-1-\ga;\del  W)}_{2,\ga;(1+\gb,0);\Omega}
+ \| \tilde \s' - \tilde s_0\|^{(-\ga;0)}_{2,\ga;(1+\gb);\R^+}
+ \| U_1^\infty - U^+_{10} \|^{(-\ga;0)}_{2,\ga;(1+\gb);[0, \infty)}\nonumber\\[2mm]
&\quad\,\,\,\, \le{}
  C_0 \left(\| U^- -U^-_0 \|_{2,\ga;(1+\b,0);\Omega^-}
     +\|\brt'\|^{(-\ga;0)}_{1,\ga; (1+\b);\R^+} \right),\label{est-U-small-pert-rot}
\end{align}
where $U^\infty_1=(u^\infty_1,  \rho^\infty)$.

\smallskip
\item
\label{thm-main-i2-b}
If $(U^-_0, U^+_0)$ corresponds to the state on arc $\wideparen{TH}$,
and
\begin{equation}\label{cond-U-small-pert2-rot}
\|U^- -U^-_0\|_{2,\ga;(\b,0);\Omega^-} + \|\brt'\|^{(-\ga -1;0)}_{2,\ga; (\b);\R^+}<\ve,
\end{equation}
then
there exists a solution $(U,\tilde \s)$ of
 {\rm {Problem \ref{pr-3.2}(ST)}} so that $U_1$ and $U_2$ defined by \eqref{var-u1-w-p-rho} satisfy
\begin{equation}\label{est-U-small-pert2-rot}
\begin{split}
&\|U_1 - U_{10}^+\|^{(-1-\ga;\pw)}_{2,\ga;(0,\gb);\Omega}
	+\|U_2 - U_{20}^+\|^{(-1-\ga; \O)}_{2,\ga;(\gb,0);\Omega}
+\|\tilde \s' - \tilde s_0\|^{(-1-\ga;0)}_{2,\ga;(\gb);\R^+}
\\
&\le{} C_0  \left(\| U^- -U^-_0\|_{2,\ga;(\b);\Omega^-}
  +\|\brt' \|^{(-1-\ga;0)}_{2,\ga; (\b);\R^+}  \right).
\end{split}
\end{equation}
\end{enumerate}
The solution $(U,\tilde \s)$ is unique within the class of solutions
for each of {\rm {Problem \ref{pr-3.2}(WT)}} and {\rm {Problem \ref{pr-3.2}(ST)}}
when the left-hand sides of
\eqref{est-U-small-pert} for {\rm {Problem \ref{pr-3.2}(WT)}} and \eqref{est-U-small-pert2}
for {\rm {Problem \ref{pr-3.2}(ST)}} are less than $C_0 \vae$ correspondingly.
\end{theorem}

Clearly, Theorem \ref{thm-main} follows from Theorem \ref{thm-main-rot} if $\varepsilon$ is small
so that, from the estimates of $\tilde \sigma$ in \eqref{est-U-small-pert-rot} or
\eqref{est-U-small-pert2-rot}, the shock remains a graph of $C^1$ function: $x_2=\s(x_1)$
after rotating the coordinates back.

\medskip
To prove Theorem \ref{thm-main-rot},
we work in the Lagrangian coordinates \eqref{def-coord} defined for the rotated coordinates $\xx=(x_1, x_2)$.
Then, as in the previous case, using the fact that the wedge boundary
$\del {W}$  is a streamline due to the slip condition \eqref{slipcon} on $\del {W}$,
we obtain that, in the present Lagrangian coordinates,
$\del {W}$ becomes the half-line:
$$
\mathcal{L}_1=\{(y_1, y_2)\,:\;y_1>0, y_2=0\}.
$$
The background shock-front $\sS_0$  is now given by
$\sS_0=\{y_2=  s_1y_1, y_1>0\}$ with $s_1=\rho_0^+u_{10}^+ \tilde s_0$.
We can assume that, in the Lagrangian coordinates, the supersonic solution
$U^-$ exists in domain $\D^-$ defined by \eqref{def-super-D}.
Shock $\sS$ is given by
$y_2= \hat\s(y_1)$ for $y_1>0$, where function $\hat \sigma$ differs from the one in Approach I
because the Lagrangian coordinates are now defined differently.
The supersonic region $\D^-_{\hat{\sigma}}$ and the subsonic region $\D_{\hat{\sigma}}$
of the solution are given by \eqref{3.11a} and \eqref{3.12a} respectively,
with the present function $\hat\sigma$.

We first present the existence and estimates of the solution in the Lagrangian coordinates:

\begin{theorem}\label{LagrCoordAppr2-thm}
Let $(U^-_0, U^+_0)$ be a constant transonic solution
for the wedge angle $\theta_{\rm w}\in (0, \theta_{\rm w}^{\rm d})$.
There are positive constants $\ga , \b,  C_0$, and $\varepsilon$ depending only on the background
states $(U^-_0, U^+_0)$  such that, if $\partial W$ in \eqref{wall-rot} and $U^-$ satisfy
\begin{enumerate}[{\rm (i)}]
  \item \eqref{cond-U-small-pert-rot}
for {\rm {Problem \ref{pr-3.2}(WT)}},

\smallskip
  \item \eqref{cond-U-small-pert2-rot}
for {\rm {Problem \ref{pr-3.2}(ST)}},
\end{enumerate}
then there exist a transonic shock $\sS_L=\{y_2=\hat\sigma(y_1),\,y_1>0\}$
and a subsonic solution $U=U(\yy)$ of \eqref{eqn-euler1}--\eqref{eqn-euler4} in
$\D_{\hat{\sigma}}$,
satisfying the Rankine-Hugoniot
conditions \eqref{con-RH1}--\eqref{con-RH4} along $\sS_L$ with $U^-$ expressed in the Lagrangian coordinates
in $\D^-_{\hat{\sigma}}$ and the slip condition $w_{|\mathcal{L}_1}=\brt'$,
as well as there exists a
function $\LU^\infty(y_2)=(u_1^\infty(y_2), 0, p_0^+, \rho^\infty(y_2))$,
such that $U(\yy)$  satisfies the following estimates:
\begin{enumerate}[{\rm (i)}]
\item For {\rm {Problem \ref{pr-3.2}(WT)}},
\begin{align}
&\qquad \|U_1 - \LU_1^\infty\|^{(-\ga;\mathcal{L}_1)}_{2,\ga;(1+\gb,0);\D_{\hat{\sigma}}}
+\|U_2 - U_{20}^+\|^{(-\ga;  \O)(-1-\ga;\mathcal{L}_1)}_{2,\ga;(1+\gb,0);\D_{\hat{\sigma}}}
+ \| \hat\s' - s_1\|^{(-\ga;0)}_{2,\ga;(1+\gb);\R^+}
 +  \| \LU_1^\infty - U^+_{10} \|^{(-\ga;0)}_{2,\ga;(1+\gb);\R^+}\nonumber\\[2mm]
&\qquad \le{}  C_0 \left(\| U^- -U^-_0 \|_{2,\ga;(1+\b,0);\D_{\hat{\sigma}}^-}
     +\|\brt'\|^{(-\ga;0)}_{1,\ga; (1+\b);\R^+} \right); \label{est-U-small-pert-Lagr}
\end{align}

\smallskip
\item For
 {\rm {Problem \ref{pr-3.2}(ST)}},
\begin{align}
&\qquad \|U_1 - \LU_1^\infty\|^{(-1-\ga;\mathcal{L}_1)}_{2,\ga;(\gb,0);\D_{\hat{\sigma}}}
	+\|U_2 - U_{20}^+\|^{(-1-\ga; \O)}_{2,\ga;(\gb,0);\D_{\hat{\sigma}}}
+\|\hat\s' - s_1\|^{(-1-\ga;0)}_{2,\ga;(\gb);\R^+}
 +  \| \LU_1^\infty - U^+_{10} \|^{(-1-\ga;0)}_{2,\ga;(\gb);\R^+}\nonumber\\[2mm]
&\qquad \le{} C_0  \left(\| U^- -U^-_0\|_{2,\ga;(\b);\D_{\hat{\sigma}}^-}
  +\|\brt'\|^{(-1-\ga;0)}_{2,\ga; (\b);\R^+}  \right), \label{est-U-small-pert2-Lagr}
\end{align}
where $\LU^\infty_1(y_2):=(u_1^\infty(y_2),  \rho^\infty(y_2))$.
\end{enumerate}
The solution $U$ is unique within the class of solutions
for each of {\rm {Problem \ref{pr-3.2}(WT)}}
and {\rm {Problem \ref{pr-3.2}(ST)}} when the left-hand sides of
\eqref{est-U-small-pert-Lagr} for {\rm {Problem \ref{pr-3.2}(WT)}}
and \eqref{est-U-small-pert2-Lagr} for {\rm {Problem \ref{pr-3.2}(ST)}}
are less than $C_0\varepsilon$ correspondingly.
\end{theorem}

We remark that function $\LU^\infty(y_2)$ can be understood as the asymptotic limit of $U(\yy)$ as $y_1\to\infty$.

\medskip
Now we describe the proof of Theorem \ref{LagrCoordAppr2-thm}, which is the main part of {\bf Approach II}.
Rewrite system  \eqref{eqn-euler1}--\eqref{eqn-euler4} into the
following nondivergence form for $U=(\uu, p,\rho)\in \R^4$:
\begin{equation}
 A(U) U_{y_1}^\top + B(U) U_{y_2}^\top =0, \label{euler-nondiv}
 \end{equation}
 where
\begin{align*}
A(U) =
\begin{bmatrix}
-\frac{1}{\rho u_1^2} & 0 &0& -\frac{1}{\rho^2 u_1}\\[1.5mm]
1 -  \frac{ p}{ \rho u_1^2} & 0 &  \frac{1}{\rho u_1} & -  \frac{ p}{ \rho^2 u_1}\\[1.5mm]
0  & 1& 0 & 0\\[1mm]
u_1 & u_2 & \frac{\g}{(\g - 1) \rho} & - \frac{\g p}{(\g - 1) \rho^2}
\end{bmatrix},\quad
B(U) =
 \begin{bmatrix}
\frac{u_2}{ u_1^2} & -\frac{1}{u_1} &0&0\\[1.5mm]
  \frac{ pu_2}{ u_1^2} & - \frac{p}{u_1} & - \frac{u_2}{ u_1} & 0 \\[1.5mm]
0  & 0& 1& 0\\[1mm]
0& 0 &0 & 0
\end{bmatrix}.
\end{align*}

Solving $\det (\gl A -B)=0$ for $\gl$, we obtain four eigenvalues:
 \begin{eqnarray*}
 \l_1 =\l_2 =0 \,\mbox{(real)},
 \qquad \l_{j}  =-\frac{c\rho}{c^2-u_1^2}\big(c u_2+(-1)^j u_1 \sqrt{c^2 -q^2} i\big)
   \,\, \mbox{for $j=3,4$ \, (complex)},
\end{eqnarray*}
where $q= \sqrt{u_1^2 + u_2^2} < c$ in the subsonic region and $i=\sqrt{-1}$.
 The corresponding left eigenvectors are
 \begin{align*}
&{\bf l}_1  = (0,0,0,1),\qquad {\bf l}_2  =(-p u_1, u_1, u_2, -1),\\
&{\bf l}_{3,4} =(\frac{p(\g p-\rho u_1^2)}{(\g -1) \rho u_1} \l_{3,4} + \frac{\g p^2 u_2}{(\g - 1) u_1},\,
-(u_1+\frac{\g p}{(\g -1) \rho u_1})\l_{3,4}-\frac{\g p u_2}{(\g -1)u_1},\, \frac{\g p}{\g -1} - u_2 \l_{3,4},\,  \l_{3,4}).
\end{align*}
Then
\begin{enumerate}
\item[(i)]
Multiplying equations \eqref{euler-nondiv} from the left by ${\bf l}_1$ leads to
the same equation \eqref{eqn-euler4}.
This, together with the Rankine-Hugoniot condition \eqref{con-RH4}, implies the Bernoulli law \eqref{eqn-bernoulli}
to be held in both supersonic and subsonic domains, as well as across the shock-front.
Therefore, $B(y_2)$ can be computed from the upstream flow
$U^-$.  If $u_1$ is a small perturbation of $u_{10}^+$, then $u_1>0$. Therefore,  we can solve \eqref{eqn-bernoulli} for $u_1$:
\begin{equation}\label{eqn-u1}
u_1 = \frac{ \sqrt{2\big(B - \frac{2\gamma p}{(\gamma -1)\rho}\big)}}{\sqrt{1+ w^2}} \qquad \mbox{with $w=\frac{u_2}{u_1}$}.
\end{equation}

\item[(ii)] Multiplying system \eqref{euler-nondiv} from the left by ${\bf l}_2$ also gives \eqref{3.30a}.

\smallskip
\item[(iii)] Multiplying  equations \eqref{euler-nondiv} from the left by ${\bf l}_3$
and separating the real and imaginary parts of the equation lead to the elliptic system:
\begin{equation}\label{eqn-wp1}
\begin{split}
&D_R w + h D_I p  =  0,   \\
&D_I w - h D_R p  = 0,
\end{split}
\end{equation}
where $D_R = \del_{y_1} + \l_R \del_{y_2}, \,  D_I = \l_I\del_{y_2},\,
\l_R =-\frac{c^2 \rho u_2}{c^2 - u_1^2 }, \,\l_I =  \frac{c\rho u_1 \sqrt{c^2 -q^2}}{c^2 - u_1^2}$,
and  $h= \frac{\sqrt{c^2-q^2}}{c\rho u_1^2}$.
\end{enumerate}
Therefore, system \eqref{eqn-euler1}--\eqref{eqn-euler4} is
decomposed into \eqref{eqn-u1}--\eqref{eqn-wp1}.

We solve this problem by iteration.
Given $U^-$ that is close to $U^-_0$ as defined in Theorem \ref{LagrCoordAppr2-thm},
we solve the problem for $U$ in the Lagrangian coordinates.
However, since $\LU^\infty$ is not known, we cannot directly
solve Problem \ref{pr-3.2}(WT) for $U$ satisfying \eqref{est-U-small-pert-Lagr},
or Problem \ref{pr-3.2}(ST) for $U$ satisfying
\eqref{est-U-small-pert2-Lagr}.
Instead, we solve Problem \ref{pr-3.2}(ST) for $U$ that is close to $U_0^+$ as in
\eqref{est-U-small-pert2}, and Problem \ref{pr-3.2}(WT) for $U$ in
similar norms with appropriate growths, but using these norms in the Lagrangian coordinates
(more precisely, the $\zz=(z_1, z_2)$--coordinates defined by \eqref{eqn-tranz} below).
Note that these norms are weaker than the ones in \eqref{est-U-small-pert-Lagr}  or
\eqref{est-U-small-pert2-Lagr} respectively; in particular, they do not determine
any limit for $U_1=(u_1, \rho)$ as $|\yy|\to\infty$ within the subsonic region.
On the other hand, these norms determine that
$(w, p)$ have the limit: $(0, p_0^+)$ at infinity within the subsonic region;
this asymptotic condition is sufficient to make the iteration problem well-defined
(in fact, we use only the asymptotic decay of $w$)
and to obtain the existence and uniqueness for the iteration problem.
After the unique solution $U$ of the problem stated in Theorem \ref{LagrCoordAppr2-thm}
is obtained by iteration,
we identify  $\LU_1^\infty=(\rho^\infty, u_1^\infty)$
and show the faster convergence of $(\rho, u_1)$ to $(\rho^\infty, u_1^\infty)$,
which lead to \eqref{est-U-small-pert-Lagr} and
\eqref{est-U-small-pert2-Lagr}, respectively.
Note that, in the estimates discussed above,
$U-U_0^+$ (rather than $U$ itself) lies in the weighted spaces \eqref{def-C}. For this
reason, it is convenient to perform the iteration in terms of
\begin{equation}\label{defDelta-iter}
 \delta U_1=U_1-U^+_{10},\qquad \delta U_2=U_2-U^+_{20},
 \qquad  \delta \hat\sigma=\hat\sigma-\hat\sigma_0=\hat\sigma-s_1y_1,
\end{equation}
where $U_1$ and $U_2$ are defined by \eqref{var-u1-w-p-rho-rot}.
Then we follow the steps below to solve this problem:

\medskip
1. {\bf Introduce a linear boundary value problem for the iteration.}
For a given shock-front $\hs$, the subsonic domain $\D^{\hs}$ is fixed, depending on $\hs$.
We make the coordinate transformation
to transform the domain from $\D^{\hs}$ to $\D$,
where $\D=\D^{\hs_0}$ with $\hs_0(y_2)=s_1y_1$ is the domain corresponding to the background solution:
\begin{equation}\label{3.12abc-0}
\D= \big\{\yy: 0 <y_2<s_1y_1\big\}
\end{equation}
with $\partial\D=\overline{\mathcal{L}_1}\cup\mathcal{L}_2$, where
\begin{equation}\label{3.12abc-1}
\mathcal{L}_1
=\{(y_1, y_2)\,:\;y_1>0, y_2=0\},
\qquad \mathcal{L}_2=\{(y_1, y_2)\,:\;y_1>0, y_2=s_1y_1\}.
\end{equation}
This transformation is:
\begin{equation}\label{eqn-tranz}
\yy=(y_1,y_2)\,\to\, \zz=(z_1,z_2):=(y_1, y_2-\d \hs(y_1)),
\end{equation}
where $\d\hs(y_1)=\hs(y_1)-\hs_0(y_1)$.
In the $\zz$--coordinates, $\mathcal{L}_1$ corresponds to $\partial W$, and
$\mathcal{L}_2$ corresponds to $\partial \sS$. Also,
$U(\yy)$ becomes $U_{\hs}(\zz)$ depending on $\hs$.
Then the upstream flow $U^-$ involves an unknown variable explicitly depending on $\hs$:
\begin{equation}\label{UdelataSigma}
U^-_{\hs}(\zz) = U^-(z_1, z_2+\d \hs(z_1)),
\end{equation}
where $U^-$ is the given upstream flow in the $\yy$--coordinates.
Equations \eqref{eqn-wp1} in $\D$ in the $\zz$-coordinates are:
\begin{equation}\label{eqn-wp-z1-0}
\begin{cases}
\widetilde{D}_R w + h \widetilde{D}_I p = 0, \\
\widetilde{D}_I w - h \widetilde{D}_R p = 0,
\end{cases}
\end{equation}
where $\widetilde{D}_R= \partial_{z_1} + (\gl_R - \gd \hs' )\partial_{z_2}$
and $\widetilde{D}_I= \gl_I \partial_{z_2}$.
Since $U_0^+$ is a constant vector and $w_0^+=0$,
the same system holds for $(\delta p, \delta w)$,
where we have used notation \eqref{defDelta-iter}.
Moreover, since the iteration:
 $(\delta U, \delta w) \to (\delta \tilde U, \delta \tilde w)$
is considered, we use $U:=U_0^++\delta U$
to determine the coefficients in \eqref{eqn-wp-z1-0} and
$(\delta \tilde p, \delta \tilde w)$ for the unknown functions.
Thus, we have
\begin{equation}\label{eqn-wp-z1}
\begin{cases}
\widetilde{D}_R \delta \tilde w + h \widetilde{D}_I \delta\tilde  p = 0,
 \\
\widetilde{D}_I \delta\tilde  w - h \widetilde{D}_R \delta\tilde  p = 0.
\end{cases}
\end{equation}
We use system \eqref{eqn-wp-z1} in $\D$ as a linear system for the iteration.

In the $\zz$--coordinates, the Rankine-Hugoniot conditions \eqref{con-RH1}--\eqref{con-RH4} keep the same form,
except that $\hs'(y_1)$ is replaced by $\hs'(z_1)$ and $U^-$ is replaced by $U^-_{\hs}$ along line $z_2=s_1 z_1$.
Among the four Rankine-Hugoniot conditions, \eqref{con-RH4} is used in the Bernoulli law.
From condition \eqref{con-RH3}, we have
\begin{equation} \label{eqn-sigmaprime}
 \hs' (z_1) = \frac{[p]}{[ u_1 w ]}(z_1, s_1z_1),
\end{equation}
which is used to update the shock-front later.
Now, because of \eqref{eqn-u1}, we can use $\bar{U}=(w, p, \rho)$ as the unknown variables along $z_2=s_1z_1$.
Using \eqref{eqn-sigmaprime} to eliminate $\hs'$ in conditions  \eqref{con-RH1}--\eqref{con-RH2} gives
\begin{eqnarray}
&&G_1(U_{\hs}^-, \bar{U}) := [p] \Big[ \frac{1}{\rho u_1 }\Big] + [w] [u_1 w] = 0, \label{con-G1}\\
&&G_2(U_{\hs}^-, \bar{U}) := [p] \Big[ u_1 + \frac{p}{\rho u_1 }\Big] + [pw] [u_1 w] = 0,\label{con-G2}
\end{eqnarray}
on $\mathcal{L}_2$. We use  conditions \eqref{con-G1}--\eqref{con-G2}
to define the linear conditions for the iteration:
$\bar U \to \widetilde{\bar{U}}$ such that, at a fixed point $\bar U= \widetilde{\bar{U}}$,
these iteration conditions imply that the original conditions
\eqref{con-G1}--\eqref{con-G2} hold. Specifically, we define the conditions:
\begin{equation}
\nabla_{\bar{U}} G_i(U_0^-, \bar{U}_0^+)\cdot \d \widetilde{\bar{U}}
= \nabla_{\bar{U}} G_i(U_0^-, \bar{U}_0^+)\cdot \d \bar{U} - G_i(U_{\hs}^-, \bar{U})
\qquad\mbox{on $\mathcal{L}_2$}
\label{con-lin-gi}
\end{equation}
for $i=1,2$,  which can be written as
\begin{equation}
b_{i1}\d \tilde{w}  + b_{i2}\d \tilde{p}  +b_{i3}\d \tilde{\rho}  =g_i (U_{\hs}^-, \bar{U})
\qquad\,\,\mbox{on $\mathcal{L}_2$},  \label{con-gi}
\end{equation}
where
$(b_{i1},b_{i2},b_{i3}) := \nabla_{\bar{U}} G_i(U_0^-, \bar{U}_0^+)$ and
$$
g_i (U_{\hs}^-, \bar{U}):= \nabla_{\bar{U}} G_i(U_0^-, \bar{U}_0^+)\cdot \d \bar{U}-G_i(U_{\hs}^-, \bar{U})
\qquad \mbox{for $i=1,2$}.
$$

Since there are two conditions in \eqref{con-gi}, $i=1,2$, we can eliminate $\d \tilde{\rho}$ to obtain
\begin{equation} \label{con-g3}
 \d \tilde{w} + b_1\d \tilde{p}=g_3
\qquad\,\,\mbox{on $\mathcal{L}_2$},
\end{equation}
where
\begin{equation}\label{5.14a}
 b_1=\frac{b_{12}b_{23}-b_{22}b_{13}}{b_{11}b_{23} - b_{21} b_{13}},  \qquad
 g_3 = \frac{b_{23}g_1-b_{13} g_2}{b_{11}b_{23} - b_{21} b_{13}}
\end{equation}
with
 \begin{eqnarray*}
b_{11}b_{23} - b_{21} b_{13}
=  (-u_{20}^-)[p_0] \left( \frac{\g p_0^+}{(\g -1)(\rho_0^+)^2 u_{10}^+} + \frac{p_0^-}{u_{10}^-}\big(\frac{1}{(\rho_0^+)^2 }
 +  \frac{\g p_0^+}{(\g -1)(\rho_0^+)^3 (u_{10}^+)^2}\big)\right) > 0.
\end{eqnarray*}

\begin{figure}
 \centering
\includegraphics[height=65mm]{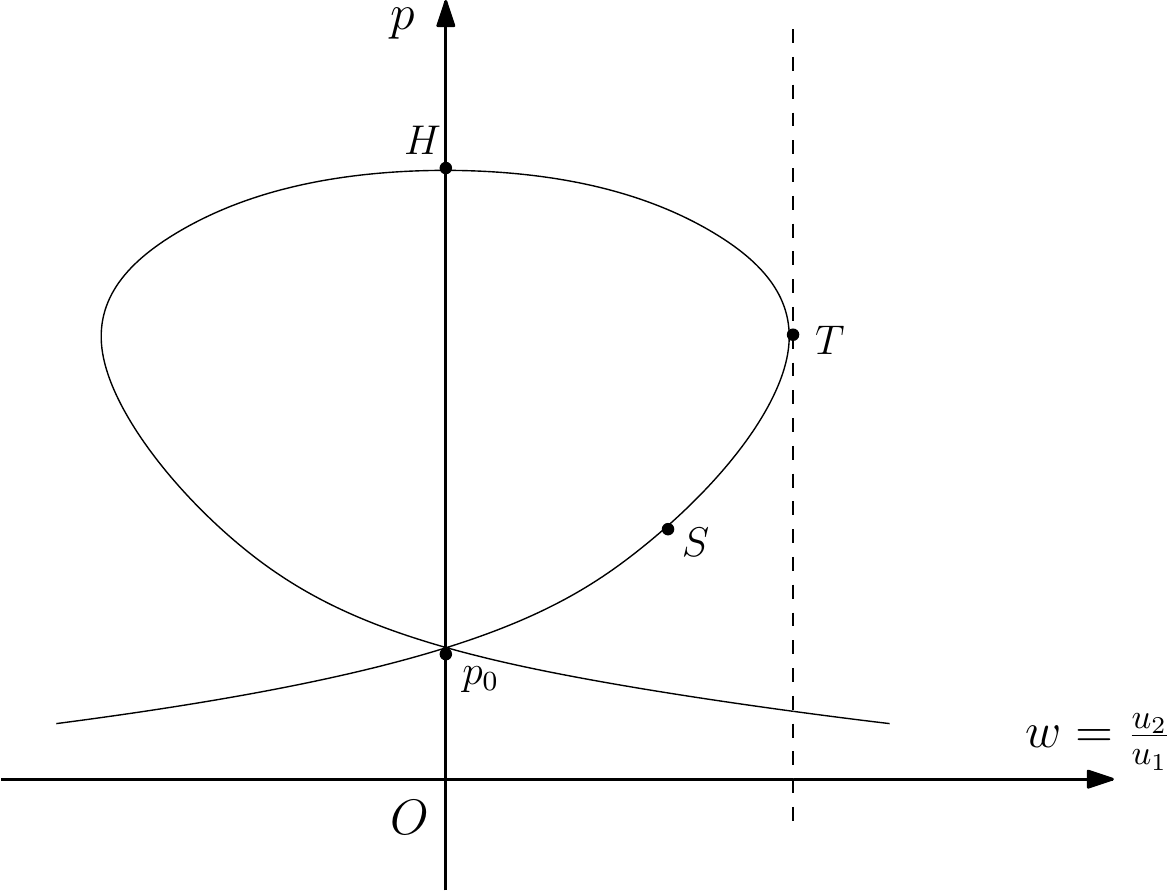}
\caption{The shock polar in the $(w, p)$--variables ({\it cf}. \cite{CCF3})}
\label{wppolar}
\end{figure}
Notice that the shock polar is a one-parameter curve determined by the Rankine-Hugoniot conditions.
If $p$ is used as the parameter, by equation \eqref{con-g3},
we obtain that $\d w= -b_1 \d p + g_3(\d p)$, which shows that $-b_1 \d p $ is the linear term and $g_3(\d p)$ is the higher order term.
We know from Fig.~\ref{wppolar} that $w(p)$ is decreasing in $p$  on arc $\wideparen{TH}$
and increasing on $\wideparen{TS}$.
Therefore, it is easy to see that
\begin{equation}\label{sgnOfCoef-weak-strong}
\mbox{\em $b_1>0$ corresponds to the state on arc $\wideparen{TH}$,
$b_1<0$ to $\wideparen{TS}$, and $b_1 =0$ at the tangent point $T$.}
\end{equation}
This difference in the sign of $b_1$  is the reason for the different rates of decay at infinity and near the origin
in the two different cases (i) and (ii) of Theorems \ref{thm-main} and \ref{LagrCoordAppr2-thm}.

It can be checked that
$$
b_{13}= -[p_0]\left(\frac{p_0^+}{(\rho_0^+)^2 u_{10}^+} +\frac{\g p_0^+}{(\g -1)(\rho_0^+)^3 (u_{10}^+)^3}\right)< 0.
$$
Thus, condition \eqref{con-gi} for $i=1$ can be rewritten as
\begin{equation} \label{con-g4}
 \d \tilde{\rho}  =g_4 - b_2 \d  \tilde{w}   -b_3  \d \tilde{p}
\qquad\mbox{on } \mathcal{L}_2,
\end{equation}
where $g_4 =\frac{g_1}{b_{13}}, b_2= \frac{b_{11}}{b_{13}}$, and $b_3=\frac{b_{12}}{b_{13}}$.

We notice that conditions \eqref{con-g3}--\eqref{con-g4} are equivalent to conditions \eqref{con-gi} for $i=1,2$.
The boundary condition on $\mathcal{L}_1$ comes from the slip condition \eqref{slipcon} on
$\partial W$.
Specifically, using \eqref{slipcon} and \eqref{wall-rot},
we obtain that $w=\brt'$ on $\del {W}$. Then, in
the $\zz$--coordinates, this must hold on $\mathcal{L}_1$.
Also, for the background solution, $\brt=0$
by \eqref{wall-rot}. Then we prescribe
\begin{equation} \label{con-L1}
 \delta \tilde w=\brt'
\qquad\mbox{on } \mathcal{L}_1.
\end{equation}

\medskip
2. {\bf Design the iteration map $\mathcal{Q}$ and prove the existence of a fixed point for $\mathcal{Q}$}.
We perform the iteration in terms of $\delta U_k$, $k=1,2$, and $\delta \hat \sigma$ as
defined by \eqref{defDelta-iter}, in the $\zz$-coordinates defined in
\eqref{eqn-tranz}. In fact, for $\hat\sigma$, we only need $\hat\sigma'$ since $\hat\sigma(0)=0$,
{\it i.e.}, the shock is attached
to the wedge vertex.  Note also that
$\delta\hat\sigma'=\hat\sigma'-s_1$.
We thus denote $V=(U_1, U_2,
\delta \hat\sigma')$ and
perform the following iteration: $\delta V\to \delta \tilde V$.
For a given $\delta V$, we determine $V:=\delta V+V_0^+$.
Then we find $\tilde V$ by solving the linear system \eqref{eqn-wp-z1}
in $\D$ with the boundary conditions \eqref{con-g3} and \eqref{con-L1},
to determine $(\tilde w, \tilde p)$
and then determine
$u_1$ from
\eqref{eqn-u1} and $\rho$ from \eqref{3.30a}
(which holds in the $\zz$--coordinates without change),
and the boundary condition
\eqref{con-g4}.
The final step is to use solution $(\delta u_1, \delta \rho, \delta w, \delta p)$ and $U^-_{\hs}$ defined by
\eqref{UdelataSigma}
on the right-hand side of \eqref{eqn-sigmaprime} to update $\delta\hat \sigma'$.
This defines the iteration map $\mathcal{Q}$ from $V$ to $\tilde V$, except that we discuss below
how the boundary value problem for \eqref{eqn-wp-z1}
with the boundary conditions \eqref{con-g3} and \eqref{con-L1} is solved
in $\D$.

As we discussed above, we perform the iteration in the spaces from
\eqref{est-U-small-pert2-Lagr} for Problem \ref{pr-3.2}(ST)
and similar norms with appropriate growths for Problem \ref{pr-3.2}(WT),
expressed in the $\zz$--coordinates
\eqref{eqn-tranz}.
We focus below on the case of Problem \ref{pr-3.2}(WT), since the other case is similar.

For $\tau>0$, define
\begin{equation}\label{spacesInZcoord}
\begin{split}
  \Sigma^\tau_1 &=
 \{v\,:\; \|v\|^{(-\ga;\mathcal{L}_1)}_{2,\ga;(0,1+\gb);\D}
+\|v_{z_1}\|^{(1-\ga;\mathcal{L}_1)}_{2,\ga;(1+\gb,1);\D}\le\tau\},\\[1mm]
 \Sigma^\tau_2 &=
 \{v\,:\; \|v\|^{(-\ga;O)(-1-\ga;\mathcal{L}_1)}_{2,\ga;(1+\gb,0);\D}\le\tau\},
 \qquad
  \Sigma^\tau_3 =
 \{v\,:\; \|v\|^{(-\ga;0)}_{2,\ga;(1+\gb);\R^+}\le\tau\},\\[1mm]
 \Sigma^\tau &=
 \{(\delta U_1, \delta U_2, \delta\hat\sigma')\,:\;
 \delta U_1\in \Sigma^\tau_1\times \Sigma^\tau_1, \;\;
 \delta U_2\in\Sigma^\tau_2\times \Sigma^\tau_2, \;\;
 \delta\hat\sigma'\in\Sigma^\tau_3\}.
\end{split}
\end{equation}
The condition on $v_{z_1}$ in $\Sigma^\tau_1$ is added for technical reasons.

It remains to discuss how we find
$(\delta \tilde w, \delta\tilde p)\in \Sigma^{C_0\varepsilon}_2\times \Sigma^{C_0\varepsilon}_2$
that solves  \eqref{eqn-wp-z1} in $\D$
with the boundary conditions \eqref{con-g3} and \eqref{con-L1}. From system \eqref{eqn-wp-z1}, we obtain
 \begin{eqnarray}
&&(\delta\tilde  p)_{z_1}=
\frac{\gl_R - \gd \hs'}{h\l_I}(\delta \tilde w)_{z_1} +
\frac{(\gl_R - \gd \hs' )^2+\l_I^2}{h\l_I}(\delta \tilde w)_{z_2},   \label{eqn-wp-z1-r}\\
&&(\delta\tilde  p)_{z_2}=
-\frac{1}{h\l_I}(\delta \tilde w)_{z_1}-
\frac{\gl_R - \gd \hs'}{h\l_I}(\delta \tilde w)_{z_2}.\label{eqn-wp-z2-r}
\end{eqnarray}
Now, differentiating and subtracting the equations, we eliminate $\delta\tilde  p$
to obtain a second-order PDE for $\delta \tilde w$ of the form:
\begin{equation}\label{eqForDeltaw}
  \sum_{i,j=1}^2(a_{ij}(\delta \tilde w)_{z_j})_{z_i}=0,
\end{equation}
where the coefficients are computed explicitly from \eqref{eqn-wp-z1-r}--\eqref{eqn-wp-z2-r}.
Note that, at the subsonic background solution  \eqref{rotatedCoor-Appr2}, we have
$$
\l_{R0}=0, \quad \l_{I0}>0, \quad h_0>0,
$$
where the left-hand sides are constants and $\delta\hat\sigma_0=0$.
Then, computing the coefficients at the
background solution, equation \eqref{eqForDeltaw} becomes
$$
\frac 1{\l_{I0}}(\delta \tilde w)_{z_1z_1}+\l_{I0}(\delta \tilde w)_{z_2z_2}=0,
$$
that is, the equation is uniformly elliptic.
Then, for the coefficients computed
at $(U^+_{10}+\delta U_1, U^+_{20}+\delta U_2, \delta \hat\sigma')$ for
$(\delta U_1, \delta U_2, \delta\hat\sigma')\in \Sigma^{C_0\varepsilon}$,
equation \eqref{eqForDeltaw} is uniformly elliptic if $\varepsilon$ is small. This
allows us to obtain the unique solution $\delta \tilde w\in \Sigma^{C_0\varepsilon}_2$ of
\eqref{eqForDeltaw} in $\D$ with the boundary conditions \eqref{con-g3} and \eqref{con-L1}.
Note that the inclusion
$\delta \tilde w\in \Sigma^{C_0\varepsilon}_2$ involves the asymptotic condition at infinity,
which makes the boundary value problem well-defined and allows us to prove the uniqueness.
After $\delta \tilde w$ is determined,
we determine $\delta\tilde p$ by the $z_2$--integration from \eqref{eqn-wp-z2-r} with the initial
condition \eqref{con-g3}, where it can be shown that $b_1 \ne 0$. Then we show
that $\delta p\in \Sigma^{C_0\varepsilon}_2$. This completes the definition of the iteration map.

The iteration set for Problem \ref{pr-3.2}(WT) is  $\Sigma^{C_0\varepsilon}$.
We first show that
$\mathcal{Q}(\Sigma^{C_0\varepsilon})\subset \Sigma^{C_0\varepsilon}$ when $\varepsilon$ is small,
and then obtain a fixed point by the Schauder fixed point theorem,
via showing that $\Sigma^{C_0\varepsilon}$ is a compact subset in the Banach space
defined by replacing $\alpha$ via $\alpha'\in (0, \alpha)$ in the norms used
in the definition of $\Sigma^\tau$
and showing that map $\mathcal{Q}$ is continuous
in this norm.

\medskip
3. {\bf  Asymptotic limit of the fixed point in the $\yy$--coordinates}.
Let $(\delta U_1, \delta U_2, \delta\hat\sigma')\in \Sigma^{C_0\varepsilon}$
be a fixed point of the iteration map, and let
$( U_1,  U_2, \hat\sigma')=(U^+_{10}+\delta U_1, U^+_{20}+\delta U_2, \delta \hat\sigma')$.

We change from the $\zz$-- to $\yy$--coordinates by inverting \eqref{eqn-tranz}:
$$
\zz:=(z_1,z_2)\,\to\, \yy:=(y_1,y_2)=(z_1, z_2+\d \hs(z_1)).
$$
Since $\delta\hat\sigma'\in \Sigma^{C_0\varepsilon}_3$, then  both \eqref{eqn-tranz} and its inverse are close
to the identity map in $C^{2,\alpha}(\D_{\hat{\sigma}}; \R^2)$
and $C^{2,\alpha}(\D; \R^2)$, respectively.
Then it follows that, in the $\yy$--coordinates,
 $(\delta U_1, \delta U_2, \delta\hat\sigma')\in \tilde\Sigma^{2C_0\varepsilon}$ if
 $\varepsilon$ is small, where
 \begin{equation}\label{spacesInZcoord-b}
\begin{split}
  \tilde \Sigma^\tau_1 &=
 \{v\,:\; \|v\|^{(-\ga;\mathcal{L}_1)}_{2,\ga;(0,1+\gb);\D_{\hat{\sigma}}}
+\|v_{z_1}\|^{(1-\ga;\mathcal{L}_1)}_{2,\ga;(1+\gb,1);\D_{\hat{\sigma}}}\le\tau\}, \quad
 \tilde \Sigma^\tau_2 =
 \{v\,:\; \|v\|^{(-\ga;O)(-1-\ga;\mathcal{L}_1)}_{2,\ga;(1+\gb,0);\D_{\hat{\sigma}}}\le\tau\},\\[1mm]
 \tilde \Sigma^\tau &=
 \{(\delta U_1, \delta U_2, \delta\hat\sigma')\,:\;
 \delta U_1\in \tilde \Sigma^\tau_1\times \tilde \Sigma^\tau_1, \;\;
 \delta U_2\in\tilde \Sigma^\tau_2\times \tilde\Sigma^\tau_2, \;\;
 \delta\hat\sigma'\in\Sigma^\tau_3\}.
\end{split}
\end{equation}
In particular, this leads to the estimates of the second and third terms on the left-hand side of
\eqref{est-U-small-pert-Lagr}.

Note that, if $v\in \tilde\Sigma^\tau_2$, then $v\to 0$ as $|\yy|\to \infty$ in $\D_{\hat{\sigma}}$,
with decay rate $|\yy|^{-(\beta+1)}$.
However, for $v\in \tilde \Sigma^\tau_1$,
no asymptotic limit  as $|\yy|\to \infty$ in $\D_{\hat{\sigma}}$ is defined.

Then, from \eqref{rotatedCoor-Appr2}--\eqref{var-u1-w-p-rho-rot},
it follows that $U_2=(w, p)\to (0, p_0^+)$ as $|\yy|\to \infty$ in $\D$;
however, for $U_1=(u_1, \rho)$, the limit is not determined by space $\Sigma^\tau_1$,
and $(u_1, \rho)$ does not converge to $(u_{10}^+, \rho_0^+)$ in general,
as we see below. Thus, we have to determine the limiting profiles
$(u_1^\infty(y_2), \rho^\infty(y_2))$.

To determine $\rho^\infty(y_2)$, we first obtain \eqref{3.30a} from  \eqref{eqn-euler1}--\eqref{eqn-euler4},
which implies \eqref{eqn-rho-p}.
Since $\hat\sigma(y_1)$ is determined,
then $A(y_2)$ in \eqref{eqn-rho-p} is determined by the upstream state
$U^-(\yy)$ from the Rankine-Hugoniot conditions \eqref{con-RH1}--\eqref{con-RH4}.
Then,
noting that $p\to p_0^\infty$, we obtain formally
$$
\rho \to \rho^\infty(y_2)=\Big(\frac{p_0^+}{A(y_2)}\Big)^{\frac 1\gamma}\qquad
\mbox{as $|\yy|\to\infty$ in $\D_{\hat{\sigma}}$}.
$$
Similarly, we use \eqref{eqn-u1} to obtain
$$
u_1\to u_1^\infty(y_2) =
  \sqrt{2\Big(B(y_2) - \frac{\g p_0^+ }{ (\g -1 )\rho^\infty(y_2) }\Big)}
\qquad\mbox{as $|\yy|\to\infty$ in $\D_{\hat{\sigma}}$}.
$$
Defining $\LU^\infty(y_2)=(u_1^\infty(y_2),0, p_0^+,\rho^\infty(y_2))$, we can show that
the estimates of the first and the last terms on the left-hand side of \eqref{est-U-small-pert-Lagr}
hold. This completes the argument for case (i) of Theorem \ref{LagrCoordAppr2-thm}.

Case (ii) is handled similarly. Note that the slower decay at infinity for case (ii), {\it i.e.},
$|\yy|^{-\beta}$, is from the elliptic estimates, even if the faster decay at infinity in
\eqref{cond-U-small-pert2} is required.
The reason for the difference in the rates for cases (i) and (ii) is
\eqref{sgnOfCoef-weak-strong}.

\medskip
4. {\bf Return to the $\xx$-coordinates}.
We obtain Theorem \ref{thm-main-rot} directly from Theorem \ref{LagrCoordAppr2-thm} by changing the coordinates.
Recall that, when the Lagrangian coordinates are defined for Theorem \ref{LagrCoordAppr2-thm},
we have used the rotated coordinates $\xx$ in \eqref{def-coord};
see the discussion in the paragraph before Theorem \ref{LagrCoordAppr2-thm}.

From the estimates in Theorem  \ref{LagrCoordAppr2-thm}, it follows that, in the Lagrangian coordinates,
$|U-U_0^+|\le C\varepsilon$ in $\D_{\hat{\sigma}}$, where $C$ depends only on $(U_0^-, U_0^+)$.
Thus, the same is true in the $\xx$--coordinates in $\Omega$.
Then it follows from \eqref{psi-potential}--\eqref{def-coord} and \eqref{rotatedCoor-Appr2} for positive $u_{10}^+$ and $\rho_0^+$
that the change of coordinates $\xx \to \yy$ given by \eqref{def-coord} is bi-Lipschitz.
Then \eqref{est-U-small-pert2-rot} follows from \eqref{est-U-small-pert2-Lagr} directly.

Similarly, the estimates of the second and third  terms on the left-hand side of
\eqref{est-U-small-pert-rot}  follow from \eqref{est-U-small-pert-Lagr} directly.
In order to obtain the estimates of the remaining terms on the left-hand side of
\eqref{est-U-small-pert-rot}, we need to identify $U^\infty(x_2)$.

Note that, on shock $\sS$, using \eqref{def-coord} and the estimate of the third term
on the left-hand side of
\eqref{est-U-small-pert-rot}, we see that, for small $\varepsilon$,
$$
\partial_{\ttau_{\sS}}\psi=\rho \uu\cdot\nnu_{\sS}\ge \rho \uu_0^+\cdot\nnu_{\sS_0}-C\varepsilon
\ge \frac 12 \rho \uu_0^+\cdot\nnu_{\sS_0}>0.
$$
Recall also that $\psi({\bf 0})=0$ by \eqref{psi-potential}.
Then, for each $y_2>0$, there exists a unique $\xx^{\rm in}(y_2)=(x_1^{\rm in}(y_2),
x_2^{\rm in}(y_2))\in \sS$ such that $\psi(\xx^{\rm in}(y_2))=y_2$ and
$$
\|\xx^{\rm in}\|_{C^{2,\alpha}([0, \infty))}\le C, \qquad\,\,
(\xx^{\rm in})'\ge C^{-1}>0 \;\;\mbox{ on $[0, \infty)$}.
$$
From this and \eqref{psi-potential}, it follows that, for each $y_2>0$,
$$
\Omega\cap\{\xx\,:\, \psi(\xx)=y_2\}=\{(x_1, x_2^*(x_1; y_2))\;:\; x_1>x_1^{\rm in}(y_2)\},
$$
where $ x_2^*(\cdot; y_2)$ is the solution of the initial value problem for the differential equation:
\begin{equation}\label{ODE-streamline}
 \begin{cases}
    \partial_{x_1} x_2^*(x_1; y_2)  =w(x_1, x_2^*(x_1; y_2)), \\[1mm]
    x_2^*(x_1^{\rm in}(y_2); y_2)=x_2^{\rm in}(y_2),
 \end{cases}
\end{equation}
where $w=\frac{u_2}{u_1}$ ({\it cf}. \eqref{var-u1-w-p-rho-rot}).
Since we have obtained the estimate
of the second term  on the left-hand side of
\eqref{est-U-small-pert}, using \eqref{rotatedCoor-Appr2}, we have
\begin{equation}\label{estWfor-streamline}
  |D^kw(\xx)|\le C_0\varepsilon (1+|\xx|)^{-1-\beta}
  \qquad\mbox{in $\Omega,\;\;$ for $k=0,1,2$}.
\end{equation}
In particular, for each $y_2\ge 0$ and $k=0, 1,2$,
\begin{equation}\label{calcIntUsingDecay}
  \int_{x_1^{\rm in}(y_2)}^\infty |D^kw(x_1, x_2^*(x_1; y_2))|\,{\rm d}x_1
\le C_0\varepsilon \int_{0}^\infty (1+x_1)^{-1-\beta}\,{\rm d}x_1\le C\varepsilon.
\end{equation}
Applying this with $k=1$, we conclude that $\lim_{x_1\to\infty} x_2^*(x_1; y_2)$ exists for each $y_2\ge 0$,
which is denoted as $x_2^\infty(y_2)$.

Recall the structure  of $\Omega$ in \eqref{domain:1-rot},
where $\brt(x_1)\to 0$ and $\tilde \sigma(x_1)\to \infty$ as $x_1\to \infty$
by \eqref{cond-U-small-pert-rot} and the estimate of the third term in \eqref{est-U-small-pert-rot}.
Differentiating \eqref{ODE-streamline} twice with respect to $y_2$ and using
the $C^2$--estimate of $\xx^{\rm in}$ and \eqref{calcIntUsingDecay},
we obtain that
$\|x_2^*(x_1; \cdot)\|_{C^2([\brt(x_1), \tilde \sigma(x_1)])}\le C$ for each $x_1>0$.
From this, we have
\begin{equation}\label{propX2-infty}
\mbox{$x_2^*(x_1; \cdot)\to x_2^\infty(\cdot)$ in $C^1$ on compact subsets on $[0, \infty)$
as $x_1\to \infty$},
\end{equation}
with $\|x_2^\infty\|_{C^2([0, \infty))}\le C$. Also, by a similar argument, using the $C^{2,\alpha}$--regularity of $\xx^{\rm in}$
and the estimate of $w$ in the second term in \eqref{est-U-small-pert-rot}, we obtain
that $x_2^\infty\in C^{2, \alpha}([0, \infty))$.

Furthermore, for the background solution,  the potential functions $\psi_0^-$ of $U_0^-$,
 $\psi_0^+$ of $U_0^+$, and $\psi_0$ of the transonic shock solution $(U_0^-, U_0^+)$ in $\{x_1>0, x_2>0\}$ are:
$$
\psi_0^-(\xx)=\rho^-_0u^-_{10}(x_1-x_2\tan\theta_{\rm w}),\quad
\psi_0^+(\xx)=\rho^+_0u^+_{10}x_1, \quad
\psi_0(\xx)=
\begin{cases}
 \psi_0^-(\xx) \;\;\;\mbox{ if } x_2< \tilde s x_1,\\[0.5mm]
  \psi_0^+(\xx) \;\;\;\mbox{ if } x_2> \tilde s x_1,
\end{cases}
$$
by using \eqref{rotatedCoor-Appr2}, where $\psi_0$ is Lipschitz. Then, estimating $\psi-\psi_0^-$ in $\Omega^-$ via
\eqref{cond-U-small-pert-rot} (where the polynomial decay is of degree $-(1+\beta)$ so
that the calculations similar to  \eqref{calcIntUsingDecay} can be used)
and using the Rankine-Hugoniot
conditions on $\sS$, we obtain
$$
|(\xx^{\rm in})'-(\xx^{\rm in}_0)'|\le C\varepsilon \qquad\mbox{ on $[0, \infty)$},
$$
where $\xx^{\rm in}_0(y_2)=\frac{y_2}{\rho^+_0u^+_{10}}(1, \tilde s_0)$
that is the corresponding function $\xx^{\rm in}$ of the background solution.

Denote by $x^*_{20}(x_1; y_2)$ the corresponding function $x_2^*(x_1; y_2)$ of the background solution.
Then
$$
x^*_{20}(x_1; y_2)=\frac{y_2}{\rho^+_0u^+_{10}}\qquad
\,\mbox{on $x_1>\frac{y_2} {\rho^+_0u^+_{10}\tilde s_0}$ for each $y_2\ge 0$}.
$$
Thus, $x^*_{20}(x_1; y_2)$ is independent of $x_1$, so that $x^*_{20}(x_1; y_2)=x^*_{20}( y_2)$.
Then, denoting
$$
g(x_1; y_2)= x^*_2(x_1; y_2)-x^*_{20}(y_2),
$$
we see that $g$ satisfies
\begin{equation}\label{ODE-streamline-dif}
 \begin{cases}
    \partial_{x_1} g(x_1; y_2)  =w(x_1, x_2^*(x_1; y_2)), \\[1mm]
     |g(x_1^{\rm in}(y_2); y_2)|\le C\varepsilon.
 \end{cases}
\end{equation}
From this and \eqref{calcIntUsingDecay}--\eqref{propX2-infty},
we obtain that $|(x_2^\infty)'-(x_{20}^\infty)'|\le C\varepsilon$,
where $(x_{20}^\infty)'(y_2)=(x^*_{20})'( y_2)=\frac{1}{\rho^+_0u^+_{10}}$.
Therefore, we have
$$
(x_2^\infty)'\ge \frac{1}{2\rho^+_0u^+_{10}}\quad\mbox{ on $[0, \infty)$},
$$
if $\varepsilon$ is small. In particular, noting that $x_2^\infty(0)=0$ since $\partial W$ is a
streamline corresponding to $\psi=0$ and $\lim_{x_1\to\infty} \brt(x_1)=0$ by \eqref{cond-U-small-pert-rot}, we obtain
$x_2^\infty([0, \infty))=[0, \infty)$.
Then there exists the inverse function $y_2^*:[0, \infty)\to [0, \infty)$ to $x_2^\infty(\cdot)$ such that
$y_2^*\in C^{2, \alpha}([0, \infty))$
with $y_2^*(0)=0$ and
$(y_2^*)'\ge \frac 1C>0$.

Finally, defining $U^\infty(x_2)=\LU^\infty( y_2^*(x_2))$, we obtain
\eqref{est-U-small-pert-rot} from \eqref{est-U-small-pert-Lagr}.

\medskip
For more details, see Chen-Chen-Feldman \cite{CCF3}.

\vspace{3pt}
\begin{remark}
For the global stability of weak transonic shocks for the $3$-D wedge problem,
see {\rm \cite{C-Fang-2,CCX}}; also see the instability phenomenon for
strong transonic shocks for the $3$-D wedge problem in {\rm \cite{LXY2}}.
For the global stability of conical shocks for the M-D conic problem,
see {\rm \cite{C-Fang}} for the transonic shock case
and {\rm \cite{CKZ21,CXY,Lie}} for the supersonic shock case.
\end{remark}

\smallskip
\section{Two-Dimensional Transonic Shocks and Free Boundary Problems
for the Self-Similar Euler Equations for Potential Flow} \label{self-simSect}

In \S 2--\S 3, we have discussed the free boundary problems for steady transonic shock solutions of the compressible Euler equations.
Now we discuss free boundary problems for time-dependent solutions.

Time-dependent solutions with shocks of the Cauchy problem for the compressible Euler system
may exhibit non-uniqueness in general; see \cite{CDK,KKMM} and the references cited therein.
On the other hand, many fundamental physical phenomena, including shock reflection/diffraction,
are determined by the time-dependent solutions of self-similar structure;
moreover, the uniqueness can be established in a carefully chosen class
of self-similar solutions with shocks.
In this section, we focus on this case; more precisely, we describe transonic shocks and free boundary problems
for self-similar shock reflection/diffraction for the Euler equations for potential flow.

The two-dimensional compressible potential flow is governed by the conservation law of mass and the Bernoulli law
for the density function $\rho$ and the velocity potential $\Phi$ ({\it i.e.}, $\uu=\nabla \Phi$):
\begin{align}
\label{1-a}
&\der_t\rho+ \nabla_{\bf x}\cdot (\rho \nabla_{\bf x}\Phi)=0,\\
\label{1-b}
&\der_t\Phi+\frac 12|\nabla_{\bf x}\Phi|^2+h(\rho)=B
\end{align}
for $t\in \R^+:=(0,\infty)$ and ${\bf x}\in \R^2$,
where $B$ is the Bernoulli constant,
and $h(\rho)$ is given by
\begin{equation}\label{1-c}
h(\rho)=\frac{\rho^{\gam-1}-1}{\gam-1} \qquad\,\, \mbox{for the adiabatic exponent $\gam>1$.}
\end{equation}

By \eqref{1-b}--\eqref{1-c},
$\rho$ can be expressed as
\begin{equation}\label{1-b1}
\rho(\der_t\Phi,\nabla_{\bf x}\Phi)=h^{-1}(B-\der_t\Phi-\frac 12|\nabla_{\bf x}\Phi|^2).
\end{equation}
Then system \eqref{1-a}--\eqref{1-b} can be rewritten as the following second-order nonlinear wave equation:
\begin{equation}
\label{1-b2}
\der_t\rho(\der_t\Phi, \nabla_{\bf x}\Phi)
+\nabla_{\bf x}\cdot\big(\rho(\der_t\Phi, \nabla_{\bf x}\Phi)\nabla_{\bf x}\Phi\big)=0
\end{equation}
with $\rho(\der_t\Phi, \nabla_{\bf x}\Phi)$ determined by \eqref{1-b1}.

Note that equation \eqref{1-b1} is invariant under the self-similar scaling:
\begin{equation}\label{4.5}
(t, {\bf x})\rightarrow (\alp t, \alp{\bf x}),\quad \Phi\rightarrow  \frac{\Phi}{\alp}\qquad\quad \text{for}\;\;\alp\neq 0,
\end{equation}
and thus it admits  self-similar solutions in the form of
\begin{equation}\label{4.6}
\Phi(t, {\bf x})=t\phi(\xxi)\qquad\quad \text{for}\;\;\xxi=\frac{{\bf x}}{t}.
\end{equation}
Then the pseudo-potential function
$$
\vphi(\xxi)=\phi(\xxi)-\frac 12|\xxi|^2
$$
satisfies
the following equation:
\begin{equation}
\label{2-1}
{\rm div}(\rho(|D\vphi|^2,\vphi)D\vphi)+2\rho(|D\vphi|^2,\vphi)=0
\end{equation}
for
\begin{equation}
\label{1-o}
\rho(|D\vphi|^2,\vphi)=
\bigl(B_0- (\gamma-1)(\frac{1}{2}|D\vphi|^2+\varphi)\bigr)^{\frac{1}{\gam-1}},
\end{equation}
where $B_0=(\gamma-1)B+1$,
and the divergence ${\rm div}$ and gradient $D$ are with respect to $\xxi\in \R^2$.

Equation \eqref{2-1}
written in the non-divergence form is
\begin{equation}\label{nondivMainEq}
(c^2-\varphi_{\cxi}^2)\varphi_{\cxi\cxi}-2\varphi_\cxi\varphi_\ceta\varphi_{\cxi\ceta}
+(c^2-\varphi_{\ceta}^2)\varphi_{\ceta\ceta}+2c^2-|D\varphi|^2=0,
\end{equation}
where the sonic speed $c=c(|D\varphi|^2,\varphi)$ is determined by
\begin{equation}\label{c-through-density-function}
c^2(|D\varphi|^2,\varphi)=
\rho^{\gamma-1}(|D\varphi|^2,\varphi)
=B_0-(\gamma-1)\big(\frac{1}{2}|D\varphi|^2+\varphi\big).
\end{equation}
Another form of \eqref{nondivMainEq}, which uses both the potential $\phi$ and
the pseudo-potential $\varphi$, is
\begin{equation}\label{equ:study}
(c^2-\varphi_{\xi_1}^2)\phi_{\xi_1\xi_1}-2\varphi_{\xi_1}\varphi_{\xi_2}\phi_{\xi_1\xi_2}+
(c^2-\varphi_{\xi_2}^2)\phi_{\xi_2\xi_2}=0.
\end{equation}

Equation \eqref{2-1} is a nonlinear PDE of mixed elliptic-hyperbolic type. It is elliptic
at $\xxi$ if and only if
\begin{equation}
\label{1-f}
|D\vphi|<c(|D\vphi|^2,\vphi)   \qquad \,\, \mbox{at $\xxi$},
\end{equation}
and is hyperbolic if the opposite inequality holds.
This can be seen more clearly from  the rotational invariance of \eqref{nondivMainEq},
by fixing $\xxi$ and choosing coordinates $(\xi_1, \xi_2)$ so that $\xi_1$ is along the direction of $D\vphi(\xxi)$.

Moreover,  from  (\ref{nondivMainEq})--(\ref{c-through-density-function}),
equation \eqref{2-1} satisfies the Galilean invariance property:
If $\varphi(\xxi)$ is a solution, then its shift $\varphi(\xxi-\xxi_0)$ for any constant vector $\xxi_0$ is also a solution.
Furthermore, $\varphi(\xxi)+{\it const.}$ is a solution of  \eqref{2-1}
with adjusted constant $B$ correspondingly in \eqref{1-o} and \eqref{c-through-density-function}.

One class of solutions of (\ref{2-1})
is that of {\em constant states} that are the solutions
with constant velocity $\mathbf{v}\in \mathbb{R}^2$.
This implies that the pseudo-potential of a constant state satisfies
$D\varphi=\mathbf{v}-\xxi$ so that
\begin{equation}\label{constantStatesForm}
\varphi(\xxi)=-\frac 12|\xxi|^2+\mathbf{v}\cdot\xxi +C,
\end{equation}
where $C$ is a constant. For such $\varphi$, the expressions
in \eqref{1-o} and \eqref{c-through-density-function}
imply that the density and sonic
speed are positive constants $\rho$ and $c$, {\it i.e.},
independent of $\xxi$.
Then, from \eqref{elliptic}
and \eqref{constantStatesForm},
the ellipticity condition for the constant state is
$$
|\xxi -\mathbf{v}|<c.
$$
Thus, for a constant state $\mathbf{v}$,
equation (\ref{2-1})
is elliptic inside the {\em sonic circle},
with center $\mathbf{v}$ and radius $c$, and hyperbolic outside this circle.

We also note that,
if density $\rho$ is a constant, then the solution is a constant state;
that is, the corresponding
pseudo-potential $\vphi$ is of form \eqref{constantStatesForm}.

Since the problem involves transonic shocks, we have to consider
weak solutions of equation \eqref{2-1},
which admit shocks.
As in \cite{ChenFeldman},
the weak solutions are defined in the distributional sense
in a domain $\Lambda$ in the $\xxi$--coordinates.

\begin{definition}\label{def:weak solution}
A function $\varphi\in W^{1,1}_{\rm loc}(\Lambda)$ is called a weak solution of
\eqref{2-1} if
\begin{enumerate}
\item[\rm (i)]
$B_0-(\gamma-1)(\frac{1}{2}|D\varphi|^2+\varphi)\geq0\quad\text{a.e. in }\Lambda$,

\smallskip
\item[\rm (ii)]
$(\r(|D\varphi|^2,\varphi),\r(|D\varphi|^2,\varphi)|D\varphi|)\in(L^1_{\rm loc}(\Lambda))^2$,

\smallskip
\item[\rm (iii)]
For every $\z\in C^{\infty}_{\rm c}(\Lambda)$,
\begin{equation}\label{def:weak solution-i3-Eqn}
  \int_{\Lambda}\big(\r(|D\varphi|^2,\varphi)D\varphi\cdot
  D\z-2\r(|D\varphi|^2,\varphi)\z\big)\mathrm{d}\xxi=0.
 \end{equation}
\end{enumerate}
\end{definition}

A shock is a curve across which $D\vphi$ is discontinuous. If $\Lambda^+$ and $\Lambda^-(:=\Lambda\setminus \ol{\Lambda^+})$
are two nonempty open subsets of a domain $\Lambda\subset \R^2$, and $\CS:=\der\Lambda^+\cap \Lambda$ is a $C^1$-curve
where $D\vphi$ has a jump, then $\vphi\in
C^1(\Lambda^{\pm}\cup \CS)\cap C^2(\Lambda^{\pm})$
is a global weak solution of \eqref{2-1} in $\Lambda$ if and only if $\vphi$ is in $W^{1,\infty}_{\rm loc}(\Lambda)$
and satisfies equation \eqref{2-1} and the Rankine-Hugoniot condition on $\CS$:
\begin{equation}
\label{1-h}
\rho(|D\vphi|^2, \vphi)D\vphi\cdot\nnu|_{\Lambda^+\cap \CS}
=\rho(|D\vphi|^2, \vphi)D\vphi\cdot\nnu|_{\Lambda^-\cap \CS}.
\end{equation}
Note that the condition $\vphi\in W^{1,\infty}_{\rm loc}(\Lambda)$ requires that
\begin{equation}
\label{1-i}
\vphi_{\Lambda^+\cap \CS}=\vphi_{\Lambda^-\cap \CS},
\end{equation}
which is consistent with ${\rm curl}(D\vphi)=0$ in the distributional sense.

A piecewise smooth solution with the discontinuities is called an {\it entropy solution}
of \eqref{2-1} if it satisfies the entropy condition: density $\rho$ increases in the pseudo-flow direction of
$D\varphi_{\Lambda^+\cap \CS}$ across any discontinuity. Then such a discontinuity is called a shock.

\subsection{The von Neumann Problem for Shock Reflection-Diffraction}
\label{RegReflProblSect}
\begin{figure}[htp]
\begin{center}
	\begin{minipage}{0.50\textwidth}
		\centering
		\includegraphics[width=0.63\textwidth]{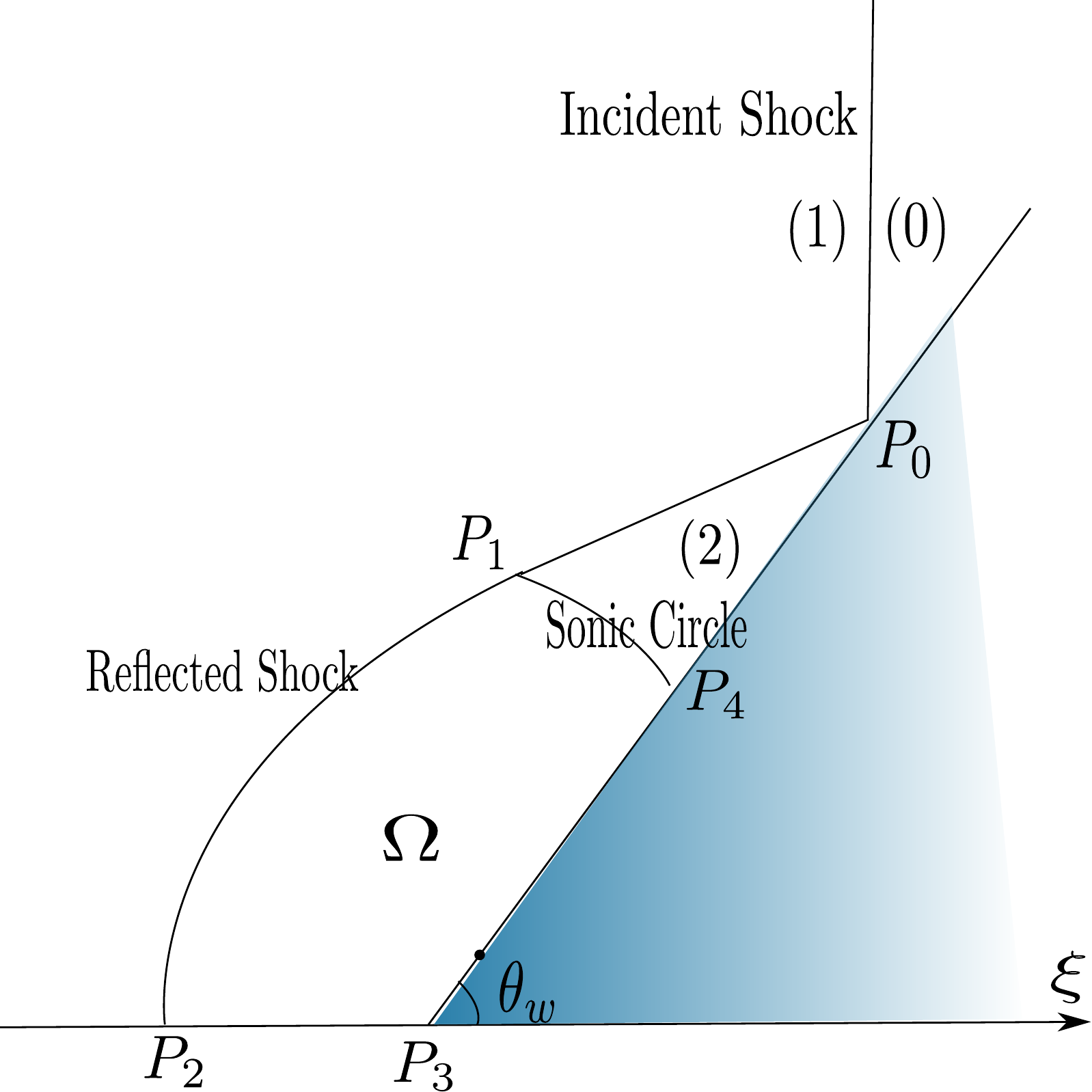}
		\caption{Supersonic regular shock reflection-diffraction configuration}
		\label{figure: free boundary problems-1}
	\end{minipage}
	\hspace{-0.5in}
	\begin{minipage}{0.49\textwidth}
		\centering
		\includegraphics[width=0.63\textwidth]{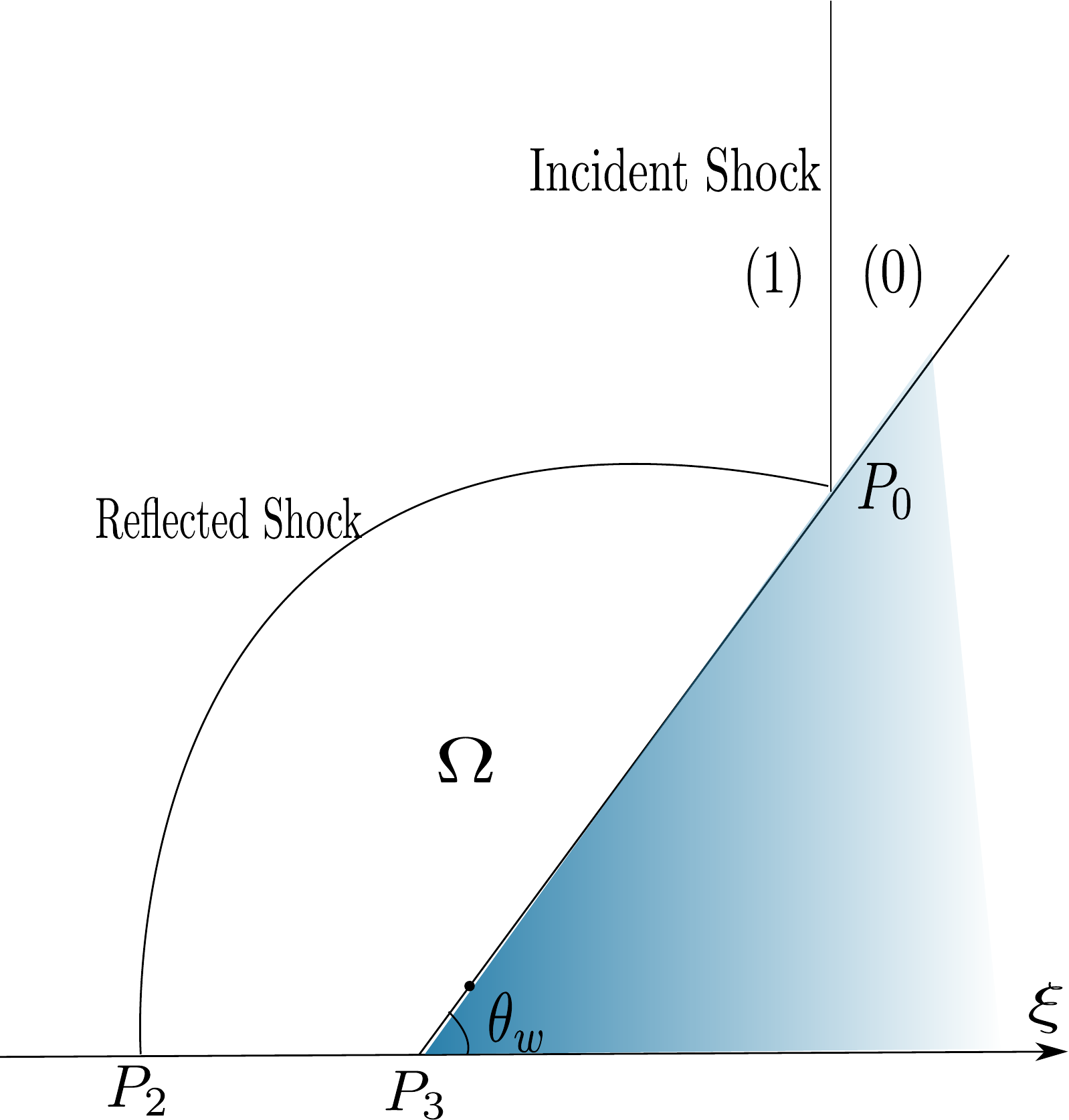}
		\caption{Subsonic$\,\,\,$  regular shock reflection-diffraction configuration}
		\label{figure: free boundary problems-2}
	\end{minipage}
\end{center}
\end{figure}

We now describe the von Neumann problem for shock reflection-diffraction, proposed for mathematical analysis first
in \cite{Neumann1,Neumann2,Neumann}.
When a vertical planar shock perpendicular to the flow
direction $x_1$ and separating two uniform states (0) and (1),
with constant velocities
$\mathbf{u}_0= (0, 0)$ and $\mathbf{u}_1=(u_1, 0)$
and constant densities $\rho_0<\rho_1$
(state (0) is ahead or to the right of the shock, and state
(1) is behind the shock), hits a symmetric wedge:
$$
W:= \big\{(x_1,x_2)\,:\, |x_2| < x_1 \tan \theta_{\rm w}, x_1 > 0\big\}
$$
head-on at time $t = 0$,
a reflection-diffraction process takes place when $t > 0$.
Then a fundamental question is what types of wave patterns of
reflection-diffraction configurations may be formed around the wedge.
The complexity of reflection-diffraction configurations was first reported
by Ernst Mach \cite{Mach} in 1878, who first observed two patterns of
reflection-diffraction configurations:
Regular reflection (two-shock configuration; see {\it e.g.}
Figs. \ref{figure: free boundary problems-1}--\ref{figure: free boundary problems-2})
and Mach reflection (three-shock/one-vortex-sheet configuration);
also see \cite{BD,CF-book2018,CF,VD}.
The issues remained dormant until the 1940s when John von Neumann \cite{Neumann1,Neumann2,Neumann},
as well as other mathematical/experimental
scientists ({\it cf.} \cite{BD,CF-book2018,CF,GlimmMajda,VD}
and the references cited therein),
began extensive research into all aspects of shock reflection-diffraction phenomena,
due to its importance in applications.
It has been found that the situations are much more complicated than
what Mach originally observed: The Mach reflection can be further
divided into more specific sub-patterns, and various other patterns of
shock reflection-diffraction configurations may occur such as the double Mach reflection,
the von Neumann reflection, and the Guderley reflection;
see \cite{BD,CF-book2018,CF,GlimmMajda,Guderley,VD}
and the references cited therein.
Then the fundamental scientific issues include:
\begin{itemize}
\item[(i)] Structures of the shock reflection-diffraction configurations;

\item[(ii)] Transition criteria among the different patterns of shock
reflection-diffraction configurations;

\item[(iii)] Dependence of the patterns upon the physical parameters such as the
wedge angle $\theta_{\rm w}$, the incident-shock-wave Mach number, and the
adiabatic exponent $\gamma>1$.
\end{itemize}

\noindent
In particular, several transition criteria among the different
patterns of shock reflection-diffraction configurations have been proposed,
including the sonic conjecture and the detachment conjecture by von Neumann
\cite{Neumann1,Neumann2,Neumann}.

A careful asymptotic analysis has been made for various reflection-diffraction
configurations in Lighthill \cite{Lighthill1,Lighthill2},
Keller-Blank \cite{KB}, Hunter-Keller \cite{HK}, Harabetian \cite{Harabetian},
Morawetz \cite{Morawetz2}, and the references
cited therein; also see Glimm-Majda \cite{GlimmMajda}.
Large or small scale numerical simulations have been also performed;
{\it cf.} \cite{BD,GlimmMajda,WC} and the references cited therein.
However, most of the fundamental issues for shock reflection-diffraction
phenomena have not been understood, especially the global structure and
transition between the different patterns of shock reflection-diffraction configurations.
This is partially because physical and numerical experiments are
hampered by many difficulties and have not yielded clear transition criteria
between the different patterns. In particular, numerical dissipation or physical
viscosity smear the shocks and cause boundary layers that interact with
the reflection-diffraction patterns and can cause spurious Mach steams; {\it cf.}
\cite{WC}. Furthermore, some different patterns occur when
the wedge angles are only fractions of a degree apart,
a resolution as yet unreachable even by
sophisticated
experiments
({\it cf.} \cite{BD,LD}).
For this reason, it is impossible to
distinguish experimentally between the sonic and detachment criteria clearly,
as pointed out in \cite{BD}.
In this regard, the necessary approach to understand
fully the shock reflection-diffraction phenomena, especially the transition criteria,
is via rigorous mathematical analysis.
To achieve this, it is essential to formulate the shock reflection-diffraction problem
as a free boundary problem and establish the global existence, regularity,
and structural stability of its solution.

\smallskip
Mathematically, the shock reflection-diffraction problem
is a two-dimensional lateral Riemann problem in domain
$\R^2\setminus \overline{W}$.

\begin{problem}[Two-Dimensional Lateral Riemann Problem]\label{ibvp-c}
{\it
Piecewise constant initial data, consisting of state $(0)$ with velocity $\mathbf{u}_0=(0,0)$
on $\{x_1>0\}\setminus \overline{W}$
and state $(1)$ with velocity $\mathbf{u}_1=(u_1, 0)$
on $\{x_1 < 0\}$ connected by a shock at $x_1=0$,
are prescribed at $t = 0$.
Seek a solution of the Euler system \eqref{1-a}--\eqref{1-b} for $t\ge 0$ subject to these initial
data and the boundary condition $\nabla\Phi\cdot\nnu=0$ on $\partial W$.
}
\end{problem}

In order to define the notion of weak solutions of {Problem \ref{ibvp-c}}, it is noted that
the boundary condition can be written as $\rho\nabla\Phi\cdot\nnu=0$ on $\partial W$,
which is spatial conormal to equation \eqref{1-b2}.
Then we have

\begin{definition}[Weak Solutions of {Problem \ref{ibvp-c}}]\label{weakSol-def-Prob1}
A function
$
\Phi\in W^{1,1}_{\rm loc}(\R_+\times (\bR^2\setminus W))
$
is called a weak
solution of  {\rm  Problem \ref{ibvp-c}}
if $\Phi$ satisfies the following properties{\rm :}
\begin{enumerate}[\rm (i)]
\item
\label{weakSol-def-i1a-Prob1}
$B-\left(\partial_t\Phi+\frac{1}{2}|\nabla_\x\Phi|^2\right)\ge h(0+)$ a.e. in
$\R_+\times (\R^2\setminus W)$.

\smallskip
\item
\label{weakSol-def-i2a-Prob1}
For $\rho(\partial_t \Phi, \nabla_\x\Phi)$ determined by \eqref{1-b1},
$$
(\rho(\partial_t\Phi, |\nabla_\x\Phi|^2), \rho(\partial_t\Phi, |\nabla_\x\Phi|^2)|\nabla_\x\Phi|)\in
(L^1_{\rm loc}(\overline\R_+\times \overline{\R^2\setminus W}))^2.
$$

\item
\label{weakSol-def-i3a-Prob1}
For every $\zeta\in C^\infty_c(\overline{\R_+}\times \bR^2)$,
\begin{align*}
&\int_0^\infty\int_{\R^2\setminus W}\Big(\rho(\partial_t\Phi, |\nabla_\x\Phi|^2)\partial_t\zeta
+
\rho(\partial_t\Phi, |\nabla_\x\Phi|^2)
\nabla\Phi\cdot\nabla\zeta\Big) \,\dr \x \dr t
+\int_{\R^2\setminus W}\rho(0,\x)\zeta(0, \x)\dr\x=0,
\end{align*}
\end{enumerate}
where
$$
\rho|_{t=0} =\begin{cases}
\rho_0 \qquad\,\, &  \mbox{for $|x_2|>x_1\tan\theta_{\rm w}$ and $x_1>0$},\\
\rho_1 \qquad\,\, &\mbox{for}\,\, x_1<0.
\end{cases}
$$
\end{definition}

\begin{remark}\label{conormCondRmk} Since $\zeta$ does not need to be zero on $\partial\Lambda$,
the integral identity in Definition {\rm \ref{weakSol-def-Prob1}}
is a weak form of equation {\rm \eqref{1-b2}}
and the boundary condition $\rho\nabla\Phi\cdot\nnu=0$ on $\partial W$.
\end{remark}

\begin{remark}\label{conormCondRmk-b}
A weak solution is called an entropy solution if it satisfies the entropy condition
that is consistent with the second law of thermodynamics
$(${\it cf}. {\rm \cite{CF-book2018,CF,Da,Lax}}$)$. In particular, a piecewise smooth solution
is an entropy solution if the discontinuities are all shocks.
\end{remark}

\medskip
Notice that {Problem \ref{ibvp-c}}
is invariant under scaling \eqref{4.5}, so it admits self-similar solutions determined by
equation \eqref{2-1} with \eqref{1-o}, along with the appropriate boundary conditions,
through \eqref{4.6}.
We now show how such solutions
in self-similar coordinates $\xxi=(\xi_1, \xi_2)=\frac{{\bf x}}{t}$
can be constructed.

First, by the symmetry of the problem with respect to the $\xi_1$--axis,
we consider only the upper half-plane
$\{\xi_2>0\}$ and prescribe the boundary condition: $\varphi_{\bn}=0$
on the symmetry line $\{\xi_2=0\}$.
Note that state (1) satisfies this condition.
Then {Problem \ref{ibvp-c}} is reformulated as a boundary value problem
in the unbounded domain:
$$
\Lambda:=\R^2_+\setminus\{\xxi\,:\,|\xi_2|\le \xi_1 \tan\theta_{\rm w}, \xi_1>0\}
$$
in the self-similar coordinates $\xxi=(\xi_1,\xi_2)$, where $\mR^2_+:=\mR^2\cap\{\xi_2>0\}$.
The incident shock in
the self-similar coordinates is the half-line $\CS_0=\{\xi_1=\xi_1^0\}\cap\Lambda$, where
\begin{equation}\label{locIncShock}
\xi_1^0=\rho_1\sqrt{\frac{2(c_1^2-c_0^2)}{(\gamma-1)(\rho_1^2-\rho_0^2)}}=\frac{\rho_1u_1}{\rho_1-\rho_0},
\end{equation}
which is determined by the Rankine-Hugoniot conditions between states (0) and (1) on $\CS_0$.
Then {Problem \ref{ibvp-c}} for self-similar solutions becomes the following problem:

\begin{problem}[Boundary Value Problem]\label{bvp-c}
{\it Seek a
solution $\varphi$ of equation \eqref{2-1}--\eqref{1-o}
in the self-similar
domain $\Lambda$ with the slip boundary condition
$D\varphi\cdot\nnu|_{\partial\Lambda}=0$
and the asymptotic boundary condition at infinity{\rm :}
$$
\varphi\to\bar{\varphi}=
\begin{cases} \varphi_0 \qquad\mbox{for}\,\,\,
                         \xi_1>\xi_1^0, \xi_2>\xi_1 \tan\theta_{\rm w},\\
              \varphi_1 \qquad \mbox{for}\,\,\,
                          \xi_1<\xi_1^0, \;\xi_2>0,
\end{cases}
\qquad \mbox{when $|\xxi|\to \infty$,}
$$
where $\pSi_0=-\frac{1}{2}|\xxi|^2$ and $\pSi_1=-\frac{1}{2}|\xxi|^2+u_1(\xi_1-\xi^0_1)$.
}
\end{problem}

A weak solution of Problem \ref{bvp-c} is obtained by the following modification
of Definition \ref{def:weak solution}:
\eqref{def:weak solution-i3-Eqn} is now required to hold
for all $\zeta\in C^\infty_c(\bR^2)$.
As discussed in Remark \ref{conormCondRmk},
with such a choice of function $\zeta$, the integral identity \eqref{def:weak solution-i3-Eqn}
includes both equation \eqref{2-1} and the boundary condition of conormal form:
$\rho D\varphi\cdot\nnu=0$ on $\partial\Lambda$.
A weak solution is called entropy solution if it satisfies
the entropy condition:
density $\rho$ increases
in the pseudo-flow direction of $D\vphi|_{\Lambda^+\cap \CS}$ across any discontinuity curve ({\it i.e.}, shock).

\smallskip
Now we describe the more detailed structure of the regular reflection-diffraction configurations as shown
in Figs. \ref{figure: free boundary problems-1}--\ref{figure: free boundary problems-2}.
If a solution has one of the regular shock reflection-diffraction configurations,
and if its pseudo-potential
$\varphi$ is $C^1$ in the subregion $\widehat{\Omega}$
between the wedge and the reflected shock, then, at $P_0$,
it should satisfy both
the slip boundary condition on the wedge and the Rankine-Hugoniot conditions
with state $(1)$ across the flat shock $\mathcal{S}_1=\{\varphi_1=\varphi_2\}$, which passes
through point $P_0$ where the incident shock meets the wedge boundary.
Define the uniform state (2) with pseudo-potential $\varphi_2(\xxi)$ such
that
$$
\varphi_2(P_0)=\varphi(P_0), \qquad D\varphi_2(\PtIncW)=
\lim_{P\to P_0,\; P\in \widehat{\Omega}} D\varphi(P).
$$
Then the constant density $\rho_2$ of state (2) is equal to $\rho(|D\varphi|^2, \varphi)(\PtIncW)$
defined by (\ref{2-1}):
$$
\rho_2=\rho(|D\varphi_2|^2, \varphi_2)(\PtIncW).
$$
From the properties of $\varphi$ discussed above, it follows that  $D\varphi_2\cdot\nnu=0$ on the wedge boundary and the Rankine-Hugoniot conditions
\eqref{1-h}--\eqref{1-i} hold on the flat shock $\CS_1=\{\varphi_1=\varphi_2\}$
between states (1) and (2), which passes through $\PtIncW$.
In particular, $\varphi_2$ satisfies the following three conditions at $P_0$:
\begin{equation}\label{condState2}
D\varphi_2\cdot\bn_{\rm w}=0, \,\,\,\,\varphi_2=\varphi_1,\,\,\,\,\r(|D\varphi_2|^2,\varphi_2)D\varphi_2\cdot\bn_{\CS_1}=
\rho_1D\varphi_1\cdot\bn_{\CS_1} \qquad
\mbox{for $\bn_{\mathcal{S}_1}=\frac{D(\varphi_1-\varphi_2)}{|D(\varphi_1-\varphi_2)|}$}.
\end{equation}
where $\bn_{\rm w}$ is the outward normal to the wedge boundary.

\smallskip
The entropy solution $\varphi$, correspondingly state (2), can be either supersonic or subsonic at $\PtIncW$.
This determines the supersonic or subsonic type of regular shock reflection-diffraction configurations.
The regular reflection solution in the supersonic region is expected to
consist of the constant states separated by straight shocks ({\it cf.} \cite[Theorem 4.1]{Serre}).
Then, when state (2) is supersonic at $\PtIncW$, it can be shown that the constant state (2), extended
up to arc $P_1P_4$ of the sonic circle of state (2) between the wall and
the straight shock $P_0P_1\subset \CS_1$ separating it from state (1),
as shown in Fig. \ref{figure: free boundary problems-1},
satisfies equation (\ref{2-1}) in the region,
the Rankine-Hugoniot condition \eqref{1-h}--\eqref{1-i} on the straight shock $P_0P_1$,
and the slip boundary condition: $D\varphi_2\cdot\bn_{\rm w}=0$ on the
wedge $P_0P_4$,
and is expected to be a part of the regular
shock reflection-diffraction configuration.
Then the supersonic regular shock reflection-diffraction configuration
on Fig. \ref{figure: free boundary problems-1}
consists of three uniform states (0), (1), (2), and a non-uniform state in domain
$\Omega=P_1P_2P_3P_4$
where equation \eqref{2-1} is elliptic.
The reflected shock $\PtIncW\PtUpL\PtLwL$ has a straight part $\PtIncW\PtUpL$.
The elliptic domain $\Omega$ is separated from the hyperbolic region $\PtIncW\PtUpL\PtUpR$
of state (2) by the sonic arc $\PtUpL \PtUpR$ which lies on the sonic circle of state (2),
and the ellipticity in $\Omega$ degenerates on the sonic arc $P_1P_4$.
The subsonic regular shock reflection-diffraction configuration
as shown in Fig. \ref{figure: free boundary problems-2}
consists of two uniform states (0) and (1), and a non-uniform
state in domain
$\Omega=P_0P_2P_3$,
where the equation is elliptic, and
$\varphi_{|\Omega}(\PtIncW)=\varphi_2(\PtIncW)$ and
$D(\varphi_{|\Omega})(\PtIncW)=D\varphi_2(\PtIncW)$.

For the supersonic regular shock reflection-diffraction configurations
in Fig. \ref{figure: free boundary problems-1},
we use $\Sonic$, $\Shock$, $\Wedge$, and $\Symm$ for the sonic arc $P_1P_4$, the curved part of the reflected shock $P_1P_2$,
the wedge boundary $P_3P_4$, and the symmetry line segment $P_2P_3$, respectively.

For the subsonic regular shock reflection-diffraction configurations
in Fig. \ref{figure: free boundary problems-2},
$\Shock$, $\Wedge$, and $\Symm$ denote $P_0P_2$, $P_0P_3$, and $P_2P_3$, respectively.
We unify the notations with the supersonic reflection case by introducing
points $\PtUpL$ and $\PtUpR$ for the subsonic reflection case as
\begin{equation}\label{P1-P4-P0-subs-Ch2}
\PtUpL\defd\PtIncW, \quad \PtUpR\defd\PtIncW, \quad \overline\Sonic\defd\{\PtIncW\}.
\end{equation}

\begin{figure}
 \centering
\includegraphics[height=47mm]{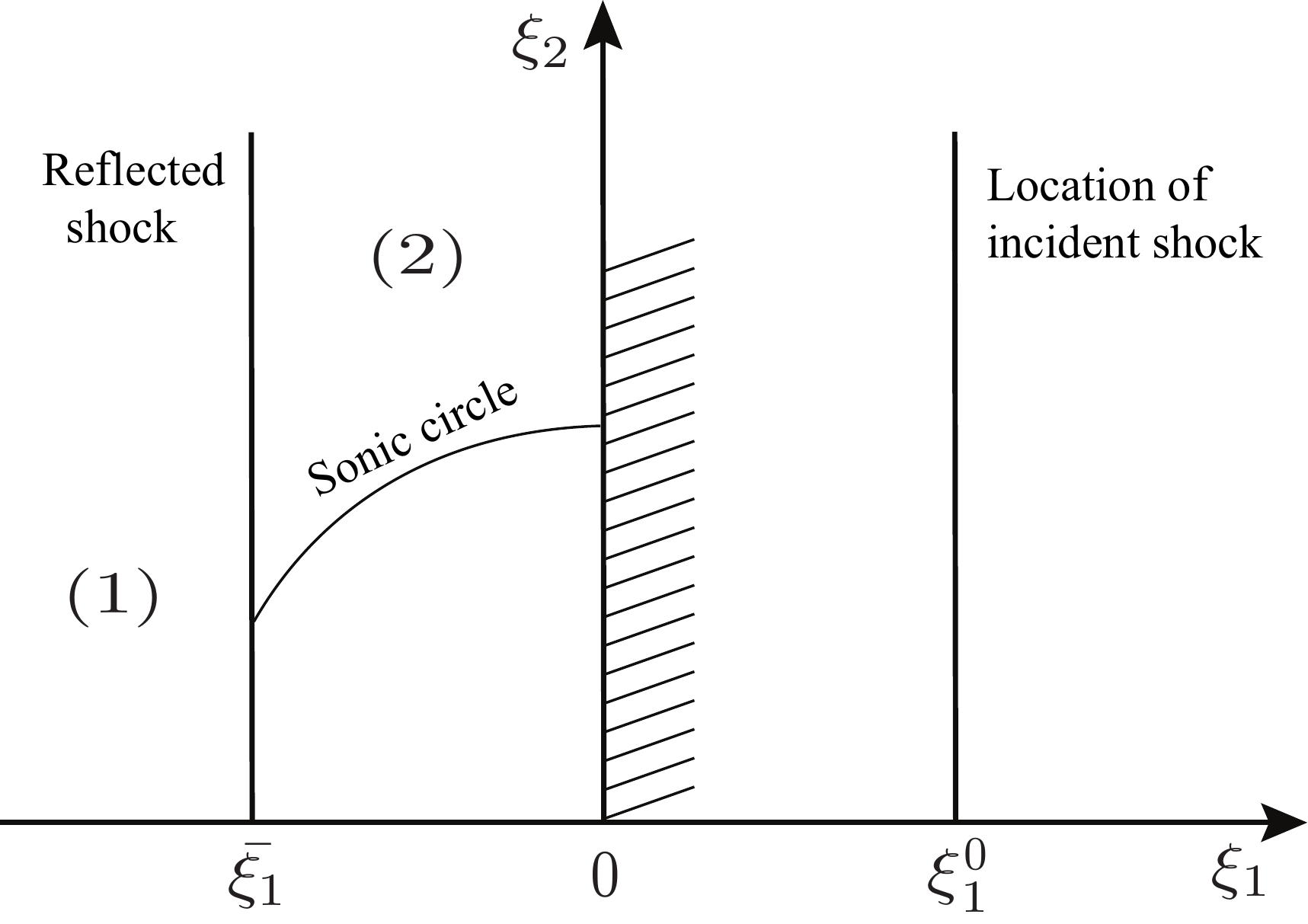}
\caption[]{Normal reflection configuration ({\it cf}. \cite{CF-book2018})}
\label{NormReflFigure}
\end{figure}

The corresponding solution for $\theta_{\rm w}=\frac\pi 2$ is called {\em normal reflection}.
In this case, the incident shock normally reflects from the flat wall; see Fig. \ref{NormReflFigure}.
The reflected shock is also a plane  $\{\xi_1=\bar{\xi}_1\}$, where $\bar{\xi}_1<0$.

From the discussion above, it follows that a necessary condition for the existence
of a regular reflection solution is the existence of the uniform state (2) with pseudo-potential $\varphi_2$ determined by
the boundary condition $D\varphi_2\cdot\nnu=0$ on the wedge  and the Rankine-Hugoniot
conditions \eqref{1-h}--\eqref{1-i}
across the flat shock $\CS_1=\{\varphi_1=\varphi_2\}$ separating it from state (1), and satisfying
the entropy condition: $\rho_2>\rho_1$.
These conditions lead to the system of algebraic equations \eqref{condState2}
for the constant velocity $\mathbf{u}_2$
and
density $\rho_2$ of state (2).
System \eqref{condState2} has solutions for some but not all of the wedge angles.
More specifically, for any fixed densities $0<\rho_0<\rho_1$ of states (0) and (1),
there exist a sonic angle $\theta_{\rm w}^{\rm s}$ and a detachment angle
$\theta_{\rm w}^{\rm d}$
satisfying
$$
0<\theta_{\rm w}^{\rm d}<\theta_{\rm w}^{\rm s}<\frac{\pi}{2}
$$
such that the algebraic system \eqref{condState2} has two solutions for
each $\theta_{\rm w}\in (\theta_{\rm w}^{\rm d}, \frac{\pi}{2})$, which become equal when
$\theta_{\rm w}=\theta_{\rm w}^{\rm d}$.
Thus, for each $\theta_{\rm w}\in (\theta_{\rm w}^{\rm d}, \frac{\pi}{2})$, there exist two
states (2), called weak and strong, with densities  $\rho_2^{\rm weak}<\rho_2^{\rm strong}$.
The weak state (2) is supersonic at the reflection point $\PtIncW(\theta_{\rm w})$ for
$\theta_{\rm w}\in (\theta_{\rm w}^{\rm s}, \frac{\pi}{2})$,
sonic for $\theta_{\rm w}=\theta_{\rm w}^{\rm s}$,
and subsonic for $\theta_{\rm w}\in (\theta_{\rm w}^{\rm d}, \hat\theta_{\rm w}^{\rm s})$
for some $\hat\theta_{\rm w}^{\rm s}\in(\theta_{\rm w}^{\rm d}, \theta_{\rm w}^{\rm s}]$.
The strong state (2) is subsonic at $\PtIncW(\theta_{\rm w})$ for all
$\theta_{\rm w}\in (\theta_{\rm w}^{\rm d}, \frac{\pi}{2})$.

There had been a long debate to
determine which of the two states (2) for
$\theta_{\rm w}\in (\theta_{\rm w}^{\rm d}, \frac{\pi}{2})$,
 weak or strong, is physical for the local theory; see
\cite{BD,CF-book2018,CF} and the references cited therein.
It was conjectured that the strong shock reflection-diffraction configuration
would be non-physical;
indeed, it is shown in Chen-Feldman \cite{ChenFeldman,CF-book2018}
that the weak shock reflection-diffraction configuration tends to the unique normal
reflection in Fig. \ref{NormReflFigure}, but the strong reflection-diffraction configuration does not,
when the wedge angle $\theta_{\rm w}$ tends to $\frac{\pi}{2}$.
The entropy condition and the definition of weak and strong states (2) imply that
$0<\rho_1<\rho_2^{\rm weak}<\rho_2^{\rm strong}$, which shows that the strength of the corresponding
reflected shock near $P_0$ in the
weak shock reflection-diffraction configuration
is relatively weak,
compared to the other shock given by the strong state (2).

\smallskip
If the weak state (2) is supersonic, the propagation speeds of the solution
are finite, and state (2) is
completely determined by the local information: state (1),
state (0), and the location of point $P_0$. That is, any information
from the reflection-diffraction region,
particularly the disturbance at corner $P_3$,
cannot travel towards the reflection point $P_0$.
However, if
it is subsonic, the information can reach $P_0$ and interact with
it, potentially
altering the subsonic reflection-diffraction configuration.
This argument motivated the following conjecture by
von Neumann in \cite{Neumann1,Neumann2}:

\medskip
{\bf The Sonic Conjecture}:
{\em There exists a supersonic regular shock reflection-diffraction
configuration when
$\theta_{\rm w}\in (\theta_{\rm w}^{\rm s}, \frac{\pi}{2})$
for $\theta_{\rm w}^{\rm s}>\theta_{\rm w}^{\rm d}$.
That is,
the supersonicity of the weak state {\rm (2)} implies the existence
of a supersonic regular reflection
solution, as shown in Fig. {\rm \ref{figure: free boundary problems-1}.}}

\medskip
Another conjecture is that the global regular shock reflection-diffraction
configuration is possible whenever the local regular reflection at the reflection
point is possible:

\medskip
{\bf The von Neumanm Detachment Conjecture}:
{\em There exists a regular shock reflection-diffraction configuration for
any wedge angle $\theta_{\rm w}\in (\theta_{\rm w}^{\rm d}, \frac{\pi}{2})$.
That is, the existence of state {\rm (2)} implies the existence
of a regular reflection solution,
as shown in Figs. {\rm \ref{figure: free boundary problems-1}--\ref{figure: free boundary problems-2}}.
}

\medskip
It is clear that the supersonic/subsonic regular shock reflection-diffraction configurations are
not possible without a local two-shock configuration at the
reflection point on the wedge, so the detachment conjecture is the weakest possible
criterion for the existence of supersonic/subsonic regular
shock reflection-diffraction configurations.

From now on, for the given wedge angle
$\theta_{\rm w}\in (\theta_{\rm w}^{\rm d}, \frac{\pi}{2})$,
state (2) represents the unique weak state (2)
and $\varphi_2$ is its pseudo-potential.
We now show how the solutions of regular shock reflection-diffraction configurations can be constructed.
This provides a solution to the von Neumann conjectures for potential flow.
Note that
state (2) is obtained from the algebraic conditions described above,
which determine
line $\CS_1$ and the sonic arc $P_1P_4$ when state (2) is supersonic at $P_0$,
and the slope of $\Shock$ at $P_0$ (arc $P_1P_4$ on the boundary of $\Omega$ becomes a corner point $P_0$)
when state (2) is subsonic at $P_0$.
Thus, the unknowns are domain $\Omega$ (or equivalently, the curved part of the reflected shock $\Shock$)
and the pseudo-potential $\varphi$ in $\Omega$.
Then, from \eqref{1-h}--\eqref{1-i}, in order to construct a solution of Problem \ref{bvp-c}
of the supersonic or subsonic regular shock reflection-diffraction configuration, it suffices
to solve the following problem:

\begin{problem}[Free Boundary Problem]\label{fbp-c}
{\it  For $\theta_{\rm w}\in (\theta_{\rm w}^{\rm d}, \frac{\pi}{2})$,
find a free boundary $($curved reflected shock$)$ $\Shock \subset \Lambda\cap \{\cxi<\xi_{1\PtUpL}\}$
{\rm (}$\Shock=P_1P_2$ on Fig. {\rm \ref{figure: free boundary problems-1}}
and $\Shock=P_0P_2$ on Fig. {\rm \ref{figure: free boundary problems-2}}$)$
and a function $\varphi$ defined in region
$\Omega$ as shown in Figs. {\rm \ref{figure: free boundary problems-1}}--{\rm \ref{figure: free boundary problems-2}}
such that
\begin{itemize}
\item[\rm (i)]
Equation \eqref{2-1} is satisfied in $\Omega$, and the equation is strictly elliptic for $\varphi$ in $\overline\Omega\setminus\overline\Sonic${\rm ,}
\item[\rm (ii)]
$\vphi=\vphi_1$ and $\rho D\vphi\cdot\nnu_{\rm s}=\rho_1 D\vphi_1\cdot\nnu_{\rm s}$ {on} the free boundary $\Shock${\rm ,}
\item[\rm (iii)]
$\vphi=\vphi_2$ and $D\vphi=D\vphi_2$ {on} $P_1P_4$
in the supersonic case as shown in Fig. {\rm \ref{figure: free boundary problems-1}}
 and at $P_0$ in the subsonic case as shown in Fig. {\rm \ref{figure: free boundary problems-1}}{\rm ,}
\item[\rm (iv)]
$D\vphi\cdot\nnu_{\rm w}=0$ {on} $\Wedge$, and
$D\vphi\cdot\nnu_{\rm sym}=0$ {on} $\Symm$,
\end{itemize}
where $\nnu_{\rm s}$, $\nnu_{\rm w}$, and $\nnu_{\rm sym}$ are the interior unit normals to $\Omega$
on $\shock$, $\Wedge$, and $\Symm$, respectively.
}
\end{problem}
Indeed, if $\varphi$ is a solution of Problem \ref{fbp-c}, we define its extension
from
$\Omega$ to $\Lambda$ by setting:
\begin{equation}\label{phi-states-0-1-2-MainThm}
\varphi=\begin{cases}
\, \varphi_0 \qquad\, \mbox{for}\,\, \xi_1>\xi_1^0 \mbox{ and } \xi_2>\xi_1\tan\theta_{\rm w},\\[1mm]
\, \varphi_1 \qquad\, \mbox{for}\,\, \xi_1<\xi_1^0
  \mbox{ and above curve} \,\, P_0\PtUpL\PtLwL,\\[1mm]
\, \varphi_2 \qquad\, \mbox{in region}\,\, P_0\PtUpL\PtUpR,
\end{cases}
\end{equation}
where we have used the notational convention {\rm (\ref{P1-P4-P0-subs-Ch2})}
for the subsonic reflection case, in which
region $P_0\PtUpL\PtUpR$ is one point
and curve $P_0\PtUpL\PtLwL$ is $P_0\PtLwL$;
see Figs. {\rm \ref{figure: free boundary problems-1}}--{\rm \ref{figure: free boundary problems-2}}.
Also, $\xi_1^0$ used in  (\ref{phi-states-0-1-2-MainThm}) is the location of the incident shock
({\textit{cf}.} (\ref{locIncShock})), and
the extension by (\ref{phi-states-0-1-2-MainThm})
is well-defined because of the requirement that $\Shock \subset \Lambda\cap \{\cxi<\xi_{1\PtUpL}\}$
in {Problem \ref{fbp-c}}.

Note that the conditions in  {Problem \ref{fbp-c}}(ii)
are the Rankine-Hugoniot conditions \eqref{1-h}--\eqref{1-i} on $\Shock$ between
$\varphi_{|\Omega}$ and $\varphi_1$.
Since $\Shock$ is a free boundary and equation  \eqref{2-1} is strictly elliptic for $\varphi$
in $\overline\Omega\setminus\overline\Sonic$,
then two conditions --- the Dirichlet and oblique derivative conditions --- on $\Shock$ are consistent with
one-phase free boundary problems for nonlinear elliptic PDEs of second order ({\it cf.} \cite{AC,ACF}).

In the supersonic case,
the conditions in {Problem \ref{fbp-c}}(iii) are the Rankine-Hugoniot conditions on $\Sonic$ between
$\varphi_{|\Omega}$ and $\varphi_2$. Indeed,
since state (2) is sonic on $\Sonic$, then it follows from
\eqref{1-h}--\eqref{1-i} that no gradient jump occurs on $\Sonic$.
Then, if $\varphi$ is a solution of {Problem \ref{fbp-c}},
its extension by \eqref{phi-states-0-1-2-MainThm} is a weak solution of {Problem \ref{bvp-c}}.
From now on, we consider a solution of {Problem \ref{fbp-c}} to be a function defined in $\Lambda$
by extension via \eqref{phi-states-0-1-2-MainThm}.

Since $\Sonic$ is not a free boundary (its location is fixed), it is not possible in general to prescribe
two conditions given in  {Problem \ref{fbp-c}}(iii) on $\Sonic$ for a second-order elliptic PDE.
In the iteration problem,
we prescribe the condition: $\varphi=\varphi_2$ on $\Sonic$, and then prove that
$D\varphi=D\varphi_2$ on $\Sonic$ by exploiting the elliptic degeneracy on $\Sonic$,
as we describe below.

\medskip
We observe that the key obstacle to prove the existence of regular shock reflection-diffraction configurations
as conjectured by von Neumann \cite{Neumann1,Neumann2} is an additional possibility that,
for some wedge angle $\theta_{\rm w}^{\rm a}\in (\theta_{\rm w}^{\rm d}, \frac{\pi}2)$, shock
$\PtIncW\PtLwL$ may attach to the wedge vertex $\PtLwR$, as observed
by experimental results ({\it cf.} \cite[Fig. 238]{VD}).
To describe the conditions of such an attachment, we note that
$$
\rho_1>\rho_0, \qquad
u_1=(\rho_1-\rho_0)
\sqrt{\frac{2(\rho_1^{\gamma-1}-\rho_0^{\gamma-1})}{\rho_1^2-\rho_0^2}},
\,\,\qquad c_1=\rho_1^{\frac{\gamma-1}2}.
$$
Then it follows from the explicit expressions above that, for each $\rho_0$,
there exists $\rho^{\rm c}>\rho_0$ such that
\begin{eqnarray*}
u_1\le c_1 \quad \mbox{if $\rho_1\in (\rho_0, \rho^{\rm c}]$}; \,\, \qquad
u_1>c_1 \quad \mbox{if $\rho_1\in (\rho^{\rm c}, \infty)$}.
\end{eqnarray*}

\smallskip
If $u_1\le c_1$, we can rule out the solution with a shock attached to the wedge vertex.
This is based on the fact that, if $u_1\le c_1$, then  the wedge vertex
$P_3=(0,0)$ lies within the sonic circle $\overline{B_{c_1}((u_1, 0))}$ of state (1),
and $\Shock$ does not intersect $\overline{B_{c_1}((u_1, 0))}$, as we show below.

\smallskip
If $u_1> c_1$, there would be a possibility that
the reflected shock could be attached to the wedge vertex
as the experiments show $(${\it e.g.}, \cite[Fig. 238]{VD}$)$.

\smallskip
Thus, in \cite{ChenFeldman,CF-book2018},
we have obtained the following results:

\smallskip
\begin{theorem}\label{mainShockReflThm}
There are two cases{\rm :}
\begin{enumerate}
\item[\rm (i)]
If  $\rho_0$ and $\rho_1$ are such that $u_1\le c_1$, then the
supersonic/subsonic regular reflection solution
exists for each wedge angle $\theta_{\rm w}\in (\theta_{\rm w}^{\rm d}, \frac{\pi}{2})$.
That is, for each $\theta_{\rm w}\in (\theta_{\rm w}^{\rm d}, \frac{\pi}{2})$,
there exists a solution $\varphi$ of {\rm Problem \ref{fbp-c}} such that
$$
\Phi(t, {\bf x}) =t\,\varphi(\frac{\bf x}{t}) +\frac{|{\bf x}|^2}{2t}
\qquad\mbox{for}\,\, \frac{\bf x}{t}\in \Lambda,\, t>0
$$
with
$$
\rho(t, {\bf x})=\Big(\rho_0^{\gamma-1}-(\gamma-1)\big(\Phi_t
      +\frac{1}{2}|\nabla_{\bf x}\Phi|^2\big)\Big)^{\frac{1}{\gamma-1}}
$$
is a global weak solution of {\rm Problem \ref{ibvp-c}} in the sense of
Definition {\rm \ref{weakSol-def-Prob1}} satisfying the entropy condition{\rm ;}
that is, $\Phi(t, {\bf x})$ is an entropy solution.

\smallskip
\item[\rm (ii)]
If  $\rho_0$ and $\rho_1$ are such that $u_1> c_1$, then there exists
$\theta_{\rm w}^{\rm a}\in [\theta_{\rm w}^{\rm d}, \frac{\pi}2)$ so that
the regular reflection solution exists for
each  wedge angle $\theta_{\rm w}\in (\theta_{\rm w}^{\rm a}, \frac{\pi}{2})$,
and the solution is of self-similar structure
described in {\rm (i)}  above.
Moreover, if $\theta_{\rm w}^{\rm a}>\theta_{\rm w}^{\rm d}$,
then, for the wedge angle $\theta_{\rm w}=\theta_{\rm w}^{\rm a}$,
there exists an {\it attached} solution, {\it i.e.},  $\varphi$ is
a solution of {\rm Problem \ref{fbp-c}} with  $\PtLwL=\PtLwR$.
\end{enumerate}
The type of regular shock reflection-diffraction
configurations $($supersonic as in Fig. {\rm \ref{figure: free boundary problems-1}}
or  subsonic as in Fig. {\rm \ref{figure: free boundary problems-2}}$)$
is determined by the type of state
{\rm (2)} at $\PtIncW${\rm :}
\begin{enumerate}[\rm (a)]
\item
For the supersonic and sonic reflection case,
the reflected shock $\PtIncW\PtLwL$ is $C^{2,\alpha}$--smooth for some $\alpha\in(0,1)$ and
its curved part $P_1P_2$ is $C^\infty$ away from $P_1$.
The solution $\varphi$ is in $C^{1,\alpha}(\overline\Omega)\cap C^\infty(\Omega)$, and is
$C^{1,1}$ across the sonic arc which is optimal{\rm ;} that is, $\varphi$ is {\em not} $C^2$ across sonic arc.
\item
For the subsonic reflection case
$($Fig. {\rm \ref{figure: free boundary problems-2}}$)$,
the reflected shock $\PtIncW\PtLwL$
and solution $\varphi$ in $\Omega$ is in $C^{1,\alpha}$
near $P_0$ and $P_3$ for some $\alpha\in(0,1)$, and $C^\infty$ away from $\{P_0,P_3\}$.
\end{enumerate}
Moreover, the regular reflection solution tends to the unique normal
reflection $($as in Fig. {\rm \ref{NormReflFigure})}
when the wedge angle $\theta_{\rm w}$ tends to $\frac{\pi}{2}$.
In addition, for both supersonic and subsonic reflection cases,
\begin{equation}\label{phi-between-in-omega}
 \varphi_2<\varphi<\varphi_1 \qquad\mbox{in $\Omega$}.
\end{equation}
Furthermore, $\varphi$ is an admissible solution in the sense of
Definition {\rm \ref{admisSolnDef}} below, so that $\varphi$ satisfies
further properties listed in
Definition {\rm \ref{admisSolnDef}}.
\end{theorem}

\medskip
Theorem \ref{mainShockReflThm} is proved by solving {Problem \ref{fbp-c}}.
The first results on the existence of global solutions of the free boundary problem ({Problem \ref{fbp-c}})
were obtained for the wedge angles sufficiently close to $\frac \pi 2$ in Chen-Feldman \cite{ChenFeldman}.
Later, in Chen-Feldman \cite{CF-book2018}, these results
were extended up to the detachment angle
as stated in Theorem \ref{mainShockReflThm}. For this extension, the
techniques developed in  \cite{ChenFeldman}, notably the estimates
near the sonic arc, were the starting point.

\medskip
{\bf Case I: The wedge angles close to $\frac \pi 2$}.
Let us first discuss the techniques in \cite{ChenFeldman}, where
we employ the approach of Chen-Feldman \cite{CF-JAMS2003} to develop an iteration scheme for
constructing a global solution of {Problem \ref{fbp-c}},
when the wedge angle $\theta_{\rm w}$ is close to $\frac \pi 2$.
For this case, the solutions are of the supersonic regular shock reflection-diffraction configuration
as in Fig. \ref{figure: free boundary problems-1}.
The general procedure is similar to the one described in \S\ref{existenceSection}, which
can be presented in the following four steps:

\medskip
1. Fix $\theta_{\rm w}$ sufficiently close to $\frac\pi 2$ so that various constants
in the argument can be controlled.
The iteration set consists of functions defined on a region ${\mathcal D}$, where ${\mathcal D}$
contains all possible $\Omega$ for the fixed $\theta_{\rm w}$.
Specifically, an important property of the regular shock reflection-diffraction configurations
is \eqref{phi-between-in-omega}, which implies that $\Omega\subset \{\varphi_2<\varphi_1\}$;
that is, $\Omega$ lies {\it below} line $\CS_1$ passing through $P_0$ and $P_1$ on
Fig. \ref{figure: free boundary problems-1}.
Note that, when $\theta_{\rm w}$ is close to $\frac\pi 2$, this line is close to the vertical reflected shock
of normal reflection on Fig. \ref{NormReflFigure}.
Then ${\mathcal D}$ is defined as a region bounded
by $\CS_1$, $\Sonic=P_1P_4$, $\Wedge=P_3P_4$, and the symmetry line $\xi_2=0$.
The iteration set is a set of functions $\varphi$ on ${\mathcal D}$,
defined by $\varphi\ge\varphi_2$ on ${\mathcal D}$ and the bound of norm of $\varphi-\varphi_2$ on ${\mathcal D}$
in the scaled and weighted $C^{2,\alpha}$ space
defined in \eqref{parabNormsApp} below.
Such functions satisfy
$$
\|\varphi-\varphi_2\|_{C^{1,\alpha}(\overline{\mathcal D})}\le C(\frac\pi 2-\theta_{\rm w}),
$$
which is small when $\frac\pi 2-\theta_w\ll 1$,
and
$$
\|\varphi-\varphi_2\|_{C^{1,1}(\overline{\mathcal D}\cap{\mathcal N}_\varepsilon(\Sonic))}<\infty.
$$
However, $\|\varphi-\varphi_2\|_{C^{1,1}(\overline{\mathcal D}\cap{\mathcal N}_\varepsilon(\Sonic))}$ is not small
even if $\frac\pi 2-\theta_{\rm w}$ is small; the reasons for that will be discussed below.

Given a function $\hat\varphi$ from the iteration set,
we define domain $\Omega(\hat \varphi):=\{\hat\varphi<\varphi_1\}$
so that the iteration free boundary is
$\Shock(\hat\varphi)=\partial \Omega(\hat\varphi)\cap {\mathcal D}$.
This is similar to \eqref{OmegaPL-f},
and the corresponding non-degeneracy similar to \eqref{nondegeneracy} in the present case is:
$$
\partial_{\xi_1}(\varphi_1-\varphi_2-\phi)\ge \frac{u_1}{2}
\qquad\,\,\mbox{in ${\mathcal D}$ \, if $\|\phi\|_{C^1(\overline{{\mathcal D}})}$ and $\frac\pi 2-\theta_{\rm w}$ are small}.
$$
Then we define the iteration equation by using form \eqref{equ:study} of equation  \eqref{2-1},
by making an elliptic truncation
(which is somewhat different from Step 1 in \S\ref{existenceSection})
and substituting $\hat\varphi$ in some terms of the coefficients of \eqref{equ:study}.
The iteration boundary condition on $\Shock(\hat\varphi)$ is an oblique derivative condition
obtained by combining two conditions in {Problem \ref{fbp-c}}(ii) and making some truncations.
On $\Sonic$, we prescribe $\varphi=\varphi_2$, {\it i.e.}, one of two conditions in {Problem \ref{fbp-c}}(iii).
On $\Wedge$ and $\Symm(\hat\varphi)$, we prescribe the conditions given in {Problem \ref{fbp-c}}(iv).
The iteration map: $\hat\varphi\to \varphi$ is defined
by solving the iteration problem to obtain $\varphi$ and then extending $\varphi$ from
$\Omega(\hat \varphi)$ to ${\mathcal D}$.

The fundamental differences between the iteration procedure in the shock reflection-diffraction problem
and the previous procedures  on transonic shocks in the steady
case in \S2--\S3 (such as \cite{CF-JAMS2003,ChenFeldman2,ChenFeldman3Arch,XY} and follow-up papers) include:
\begin{enumerate}[{\rm (i)}]
\item \label{numRegRefl-1}
The procedures on steady transonic shocks in \S 2--\S3 are for the perturbation case.
In particular, the ellipticity of the iteration equation and the removal of the elliptic cutoff are achieved
by making the iteration set sufficiently close to the background solution in $C^1$ or a stronger norm.
For the regular reflection problem, this cannot be done because of the elliptic degeneracy near the sonic arc.

\item  \label{numRegRefl-2}
Only one condition  on $\Sonic$ can be prescribed; however, both
 $\varphi=\varphi_2$ and $D\varphi=D\varphi_2$ on $\Sonic$ are needed to be matched to obtain a global entropy solution.
This is resolved by exploiting the elliptic degeneracy on $\Sonic$.
\end{enumerate}

\smallskip
2. In order to see the elliptic degeneracy on $\Sonic$ more explicitly,
we fix the wedge angle $\theta_{\rm w}$ and the corresponding pseudo-potential
$\varphi_2=\varphi_2^{(\theta_{\rm w})}$ of the weak state (2), and rewrite
equation \eqref{nondivMainEq} in terms of the function:
$$
\psi=\varphi-\varphi_2
$$
in the following coordinates flattening $\Sonic$:
\begin{equation}\label{x-y}
x=c_2-r,\qquad  y=\theta-\theta_{\rm w},
\end{equation}
where $(r, \theta)$ are the polar coordinates centered at
$O_2=\mathbf{u}_2$
of the sonic circle of state (2).
Then
$$
\mbox{$\Omega_\varepsilon:=\Omega\cap {\mathcal N}_\varepsilon(\Sonic)\subset \{x>0\}\,\,$ for small $\varepsilon>0$,
$\qquad\,\,$ $\Sonic\subset \{x=0\}$.}
$$
In what follows, we always assume that $\varphi\in C^{1,1}(\overline\Omega_\varepsilon)$ as
in Theorem \ref{mainShockReflThm} for the supersonic case.
Then, by the conditions in Problem \ref{fbp-c}(iii)  and the definition of $\psi$,
\begin{align}
&\psi=0  \qquad\quad\,\mbox{on $\Sonic$},\label{psiZeroOnSonic}\\
&D\psi=0 \qquad\,\,\mbox{on $\Sonic$}. \label{D-psiZeroOnSonic}
\end{align}
Moreover, we {\it a priori} assume that
solution $\varphi$ satisfies \eqref{phi-between-in-omega} in $\Omega$
to derive the required estimates of the solution;
with these estimates, we then construct such a solution and verify that it satisfies \eqref{phi-between-in-omega}.
The heuristic motivation of \eqref{phi-between-in-omega} is the following:
From Figs. \ref{figure: free boundary problems-1}--\ref{figure: free boundary problems-2},
it appears that $\Shock$ (and hence $\Omega$) is located {\it below} line $\CS_1$, {\it i.e.},
in the half-plane $\{\varphi_1>\varphi_2\}$.
Thus, $\varphi=\varphi_1>\varphi_2$ on $\Shock$, and
$\varphi_1>\varphi_2=\varphi$ on $\Sonic$.
Also,
the potential functions
$\phi_1$ and $\phi_2$
of states (1) and (2)
are linear functions,  thus they
satisfy equation \eqref{equ:study} with coefficients determined by $\varphi$,
considered as a linear equation for $\phi$.
Taking into account the inequalities on
$\Shock$ and $\Sonic$ noted above, and the oblique boundary conditions on $\Wedge$ and $\Symm$,
we obtain \eqref{phi-between-in-omega} by the maximum principle.
Then, from \eqref{phi-between-in-omega}, we have
\begin{equation}\label{psiPositiveInOmega}
 \psi>0\qquad\mbox{in $\Omega$}.
\end{equation}
Even though the previous argument is heuristic, the fact that it comes from the structure of the problem
allows us to include the condition that $\psi\ge 0$ in the definition of the iteration set
and close the iteration argument for constructing the solutions within this set.

\medskip
Equation \eqref{nondivMainEq} in $\Omega\cap {\mathcal N}_\varepsilon(\Sonic)$ for $\psi$
in the $(x,y)$--coordinates \eqref{x-y} is
\begin{equation}
\big(2x-(\gamma+1)\psi_x+O_1 \big)\psi_{xx} +O_2\psi_{xy} + ({1\over
c_2}+O_3 )\psi_{yy} -(1+O_4)\psi_{x} +O_5\psi_{y}=0,
\label{equationForPsi-sonicStruct}
\end{equation}
where
\begin{equation}\label{erTerms-xy-nontrunc}
\begin{split}
O_1(D\psi,\psi,x) &= -\frac{x^2}{c_2}+{\gamma+1\over
2c_2}(2x-\psi_x)\psi_x -{\gamma-1\over c_2}\Big(\psi+{1\over
2(c_2-x)^2}\psi_y^2\Big),
\\
O_2(D\psi,\psi, x)&=-\frac{2(\psi_x+c_2-x)\psi_y}{c_2(c_2-x)^2},
\\
O_3(D\psi,\psi, x) &={1\over c_2(c_2-x)^2}\Big(x(2c_2-x)-
(\gamma-1)\big(\psi+(c_2-x)\psi_x+{1\over 2}\psi_x^2\big)
-\frac{\gamma+1}{2(c_2-x)^2}\psi_y^2\Big),
\\
O_4(D\psi,\psi, x)
&=\frac{1}{c_2-x}\Big(x- {\gamma-1\over c_2}\big(\psi+(c_2-x)\psi_x+{1\over 2}\psi_x^2
  +\frac{(\gamma+1)\psi_y^2}{2(\gamma-1)(c_2-x)^2}\big)\Big),
\\
O_5(D\psi,\psi, x)&=
-\frac{2(\psi_x+c_2-x)\psi_y}{c_2(c_2-x)^3}.
\end{split}
\end{equation}
Since $\psi\in  C^{1,1}(\overline\Omega_\varepsilon)$,
it follows from \eqref{psiZeroOnSonic}--\eqref{D-psiZeroOnSonic}
that $|\psi(x,y)|\le Cx^2$ and
\begin{equation}\label{gradEstCx}
 | D\psi(x,y)|\le Cx\qquad\mbox{in $\Omega_\varepsilon$},
\end{equation}
so that
\begin{equation}\label{Ok-termEst}
 | O_1(D\psi, \psi, x)|\le  N|x|^2, \quad |
O_k(D\psi, \psi, x)|\le N|x|\qquad\,\,\, \mbox{for $k=2,\dots,5$}.
\end{equation}
Using \eqref{Ok-termEst}, we can show that
$O_k(D\psi,\psi, x)$ are small perturbations of the
leading terms of equation (\ref{equationForPsi-sonicStruct})
in $\Omega_\varepsilon=  \Omega\cap {\mathcal N}_\varepsilon(\Sonic)$.
Also, if \eqref{gradEstCx} holds, equation
\eqref{equationForPsi-sonicStruct} is strictly elliptic in $\overline\Omega_\varepsilon\setminus\overline{\Sonic}$ if
\begin{equation}\label{elliptCond}
\psi_x(x,y)\le \frac {2\mu}{\gamma+1}x
\end{equation}
for $\mu\in (0, 1)$, when $\varepsilon=\varepsilon(\mu, N)$ is small.
For $\theta_{\rm w}$ close to $\frac\pi 2$,
it can be shown that any solution of {Problem \ref{fbp-c}} (with some natural regularity properties)
satisfies that, for any small $\delta>0$,
\begin{equation}\label{elliptCond-actual}
|\psi_x(x,y)|\le \frac {1+\delta}{\gamma+1}x
\qquad\mbox{in $\Omega_\varepsilon$  for small $\varepsilon=\varepsilon(\delta)$,}
\end{equation}
which verifies \eqref{elliptCond} with any $\mu\in (\frac 12, 1)$ ({\it e.g.}, with $\mu=\frac 23$)
if $\delta$ is correspondingly small.

\smallskip
3. The iteration equation near $\Sonic$ is defined based on the above facts.
The iteration set $\setK$ used in \cite{ChenFeldman} is such that every
$\hat\psi= \hat\varphi-\varphi_2\in\setK$ satisfies \eqref{psiZeroOnSonic} and \eqref{gradEstCx}
for some $C, \varepsilon>0$.
Then the iteration equation for $\psi$ is
 \begin{equation}
\big(2x-(\gamma+1)x\eta(\frac{\psi_x}x)+ O_1^{(\hat\psi)} \big)\psi_{xx} +O_2^{(\hat\psi)}\psi_{xy} + ({1\over
c_2}+O_3^{(\hat\psi)} )\psi_{yy} -(1+O_4^{(\hat\psi)})\psi_{x} +O_5^{(\hat\psi)}\psi_{y}=0,
\label{equationForPsi-sonicStruct-Iter}
\end{equation}
where the cutoff function $\eta\in C^\infty(\bR)$ satisfies
$|\eta|\le \frac 5{3(\gamma+1)}$, $\eta'\ge 0$,
and $\eta(s)=s$ if $|s|\le \frac 4{3(\gamma+1)}$, and some other technical conditions.
The terms $O_k^{(\hat\psi)}$, $k=1,\dots, 5$, are obtained from $O_k$
by substituting $\hat\psi$ into certain terms in \eqref{erTerms-xy-nontrunc}
and performing the
cutoff in the remaining terms, so that estimates
\eqref{Ok-termEst} hold.
Then \eqref{equationForPsi-sonicStruct-Iter} is strictly
elliptic in $\overline\Omega_\varepsilon\setminus\overline{\sonic}$ for small $\varepsilon$,
and its ellipticity degenerates on $\sonic$.
Since the solution of  {Problem \ref{fbp-c}}  satisfies equation
\eqref{equationForPsi-sonicStruct}
and inequality \eqref{elliptCond-actual} with $\delta=\frac{1}{3}$
in $\Omega_\varepsilon$ for small $\varepsilon$,
then it satisfies equation
\eqref{equationForPsi-sonicStruct-Iter} in $\Omega_\varepsilon$ with $\hat\psi=\psi$.
Indeed, we have the estimate: $|\psi_x|\le \frac 4{3(\gamma+1)}x$, so that $x\eta(\frac{\psi_x}x)=\psi_x$;
and the cutoffs in the terms of $O_k^{(\hat\psi)}$ are removed similarly.

We also note that the degenerate ellipticity structure of equation \eqref{equationForPsi-sonicStruct-Iter} is the following:
Writing \eqref{equationForPsi-sonicStruct-Iter} in the form
\begin{align}\label{potFlow-xy-eqn}
&\sum_{i,j=1}^2  A_{ij}(D\psi, \psi, x)D_{ij}\psi
+\sum_{i=1}^2  A_i(D\psi, \psi, x)D_{i}\psi=0
\end{align}
with $A_{12}=A_{21}$, we see that,  for any $\xxi=(\xi_1,\xi_2)\in \bR^2$,
\begin{equation}\label{ellipt-potFlow-xy-eqn}
 \lambda |\xxi|^2\le
A_{11}({\bf p}, z, x)\frac{\xi_1^2}{x}
+2 A_{12}({\bf p}, z, x)\frac{\xi_1\xi_2}{x^{1/2}}
+A_{22}({\bf p}, z, x)\xi_2^2
\le \frac{1}{\lambda} |\xxi|^2
\end{equation}
for all $({\bf p}, z)\in \bR^2\times\bR$ and $\x\in (0, \varepsilon)$.

\smallskip
We consider solutions of \eqref{equationForPsi-sonicStruct-Iter} in $\Omega_\varepsilon$
satisfying \eqref{psiZeroOnSonic} and \eqref{psiPositiveInOmega}.
Since condition \eqref{D-psiZeroOnSonic} can not be prescribed
in the iteration problem as discussed above, we have to obtain
\eqref{D-psiZeroOnSonic} from the estimates of the solutions by exploiting the elliptic degeneracy.
The estimates of the positive solutions of \eqref{equationForPsi-sonicStruct-Iter} with
\eqref{psiZeroOnSonic}
in $\Omega_\varepsilon$ are based on the fact that, for any $\delta>0$, the function:
$$
w_\delta(x,y)=\frac{1+\delta}{2(\gamma+1)} x^2
$$
is a supersolution of \eqref{equationForPsi-sonicStruct-Iter} in $\Omega_\varepsilon$
if $\varepsilon=\varepsilon(\delta)$ is small; that is, ${\mathcal N}(w_\delta)<0$
in $\Omega_\varepsilon$,
where ${\mathcal N}(\cdot)$ denotes the operator determined by the left-hand side of
\eqref{equationForPsi-sonicStruct-Iter}.
Using this, the boundary conditions on $\Shock$ and $\Wedge$,
and \eqref{psiPositiveInOmega}, we obtain by the comparison principle that
\begin{equation}\label{quadrGrowthNearSonic}
0\le\psi\le Cx^2\qquad\,\,\mbox{in $\Omega_\varepsilon$},
\end{equation}
where $\varepsilon$ and $C$ are uniform with respect to the wedge angles near $\frac\pi 2$.
Note that $-w_\delta$ is {\em not} a subsolution of \eqref{equationForPsi-sonicStruct} so that it cannot be used
to bound $\psi$ from below.
Thus, property \eqref{psiPositiveInOmega}, which is derived from the global structure of the solution,
is crucially used in this argument.
Then, in
\eqref{quadrGrowthNearSonic}, the upper bound is from the local estimates near $\Sonic$,
while the lower bound is from the global structure of the problem.

In particular, \eqref{quadrGrowthNearSonic}
implies
that $D\psi=0$ on $\Sonic$, which resolves
the issue described in
\eqref{numRegRefl-2} above.
Furthermore, from \eqref{quadrGrowthNearSonic}, using the non-isotropic
{\it parabolic} rescaling corresponding to the elliptic degeneracy
\eqref{ellipt-potFlow-xy-eqn} of equation \eqref{equationForPsi-sonicStruct-Iter} near $x=0$,
we obtain the estimates in the appropriately weighted and scaled H\"{o}lder norm in $\Omega_\varepsilon$,
which also imply
the uniform $C^{1,1}$ estimates:
\begin{equation}\label{secodDerivBdParab}
|D^2\psi|\le C \qquad\mbox{in $\Omega_\varepsilon$}.
\end{equation}

More precisely, we denote  this norm by
$\|\psi\|_{2,\alpha,{\Omega_\varepsilon}}^{\rm (par)}$ and define it as follows: Denote
$\zz=(x,y)$ and $\tilde \zz=(\tilde x,\tilde y)$ with $x, \tilde x\in(0,
2\varepsilon)$ and
$$
\delta^{\rm (par)}_\alpha(\zz, \tilde \zz)\defd \left(|x-\tilde x|^2+
\max(x,\tilde x)|y-\tilde y|^2\right)^{\alpha/2}.
$$
Then, for $\psi\in C^2({\Omega_\varepsilon})\cap C^{1,1}(\overline{\Omega_\varepsilon})$ written
in the $(x,y)$--coordinates, we define
\begin{equation}\label{parabNormsApp}
\begin{split}
&\qquad\,\, \|\psi\|^{\rm (par)}_{2,0,{\Omega_\varepsilon}} :=\sum_{0\le k+l\le 2}
\sup_{\zz\in{\Omega_\varepsilon}}\left(x^{k+l/2-2}|\partial_x^k\partial_y^l\psi(\zz)|\right),\\
&\qquad\,\, [\psi]^{\rm (par)}_{2,\alpha,{\Omega_\varepsilon}} :=\sum_{k+l=2}\sup_{\zz,
\tilde \zz\in{\Omega_\varepsilon}, \zz\ne \tilde \zz}
 \bigg(\min(x^{k+l/2-2}, \tilde x^{k+l/2-2})
 \frac{|\partial_x^k\partial_y^l\psi(\zz)-\partial_x^k\partial_y^l\psi(\tilde \zz)|}
 {\delta^{\rm (par)}_\alpha(\zz,\tilde \zz)}\bigg),
  \\
&\qquad\,\,\|\psi\|^{\rm (par)}_{2,\alpha,{\Omega_\varepsilon}}
:=\|\psi\|^{\rm (par)}_{2,0,{\Omega_\varepsilon}} +[\psi]^{\rm (par)}_{2,\alpha,{\Omega_\varepsilon}}.
\end{split}
\end{equation}

Now we obtain the required estimates in the norm in \eqref{parabNormsApp},
under the assumption that \eqref{quadrGrowthNearSonic} holds in $\Omega_{2\varepsilon}$.
For every $\zz_0=(x_0, y_0)\in \overline\Omega_\varepsilon\setminus\Sonic$ (so that $x_0\in (0, \varepsilon]$),
we define
\begin{equation}\label{parabRectangles}
R_{\zz_0}=\Big\{(x,y)\;\;:\;\; |x-x_0|<\frac{x_0}{10},\;
|y-y_0|<\frac{\sqrt{x_0}}{10}\Big\}\cap \Omega.
\end{equation}
Note that $\dist(R_{\zz_0}, \sonic)=\frac{9}{10}x_0>0$.
We rescale the rectangle in \eqref{parabRectangles} to the unit square
$Q_1=(-1,1)^2$:
\begin{equation}\label{rescaled-parabRectangles}
Q_1^{(\zz_0)}:=\Big\{(S,T)\in Q_1\; : \; (x_0+\frac{x_0}{10}S,\;
y_0+\frac{\sqrt{x_0}}{10}T)\in \Omega\Big\},
\end{equation}
and define the scaled version of $\psi$ in the $(S,T)$--coordinates in $Q_1^{(\zz_0)}$:
\begin{equation}\label{parabRescaling}
\psi^{(\zz_0)}(S, T):=\frac{1}{x_0^2}\psi(x_0+\frac{x_0}{10}S,\; y_0+\frac{\sqrt{x_0}}{10}T)
\qquad\mbox{for $(S, T)\in Q_1^{(z_0)}$}.
\end{equation}
Note that this rescaling is non-isotropic with respect to the two variables $x$ and $y$.
By \eqref{quadrGrowthNearSonic}, we have
\begin{equation}\label{parabRescalingL-infty}
\|\psi^{(\zz_0)}\|_{L^\infty(\overline{Q_1^{(\zz_0)}})}\le C\qquad
\text{for any $\zz_0=(x_0, y_0)\in \overline\Omega_\eps\setminus\overline{\Sonic}$}.
\end{equation}
Rewriting equation \eqref{equationForPsi-sonicStruct-Iter} in terms of
$\psi^{(\zz_0)}$ in the $(S, T)$--coordinates and noting the degenerate ellipticity
structure \eqref{ellipt-potFlow-xy-eqn},
we find that $\psi^{(\zz_0)}$ satisfies a
uniformly elliptic equation in $Q_1^{(\zz_0)}$ with the ellipticity constants and
certain H\"{o}lder norms of the coefficients independent of $\zz_0$.
We also rescale  the boundary conditions on $\Shock\cap\partial\Omega_\varepsilon$
and $\Wedge\cap\partial\Omega_\varepsilon$ in a
similar way, when $\zz_0$ is on the corresponding part of the boundary.
Then we apply the local elliptic $C^{2,\alpha}$--estimates for $\psi^{(\zz_0)}$
in $Q_1^{(\zz_0)}$
in the following cases:
\begin{enumerate}
\item[\rm (i)] Interior rectangles $R_{\zz_0}$, {\it i.e.}, all $\zz_0$ such that $Q_1^{(\zz_0)}=Q_1$ holds,

\smallskip
\item[\rm (ii)] Rectangles  $R_{\zz_0}$ centered on the shock: $\zz_0\in\Shock\cap\partial\Omega_\varepsilon$,

\smallskip
\item[\rm (iii)] Rectangles  $R_{\zz_0}$ centered on the wedge: $\zz_0\in\Wedge\cap\partial\Omega_\varepsilon$,
\end{enumerate}
where, in the last two cases, we use the local estimates for the corresponding boundary value problems.
Using \eqref{parabRescalingL-infty}, we obtain
$$
\|\psi^{(\zz_0)}\|_{C^{2,\alpha}(\overline{Q_{1/2}^{(\zz_0)}})}\le C \qquad\,\,\mbox{with $C$ independent of $\zz_0$},
$$
where $Q_{1/2}^{(\zz_0)}=Q_{1}^{(\zz_0)}\cap (-\frac{1}{2},\;\frac{1}{2})^2$.
Rewriting in terms of $\psi$ in the $(x,y)$--coordinates and combining the estimates for all $\zz_0$ as above,
we obtain the estimate:
$\|\psi\|_{2,\alpha,{\Omega_\varepsilon}}^{\rm (par)}\le C$
in norm \eqref{parabNormsApp},
which also implies the $C^{1,1}$--estimates \eqref{secodDerivBdParab}.

\begin{remark}\label{est2ndDer-rmk}
 Note that $\psi^{(\zz_0)}_{SS}(S, T)=\frac{1}{100}\psi_{xx}(x_0+\frac{x_0}{10}S,\; y_0+\frac{\sqrt{x_0}}{10}T)$. It follows that
  $\|D^2\psi\|_{L^\infty}$ cannot be made small by choosing the
parameters, e.g., choosing $\varepsilon$ small or $\theta_{\rm w}$ close to $\frac\pi 2$.
\end{remark}

\begin{remark}\label{apriri-est-positiveSolRmk}
The above argument, starting from \eqref{parabRectangles}, is also used for the {\it a priori} estimates
of the positive solutions of \eqref{equationForPsi-sonicStruct}--\eqref{erTerms-xy-nontrunc}
 with condition \eqref{psiZeroOnSonic},
satisfying \eqref{gradEstCx} and the ellipticity condition \eqref{elliptCond} with some $\mu\in (0, 1)$.
Note that \eqref{psiZeroOnSonic}, \eqref{gradEstCx},
and $\psi\ge 0$  imply \eqref{quadrGrowthNearSonic},
which is used in the argument.
\end{remark}

\begin{remark}\label{positiveSolRmk}
Remark {\rm \ref{apriri-est-positiveSolRmk}} applies only to the positive solutions
of  \eqref{equationForPsi-sonicStruct} with condition \eqref{psiZeroOnSonic}.
For the negative solutions of \eqref{equationForPsi-sonicStruct} with condition \eqref{psiZeroOnSonic},
the equation is {\em uniformly} elliptic up to $\{x=0\}$ and,
similar to Hopf's lemma, the negative solutions have {\em linear} growth{\rm :} $|\psi(x,y)|\ge \frac 1C x$,
in a contrast with \eqref{quadrGrowthNearSonic}.
This feature is used in the proof of certain geometric properties of the free boundary for the wedge angles
away from $\frac\pi 2$, where we note that $\varphi-\varphi_1<0$ by \eqref{phi-between-in-omega}.
\end{remark}

\smallskip
{4}. In order to remove the ellipticity cutoff in \eqref{equationForPsi-sonicStruct-Iter},
{\it i.e.}, to show that the fixed point solution of
\eqref{equationForPsi-sonicStruct-Iter} ({\it i.e.}, with   $\psi=\hat\psi$) actually satisfies
\eqref{equationForPsi-sonicStruct},
we need to show that
$|\psi_x|\le \frac 4{3(\gamma+1)}x$, as we have discussed right after \eqref{equationForPsi-sonicStruct-Iter}.
Combining \eqref{secodDerivBdParab} with $D\psi=0$ on $\Sonic$, we obtain that
$|D\psi(x,y)|\le Cx$ in $\Omega_\varepsilon$, which does not remove the ellipticity cutoff,
unless we show the explicit bound $C\le\frac 4{3(\gamma+1)}$.
However, this bound does not follow from the estimates discussed above ({\it cf.} Remark \ref{est2ndDer-rmk}).

Note that the
only explicit solution we know is the normal reflection for $\theta_{\rm w}=\frac \pi 2$, for which
$\varphi=\varphi_2^{(\frac\pi 2)}$, {\it i.e.}, $\psi=0$ in $\Omega$.
Also, the analysis in Bae-Chen-Feldman \cite{BaeChenFeldman} has shown that the solutions of
{\rm Problem \ref{fbp-c}} for the supersonic regular shock reflection-diffraction configuration
satisfy that, for small $\varepsilon$,
$$
\psi_x\thicksim \frac x{\gamma+1} \qquad \mbox{in $\Omega_\varepsilon\cap\{(x,y)\;:\; \mbox{dist}\,((x,y), \Shock)>\sqrt{x}\}$,}
$$
but
$$
D\psi=o(x) \qquad \mbox{in $\Omega_\varepsilon\cap\{(x,y)\;:\; \mbox{dist}\,((x,y), \Shock)<x^2\}$}.
$$
This shows that the convergence of
solutions $\varphi^{(\theta_{\rm w})}$ of {\rm Problem \ref{fbp-c}}
to $\varphi^{(\frac\pi 2)}$ as $\theta_{\rm w}\to {\frac\pi 2}^-$
does not hold
in $C^2$ up to the sonic arc $\Sonic$ (but holds in $C^{1,\alpha}$) after mapping $\Omega^{(\theta_{\rm w})}$
to a fixed domain for all $\theta_{\rm w}$.
Moreover, the difference between the behaviors of $D\psi$ near $\Shock$ and away from $\Shock$ within $\Omega_\varepsilon$ shows that
there is no clear background solution such that the appropriate iteration set would lie in its small neighborhood
in the norm sufficiently strong to remove the
ellipticity cutoff  in \eqref{equationForPsi-sonicStruct-Iter} by the smallness of the norm.
Then, in order to remove the ellipticity cutoff for the fixed point of the iteration,
we derive an equation for $\psi_x$  in $\Omega_\varepsilon$ and boundary conditions
on $\Shock\cap\{x<\varepsilon\}$ and $\Wedge\cap\{x<\varepsilon\}$,
and prove that
$$
\psi_x\le
\frac 4{3(\gamma+1)}x
$$
from this boundary value problem,
if the wedge angle $\theta_{\rm w}$ is sufficiently close to $\frac \pi 2$.
The estimate from below:
$$
\psi_x\ge
-\frac 4{3(\gamma+1)}x
$$
is proved from the global setting of {\rm Problem \ref{fbp-c}} under the same condition on $\theta_{\rm w}$.
This use of the local and global structure is similar to that in the proof of \eqref{quadrGrowthNearSonic}.

Note that, in this argument for the wedge angles near $\frac\pi 2$, the non-perturbative nature
of the problem is seen only in the estimates of the solution near $\Sonic$,
specifically in the fact that $D^2\psi$ on $\Sonic$ does not tend to zero as $\theta_{\rm w}\to\frac{\pi}{2}$.
The free boundary
$\Shock$  in this case is near $\CS_1(\theta_{\rm w})$, and also close to the
reflected shock of the normal reflection as in Fig. \ref{NormReflFigure},
which is the  vertical line $\CS_1(\frac\pi 2)$.
Also,
$\|\varphi-\varphi_2^{(\theta_{\rm w})}\|_{C^1(\Omega)}\le C(\frac\pi 2 -\theta_{\rm w})$, which is small.
Thus, away from $\Sonic$, the argument is perturbative for the wedge angles near $\frac\pi 2$.
In the case of general wedge angles in Theorem \ref{mainShockReflThm},
the free boundary $\Shock$ is no longer close to a line, its structure is not known {\it a priori}, thus the
study of geometric properties of the free boundary is a part of the argument.

\bigskip
{\bf Case II. General wedge angles up to the detachment angle}.
For the general case and the proof of Theorem \ref{mainShockReflThm},
we follow the approach introduced in
Chen-Feldman \cite{CF-book2018}.
Similar to the case of wedge angles near $\frac\pi 2$ where we have restricted our consideration to the class of
solutions satisfying $\psi\ge 0$ in $\Omega$ and established the existence of such solutions,
for the general case,
we define a class of {\em admissible solutions}, make the necessary {\it a priori} estimates of such solutions,
and then employ these estimates to prove the existence of solutions in this class.
Our motivation for the definition of admissible solutions is from the following properties of
supersonic regular reflection solutions $\varphi$ for the wedge angles close to $\frac\pi 2$;
or more generally, for the supersonic regular reflection solutions $\varphi$ satisfying that
$\|\varphi-\varphi_2^{(\theta_{\rm w})}\|_{C^1(\Omega)}$ is small: If (\ref{2-1}) is strictly elliptic
for $\varphi$ in $\overline\Omega\setminus \overline\Sonic$, then $\varphi$ satisfies \eqref{phi-between-in-omega}
and the monotonicity properties:
\begin{equation}\label{coneOfMonotRegRefl}
\partial_{\xi_2}(\varphi_1-\varphi)\le 0, \quad D(\varphi_1-\varphi)\cdot {\bf e}_{\CS_1}\le 0
\qquad\,\,\,\, \mbox{in $\Omega\,\,$ for ${\bf e}_{\CS_1}=\frac{{P_0P_1}}{|P_0P_1|}$}.
\end{equation}

We now present the outline of the proof of Theorem \ref{mainShockReflThm} in the following four steps:

\medskip
{1}. Motivated by the discussion above,
for the general case, we define the  admissible solutions as the solutions
of {Problem \ref{fbp-c}}
(thus the solutions with weak regular reflection-diffraction configuration
of either supersonic or subsonic type)
satisfying the following properties:
\begin{definition}\label{admisSolnDef}
Let $\theta_{\rm w}\in (\theta_{\rm w}^{\rm d},\frac{\pi}2)$.
A function $\varphi\in C^{0,1}(\overline\Lambda)$ is an admissible solution of the regular reflection
problem if $\varphi$ is a solution of {\rm Problem \ref{fbp-c}} extended to $\Lambda$
by \eqref{phi-states-0-1-2-MainThm} $($where
$P_0P_1P_4$
is a point in the subsonic and
sonic cases$)$ and satisfies the following properties{\rm :}
\begin{enumerate}[\rm (i)]
\item\label{RegReflSol-Prop0}
The structure of solutions{\rm :}
		
\smallskip
\begin{itemize}
\item If $|D\varphi_2(P_0)|>c_2$, then $\varphi$ is of the {\em supersonic} regular shock
reflection-diffraction configuration shown on Fig. {\rm \ref{figure: free boundary problems-1}}	
and satisfies that the curved part of reflected-diffracted shock $\Gsh$ is  $C^{2}$ in its relative interior{\rm ;}
curves $\Gsh$,  $\Sonic$, $\Wedge$, and $\Symm$ do not have common points except their endpoints{\rm ;}
$\varphi\in C^{0,1}(\Lambda)\cap C^1(\Lambda\setminus ({\CS_0}\cup \overline{P_0P_1P_2}))$ and
$\varphi\in C^{1}(\overline{\Omega})\cap C^{3}(\overline\Omega\setminus(\overline\Gso\cup\{P_2, P_3\}))$.

\item If $|D\varphi_2(P_0)|\le c_2$, then $\varphi$ is of the {\em subsonic} regular shock
reflection-diffraction configuration  shown on Fig. {\rm \ref{figure: free boundary problems-2}}			
and satisfies that the reflected-diffracted shock $\Gsh$ is  $C^{2}$ in its relative interior{\rm ;}
curves $\Gsh$,  $\Wedge$, and $\Symm$ do not have common points except their endpoints{\rm ;}
$\varphi\in C^{0,1}(\Lambda)\cap C^1(\Lambda\setminus ({\CS_0}\cup\overline\Gsh))$
and $\varphi\in  C^{1}(\overline{\Omega})\cap
				C^3(\overline\Omega\setminus\{P_0, P_3\})$.
\end{itemize}
Moreover, in both the supersonic and subsonic cases,
the extended curve
$\Gsh^{\rm ext}:=\Gsh\cup \{P_0\}\cup\Gsh^-$ is $C^1$ in its relative interior,
where $\Gsh^-$ is the reflection of $\Gsh$ with respect to the $\xi_1$--axis.
		
\smallskip		
\item\label{RegReflSol-Prop1}
Equation \eqref{2-1} is strictly elliptic in
		$\overline\Omega\setminus\,\overline{\Gso}$, {\it i.e.},
		$|D\varphi|<c(|D\varphi|^2, \varphi)$
		 in $\overline\Omega\setminus\,\overline{\Gso}$.
		
\smallskip
\item\label{RegReflSol-Prop1-1}
$\partial_{\bn}\varphi_1>\partial_{\bn}\varphi>0$ on $\Gsh$,
where $\bn$ is the normal to $\Gsh$, pointing to the interior of $\Omega$.
		
\smallskip
\item \label{RegReflSol-Prop1-1-1}  Inequalities hold{\rm :}
\begin{equation}\label{phi-between-in-omega-nonSt}
\varphi_2\le\varphi\le \varphi_1 \qquad\mbox{in $\Omega$}.
\end{equation}
		
\smallskip
\item\label{RegReflSol-Prop2}
\eqref{coneOfMonotRegRefl} is satisfied, where vector ${\bf e}_{\CS_1}$
is defined as the unit vector parallel to $\CS_1$ and pointing into $\Lambda$ at $P_0$
for the general case.
\end{enumerate}
\end{definition}

 Note that \eqref{coneOfMonotRegRefl} implies that
\begin{equation}\label{coneOfMonotRegRefl-cone}
  D(\varphi_1-\varphi)\cdot {\bf e}\le 0 \qquad\,\, \mbox{in $\overline\Omega\,\,$ for
  all ${\bf e}\in \overline{Cone({\bf e}_{\xi_2}, {\bf e}_{\CS_1})}$},
\end{equation}
where
$Cone({\bf e}_{\xi_2}, {\bf e}_{\CS_1})=\{a\,{\bf e}_{\xi_2}+b \,{\bf e}_{\CS_1}\;:\; a, b>0\}$
with ${\bf e}_{\xi_2}=(0,1)$.
Notice that ${\bf e}_{\xi_2}$ and ${\bf e}_{\CS_1}$ are not parallel if $\theta_{\rm w}\ne \frac\pi 2$.

\medskip
{2}.
To prove the existence of admissible solutions for each wedge angle in Theorem \ref{mainShockReflThm},
we derive uniform {\it a priori} estimates for admissible solutions with
any wedge angle $\theta_{\rm w} \in [\theta_{\rm w}^{\rm d}+\sigma, \frac\pi 2]$
for each small $\sigma>0$,
show the compactness of this subset of admissible solutions
in the appropriate norm,
and then apply the degree theory to establish
the existence of admissible solutions for each $\theta_{\rm w} \in [\theta_{\rm w}^{\rm d}+\sigma, \frac\pi 2]$,
starting from the unique normal reflection solution for $\theta_{\rm w}=\frac\pi 2$.
To derive the {\it a priori} estimates for admissible solutions,
we first obtain the required estimates related to the geometry of shock $\Gamma_{\rm shock}$ and domain $\Omega$, as well as
the basic estimates of solution $\varphi$.
We prove:
\begin{enumerate}[{\rm (a)}]
\item \label{StepsOfProof-admisSol-1}
The inequality in \eqref{coneOfMonotRegRefl-cone} is strict
for any ${\bf e}\in Cone({\bf e}_{\xi_2}, {\bf e}_{\CS_1})$.
Combined with the first inequality in \eqref{phi-between-in-omega-nonSt} and
the fact that $\varphi=\varphi_1$ on $\Shock$, this implies that $\Shock$
is a Lipschitz graph with a uniform Lipschitz estimate for all admissible solutions.
\item \label{StepsOfProof-admisSol-2}
The uniform bounds on $\text{diam}(\Omega), \|\varphi\|_{C^{0,1}(\Omega)}$,
and the directional monotonicity of $\varphi-\varphi_2$ near the sonic arc for a cone of directions.
\item \label{StepsOfProof-admisSol-3}
The uniform positive lower bound for the distance between the shock and the wedge,
and the uniform separation of the shock and the symmetry line (that is,
$\shock$ is away from a uniform conical neighborhood of $\Symm$ with vertex at their common endpoint $P_2$);

\item  \label{StepsOfProof-admisSol-4}
The uniform positive lower bound for the distance between the shock and the sonic circle ${B_{c_1}((u_1, 0))}$ of
state (1), by using the properties described in
Remark \ref{positiveSolRmk}.
This allows us to estimate the ellipticity of (\ref{2-1}) for $\varphi$ in $\Omega$
(depending on the distance to the sonic arc $P_1P_4$ for the supersonic regular shock reflection-diffraction configuration
and to $\PtIncW$ for the subsonic regular shock reflection-diffraction configuration).
\item \label{StepsOfProof-admisSol-5}
Estimate \eqref{gradEstCx} holds in the supersonic case, by using
the monotonicity of $\psi=\varphi-\varphi_2$ near the sonic arc in a cone of directions
shown in \eqref{StepsOfProof-admisSol-2} and the conditions on $\Sonic$ in {\rm Problem \ref{fbp-c}}.
\end{enumerate}

The results of \eqref{StepsOfProof-admisSol-1}--\eqref{StepsOfProof-admisSol-3} are obtained via the maximum principle,
by considering equation \eqref{equ:study} as a linear elliptic equation for $\phi$ and using the boundary conditions
on $\Shock$, $\Sonic$, $\Wedge$, and $\Symm$ in  {Problem \ref{fbp-c}} and
\eqref{phi-between-in-omega-nonSt}--\eqref{coneOfMonotRegRefl-cone}.
The results of \eqref{StepsOfProof-admisSol-3}, combined with \eqref{StepsOfProof-admisSol-1},
show the structure of $\Omega$ which allows us to perform the uniform local elliptic estimates in various parts of $\Omega$:
the interior,  near a point $P$ in a relative interior of $\Shock$, $\Wedge$, and $\Symm$,
and locally near corners $P_2$ and $P_3$.

\smallskip
Based on estimates \eqref{StepsOfProof-admisSol-1}--\eqref{StepsOfProof-admisSol-4},  we show
the uniform regularity estimates for the solution and the free boundary in the weighted and scaled $C^{2, \alpha}$
norms away from the sonic arc in the supersonic case and away from $P_0$ in the subsonic case,
{\it i.e.}, in $\Omega\setminus\Omega_\varepsilon$, for any small $\varepsilon>0$, for some $\alpha\in(0,1)$.
The equation is uniformly elliptic in this region, with the ellipticity constant depending on $\varepsilon$.
Thus, the estimates depend on $\varepsilon$.

\medskip
{3}. Below we discuss the estimates near $\Sonic$ (resp. near $P_0$ in the subsonic/sonic case), {\it i.e.},
the estimates in $\Omega_{2\varepsilon}$ for some $\varepsilon$, independent of
$\theta_{\rm w} \in [\theta_{\rm w}^{\rm d}+\sigma, \frac\pi 2]$,
which allows us to complete the uniform {\it a priori} estimates for admissible solutions
with wedge angles $\theta_{\rm w} \in [\theta_{\rm w}^{\rm d}+\sigma, \frac\pi 2]$.
We obtain the estimates near $\Sonic$ (or $\PtIncW$
for the subsonic reflection), {\it i.e.}, in $\Omega_{2\varepsilon}$,
in scaled and weighted $C^{2,\alpha}$ for $\varphi$
and the free boundary $\Shock\cap\partial\Omega_{2\varepsilon}$, by considering separately
four cases depending on $\frac{|D\varphi_2|}{c_2}$ at $\PtIncW$:

\begin{itemize}
\item[(i)] Supersonic:
$\frac{|D\varphi_2|}{c_2} \ge 1+\delta$,

\smallskip
\item[(ii)]
 Supersonic (almost sonic): $1<
\frac{|D\varphi_2|}{c_2} < 1+\delta$,

\smallskip
\item[(iii)]
Subsonic (almost sonic,  including sonic): $1-\delta\le \frac{ |D\varphi_2|}{c_2} \le 1$,

\smallskip
\item[(iv)]
Subsonic: $\frac{|D\varphi_2|}{c_2} \le 1-\delta$,
\end{itemize}
for small $\delta>0$ chosen so that the estimates can be obtained.
The choice of $\delta$ determines $\varepsilon$.

For cases (i)--(ii), equation (\ref{2-1}) is degenerate elliptic in $\Omega$
near  $\PtUpL\PtUpR$ on Fig. \ref{figure: free boundary problems-1}.
For case (iii), except the sonic case $ \frac{ |D\varphi_2(P_0)|}{c_2} = 1$, the equation is uniformly elliptic in $\overline\Omega$,
but the ellipticity constant
is small and tends to zero near $\PtIncW$ on Fig. \ref{figure: free boundary problems-2} as
$\frac{ |D\varphi^{(\theta_{\rm w})}_2(P_0)|}{c_2}\to 1^-$,
{\it i.e.}, as the subsonic angles $\theta_{\rm w}$ tend to the sonic  angle.
Thus, for cases (i)--(iii), we exploit the local elliptic degeneracy,
which allows us to find a comparison function in each case, to show
the appropriately fast decay of $\varphi-\varphi_2$ near $\PtUpL\PtUpR$ for cases (i)--(ii)
and near $\PtIncW$ for case (iii);
furthermore, combining with appropriate local non-isotropic rescaling to obtain the uniform ellipticity,
we obtain the {\it a priori} estimates in the weighted and scaled $C^{2,\alpha}$--norms.
In cases (i)--(ii), the norms are \eqref{parabNormsApp}.
For case (iii), we use the different norms to obtain the estimates that imply the standard
$C^{2,\alpha}$--estimates.
To obtain these estimates, for case (i), we use
the argument developed in Chen-Feldman \cite{ChenFeldman}
and described above (see Remark \ref{apriri-est-positiveSolRmk}), where
the ellipticity estimate \eqref{elliptCond} follows from the estimates
described in \eqref{StepsOfProof-admisSol-4} above and \eqref{gradEstCx} obtained
in \eqref{StepsOfProof-admisSol-5}.
These estimates hold in $\Omega_\varepsilon$ with
$\varepsilon \lesssim (\mbox{length}(\Sonic))^2$ because the {\it rectangles}
$R_{(x_0, y_0)}$ defined by \eqref{parabRectangles} do not fit into $\Omega$ for larger $x_0$,
which means, for example, that $R_{(x_0, y_0)}\cap\Wedge\ne\emptyset$ for $(x_0, y_0)\in \shock\cap\partial\Omega_\varepsilon$
with $x_0\ge C(\mbox{length}(\Sonic))^2$
if $C$ is fixed and $\mbox{length}(\Sonic)$ is small, because the length of the $y$-side of
$R_{(x_0, y_0)}$ is $\frac{\sqrt{x_0}}{10}$, and $\Shock$ and $\Wedge$ are smooth curves
that intersect $\Sonic$ transversally.
However, $\mbox{length}(\Sonic)$
tends to zero, as $\displaystyle \frac{|D\varphi^{(\theta_{\rm w})}_2(P_0)|}{c_2}\to 1^+$,
{\it i.e.}, when the supersonic wedge angle tends to the sonic angle.
Thus, a different argument, involving an appropriate scaling, is employed for case (ii) in order
to keep $\varepsilon$ uniform for all $\theta_{\rm w} \in [\theta_{\rm w}^{\rm d}+\delta, \frac\pi 2]$.
Another version of that argument (with a different scaling) is applied for case (iii).
For both cases (ii)--(iii), we need to use smaller rectangles than those for case (i),
but this requires stronger growth estimates than \eqref{quadrGrowthNearSonic}
to obtain a bound in $C^{1,1}$ from the corresponding weighted and scaled estimates.
We obtain such growth estimates by using the conditions of cases (ii)--(iii)
for sufficiently small $\delta$.
 For  case (iv), the equation is uniformly elliptic
in $\Omega$ for the admissible solution, where the ellipticity constant is not small,
and the estimates are more technically challenging than those for cases (i)--(iii).
This can be seen as follows:
For all cases (i)--(iv),
the free boundary has a lower {\it a priori} regularity in
the sense that only the Lipschitz estimate of $\Shock$ is obtained in
\eqref{StepsOfProof-admisSol-1} above; however, for case (iv),
the uniform ellipticity
combined with oblique boundary conditions does not allow
a comparison function that leads to the fast decay of $|\varphi-\varphi_2|$ near $\PtIncW$.
Thus, we prove the $C^{\alpha}$--estimates of $D(\varphi-\varphi_2)$ near $P_0$,
by deriving the equations and boundary conditions for two directional derivatives of
$\varphi-\varphi_2$ near $P_0$,
and performing the hodograph transform
to flatten the free boundary.

\medskip
{4.} In order to prove the existence of solutions, we perform an iteration,
which is an extension of the iteration process used in Chen-Feldman \cite{ChenFeldman}.
First, given an admissible solution $\varphi$ for the wedge angle $\theta_{\rm w}$,
we map its elliptic domain $\Omega(\varphi, \theta_{\rm w})$ to a unit square $Q=(0, 1)^2$
so that, for the supersonic case,
the boundary parts $\Shock$, $\Sonic$, $\Wedge$, and $\Symm$ are mapped to the respective sides of $Q$,
and the other properties of this map are satisfied.
For the subsonic case,
the map is discontinuous at $P_0=\overline{\Sonic}$ (mapping the triangular domain to a square).
Moreover, we define a function $u$ on $Q$ by expressing $u:=\varphi-\tilde\varphi_2^{(\theta_{\rm w})}$
in the coordinates on $Q$, where $\tilde\varphi_2^{(\theta_{\rm w})}$ is a function determined
by $\theta_{\rm w}$ and equals to $\varphi_2$ near $\overline{\Sonic}$;
we skip the complete technical definition here.
For appropriate functions $u$ on $Q$ and the wedge angle $\theta_{\rm w}$, this map can be inverted,
{\it i.e.}, the elliptic domain $\Omega(u, \theta_{\rm w})$ and the iteration free boundary
$\Shock(u, \theta_{\rm w})$ can be determined, and a function $\varphi^{(u, \theta_{\rm w})}$ on
$\Omega(u, \theta_{\rm w})$ is defined by expressing $u$ in the coordinates on
$\Omega(u, \theta_{\rm w})$ and adding $\tilde\varphi_2^{(\theta_{\rm w})}$ so that,
if $u$ is obtained from the admissible solution $\varphi$ with the elliptic domain $\Omega$ as described above,
then $\Omega(u, \theta_{\rm w})=\Omega$ and $\varphi^{(u, \theta_{\rm w})}=\varphi$ in $\Omega$.
Moreover, the map: $\Omega(u, \theta_{\rm w})\to Q$ and its inverse
satisfy certain continuity properties with respect to $(u, \theta_{\rm w})$.
The iteration is performed in terms of the functions defined on $Q$.
The iteration set consists of pairs $(u, \theta_{\rm w})$, where $u$ is in a weighted and scaled
$C^{2,\alpha}$ space on $Q$, denoted as $C^{2,\alpha}_{**}$ (its definition is technical, so we skip it here),
and satisfy
\begin{enumerate}[{\rm (i)}]
\item $\|u\|_{C^{2,\alpha}_{**}}\le M(\theta_{\rm w})$, where $M(\theta_{\rm w})$ is defined explicitly,
   based on the {\it a priori} estimates discussed above;
\item $\Omega(u, \theta_{\rm w})$, $\Shock(u, \theta_{\rm w})$, and
$\varphi^{(u, \theta_{\rm w})}$
on $\Omega(u, \theta_{\rm w})$ satisfy some geometric and analytical properties.
\end{enumerate}
The iteration map: $(\hat u, \theta_{\rm w})\to (u, \theta_{\rm w})$ is defined by solving the iteration
problem on $\Omega(u, \theta_{\rm w})$ and then mapping its solution $\varphi$
to a function $u$ on $Q$. This mapping includes additional steps, compared to the one described above.
Specifically, we modify the iteration free boundary by using the solution $\varphi$ of the iteration problem so that,
in the mapping: $(\varphi, \theta_{\rm w})\to u$,
the resulting function $u$ on $Q$ keeps
the regularity obtained from solving the iteration problem.
This yields the compactness of the iteration map.
We prove that, for a fixed point $(u, \theta_{\rm w})$ of the iteration map,
$\varphi^{(u, \theta_{\rm w})}$ on $\Omega(u, \theta_{\rm w})$ is an admissible solution.
We use the degree theory to establish
the existence of admissible solutions as fixed points of the iteration map
for each $\theta_{\rm w} \in [\theta_{\rm w}^{\rm d}+\delta, \frac\pi 2]$,
starting from the unique normal reflection solution for $\theta_{\rm w}=\frac\pi 2$.
The compactness of the iteration map described above is necessary for that.
The {\it a priori} estimates of admissible solutions discussed above are used in the degree theory argument
in order to define the iteration set such that a fixed point of the iteration map ({\it i.e.}, admissible solution)
cannot occur on the boundary of the iteration set, since that would contradict the {\it a priori} estimates.
With all of these arguments, we complete the proof of Theorem \ref{mainShockReflThm}.
This provides a solution to the von Neumann's conjectures.

\smallskip
More details can be found in Chen-Feldman \cite{CF-book2018};
also see \cite{ChenFeldman}.

\subsection{The Prandtl-Meyer Problem for Unsteady
Supersonic Flow onto Solid Wedges}
\label{PrandtlMeyerProblSect}

As we discussed in \S 2--\S3, steady shocks appear when a steady supersonic flow hits a straight wedge;
see Figure \ref{Figure2}.
Since both weak and strong steady shock solutions are stable in the steady regime,
the static stability analysis alone
is not able to single out one of them in this sense,
unless an additional condition is posed on the speed of the downstream flow at infinity.
Then the dynamic stability analysis becomes more significant
to understand the non-uniqueness issue of the steady oblique shock solutions.
However, the problem for the dynamic stability of the steady shock solutions for
supersonic flow past solid wedges involves several additional
difficulties.
The recent efforts have been focused on the construction of the global Prandtl-Meyer reflection
configurations in the self-similar coordinates for potential flow.

As we discussed earlier, if a supersonic flow
with a constant density $\irho >0$ and a velocity ${\bf u}_{0}=(\iu, 0)$,
$\iu> c_0:=c(\rho_0)$,
impinges toward wedge $W$ in \eqref{wedge-1},
and if $\theta_{\rm w}$ is less
than the detachment angle $\theta_{\rm w}^{\rm d}$,
then the well-known {\emph{shock polar analysis}} shows that there are two different
steady weak solutions:
{\emph{the steady weak shock solution}} $\bar{\Phi}$ and {\emph{the steady strong shock solution}},
both of which satisfy the entropy condition and the slip boundary condition (see Fig.  \ref{Figure2}).

Then the dynamic stability of the weak transonic shock solution for potential flow can be
formulated as the following problem:

\begin{problem}[Initial-Boundary Value Problem]
\label{problem-1}
Given $\gam>1$, fix $(\irho, \iu)$ with $\iu>c_0$.
For a fixed $\theta_{\rm w}\in (0,\theta_{\rm w}^{\rm d})$,
let $W$ be given by \eqref{wedge-1}.
Seek a global entropy solution $\Phi\in W^{1,\infty}_{\rm loc}(\R_+\times (\R^2\setminus W))$
of Eq. \eqref{1-b2} with $\rho$ determined by \eqref{1-b1}
and $B=\frac{\iu^2}{2}+h(\irho)$ so that $\Phi$ satisfies the initial condition at $t=0${\rm :}
\begin{equation}\label{1-d}
(\rho,\Phi)|_{t=0}=(\irho, \iu x_1) \qquad\,\,\, \text{for}\;\;\xx\in \R^2\setminus W,
\end{equation}
and the slip boundary condition along the wedge boundary $\der W${\rm :}
\begin{equation}\label{1-e}
\nabla_{\bf x}\Phi\cdot \nnu_{\rm w} |_{\der W}=0,
\end{equation}
where $\nnu_{\rm w}$ is the exterior unit normal to $\der W$.

In particular, we seek a solution $\Phi\in W^{1,\infty}_{\rm loc}(\R_+\times (\R^2\setminus W))$
that converges to
the steady weak oblique shock solution $\bPhi$ corresponding to
the fixed parameters $(\irho, \iu, \gam, \theta_{\rm w})$
with $\bar{\rho}=h^{-1}(B-\frac 12|\nabla \bPhi|^2)$,
when $t\to \infty$, in the following sense{\rm :}
For any $R>0$, $\Phi$ satisfies
\begin{equation}
\label{time-asymp-lmt}
\lim_{t\to \infty} \|(\nabla_{\bf x}\Phi(t,\cdot)-\nabla_{\bf x}\bPhi,
\rho(t,\cdot)-\bar{\rho})\|_{L^1(B_R({\bf 0})\setminus W)}=0
\end{equation}
for $\rho(t,{\bf x})$ given by \eqref{1-b1}.
\end{problem}

\smallskip
Since the initial data functions in \eqref{1-d} do not satisfy the boundary condition \eqref{1-e},
a boundary layer is generated along the wedge boundary starting at $t=0$,
which forms the Prandtl-Meyer reflection configurations; see
Bae-Chen-Feldman \cite{BCF-14} and the references cited therein.

\smallskip
Notice that the initial-boundary value problem, {Problem \ref{problem-1}},
is invariant under scaling \eqref{4.5}. Thus,
we study the existence of self-similar solutions
determined by
equation \eqref{2-1} with \eqref{1-o} through \eqref{4.6}.

As the upstream flow has the constant velocity $(\iu,0)$,
noting the choice of $B$ in Problem \ref{problem-1},
the corresponding pseudo-potential
$\ivphi$ has the expression of
\begin{equation}
\label{1-m}
\ivphi=-\frac 12|\xxi|^2+\iu\xi_1
\end{equation}
in self-similar coordinates $\xxi=\frac{{\bf x}}{t}$,
as shown directly from \eqref{constantStatesForm}.
Notice also the symmetry of the domain and the upstream flow in {Problem \ref{problem-1}}
with respect to the $x_1$--axis.
{Problem \ref{problem-1}} can then be reformulated as the following boundary value problem
in the domain:
$$
\Lambda:=\R^2_+\setminus\{\xxi\,: \,\xi_2\le \xi_1\tan\theta_{\rm w},\, \xi_1\ge 0\}
$$
in the
self-similar coordinates $\xxi$,
which corresponds to domain $\{(t, {\bf x})\, :\, {\bf x}\in \R^2_+\setminus W,\, t>0\}$
in the $(t, {\bf x})$--coordinates, where $\R^2_+=\{\xxi\,: \,\xi_2>0\}$.

\medskip
\begin{problem}[Boundary Value Problem]
\label{problem-5.2}
Seek a solution $\vphi$ of equation \eqref{2-1} in the self-similar domain $\Lambda$
with the slip boundary condition{\rm :}
\begin{equation}
\label{1-k}
D\vphi\cdot \nnu|_{\partial \Lambda}=0
\end{equation}
and the asymptotic boundary condition{\rm :}
\begin{equation}\label{1-k-b}
\vphi-\vphi_0\longrightarrow 0
\end{equation}
along each ray $R_\theta:=\{ \xi_1=\xi_2 \cot \theta, \xi_2 > 0 \}$
with $\theta\in (\theta_{\rm w}, \pi)$ as $\xi_2\to \infty$
in the sense that
\begin{equation}\label{1-k-c}
\lim_{r\to \infty} \|\varphi  - \varphi_0\|_{C(R_{\theta}\setminus B_r(0))} = 0.
\end{equation}
\end{problem}

\begin{figure}
 	\centering
 	\includegraphics[height=35mm]{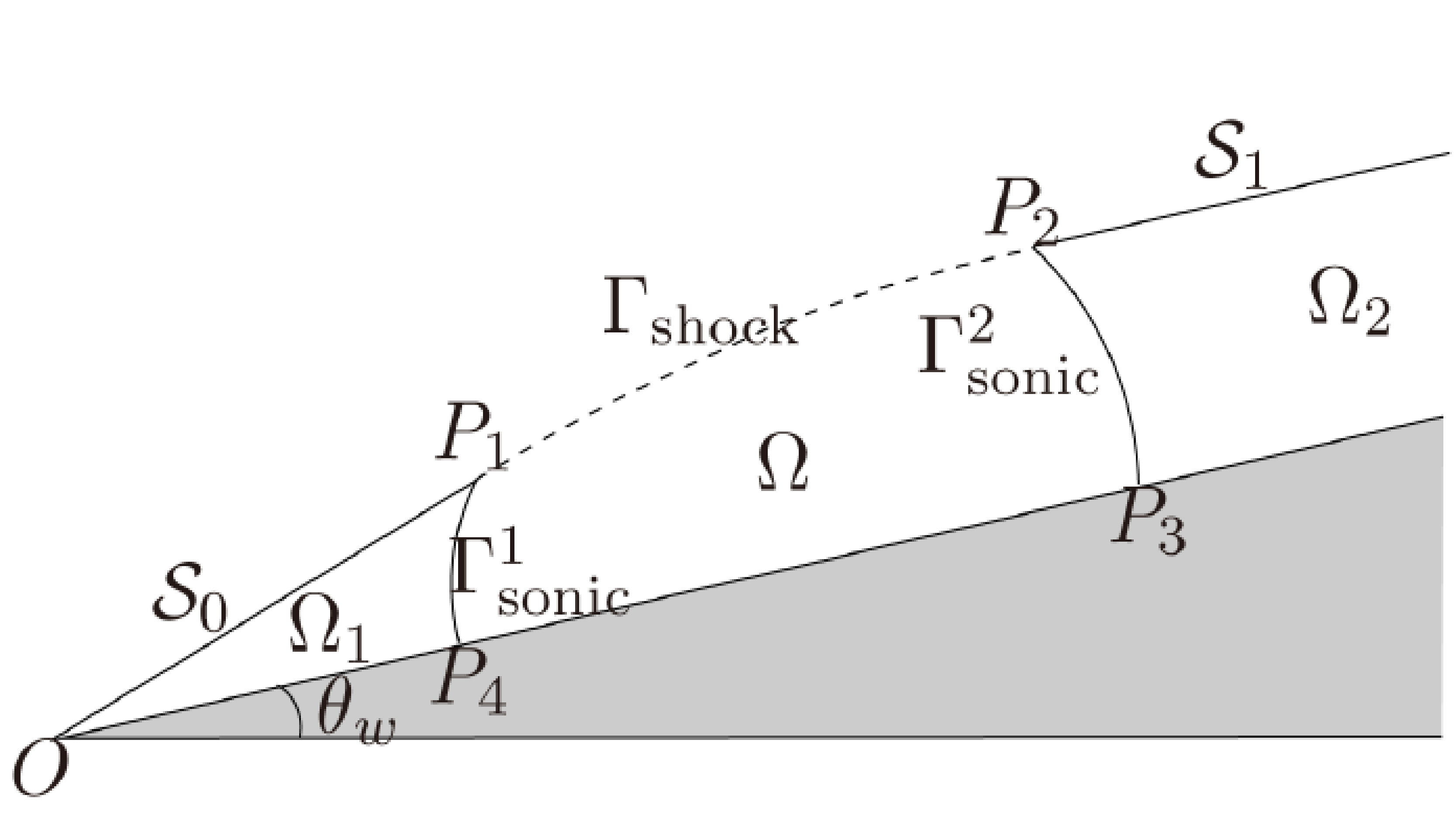}
		\caption{Self-similar solutions for $\theta_{\rm w}\in (0, \theta_{\rm w}^{\rm s})$
in the self-similar coordinates $\xxi$ ({\it cf}.  \cite{BCF-14})}\label{fig:global-structure-1}
 \end{figure}

\begin{figure}
\centering
\includegraphics[height=35mm]{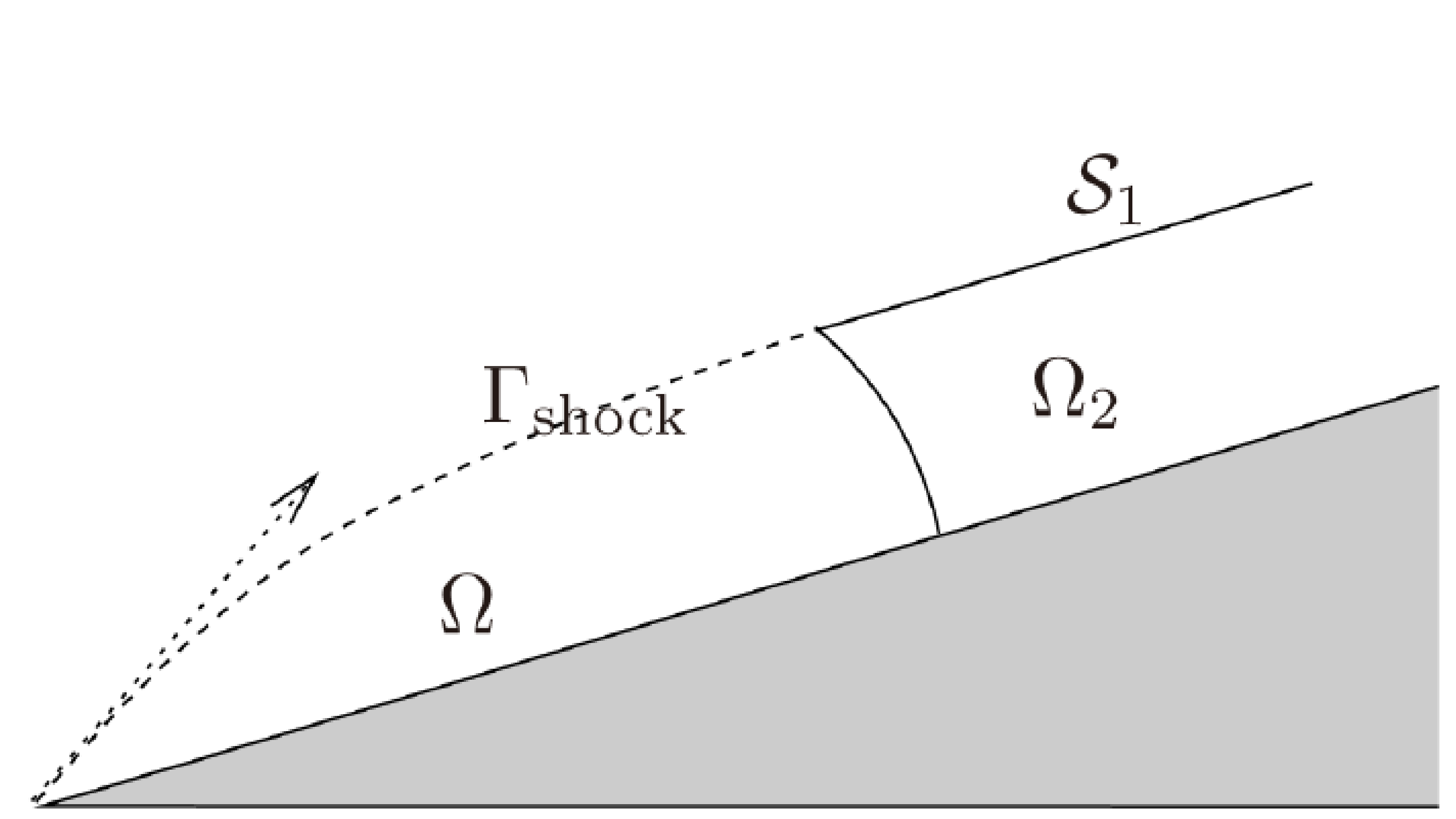}
\caption{Self-similar solutions for $\theta_{\rm w}\in [\theta_{\rm w}^{\rm s},\theta_{\rm w}^{\rm d})$
in the self-similar coordinates $\xxi$ ({\it cf}. \cite{BCF-14})}\label{fig:global-structure-2}
\end{figure}

In particular, we seek a global entropy solution of {Problem \ref{problem-5.2}}
with two types of Prandtl-Meyer reflection configurations whose occurrence
is determined by the wedge angle $\theta_{\rm w}$ for the two different cases:
One contains a straight weak oblique shock ${\mathcal S}_0$ attached to the wedge vertex $O$ and
connected to a normal shock ${\mathcal S}_1$ through a curved shock  $\Gamma_{\rm shock}$
when $\theta_{\rm w}<\theta_{\rm w}^{\rm s}$,
as shown in Fig. \ref{fig:global-structure-1};
the other contains a curved shock  $\Gamma_{\rm shock}$ attached
to the wedge vertex
and connected to a normal shock ${\mathcal S}_1$
when $\theta_{\rm w}^{\rm s}\le \theta_{\rm w}<\theta_{\rm w}^{\rm d}$,
as shown in Fig. \ref{fig:global-structure-2},
in which the curved shock $\Gamma_{\rm shock}$ is tangential to
the straight weak oblique shock  ${\mathcal S_0}$ at the wedge vertex.

To seek a global entropy solution of {Problem \ref{problem-5.2}} with the structure
of Fig. \ref{fig:global-structure-1} or Fig. \ref{fig:global-structure-2},
one needs to compute the pseudo-potential function $\vphi$ below ${\mathcal S_0}$.

Given $M_0>1$, $\rho_1$ and $\uu_1$ are determined by using the shock polar as in Fig. \ref{Figure2}
for steady potential flow
(note that the shock polar is now different from the one for
the full Euler system but has the same shape as in Fig. \ref{Figure2}).
Similar to those in \S 3.1, in the potential flow case,
for any wedge angle $\theta_{\rm w}\in (0,\theta_{\rm w}^{\rm s})$,
line $u_2=u_1\tan\theta_{\rm w}$ and the shock polar
intersect at a point $\uu_1$ with $|\uu_1| >c_1$ and $u_{11}<\iu$;
while, for any wedge angle $\theta_{\rm w}\in [\theta_{\rm w}^{\rm s}, \theta_{\rm w}^{\rm d})$,
they intersect at a point $\uu_1$ with $u_{11}>u_{\rm 1d}$
and $|\uu_1|<c_1$, where $u_{\rm 1d}$ is the $u_1$--component
of the unique detachment state $\uu_{\rm d}$ when $\theta_{\rm w}=\theta_{\rm w}^{\rm d}$.
The intersection state $\uu_1$ is the velocity for steady potential flow
behind an oblique shock $\mathcal{S}_0$ attached to the wedge vertex with angle $\theta_{\rm w}$.
The strength of shock $\mathcal{S}_0$ is relatively weak compared to the other shock
given by the other intersection point on the shock polar,
hence we call $\mathcal{S}_0$
a \emph{weak oblique shock},
and the corresponding state $\uu_1$ is a \emph{weak state}.

We also note that states
$\uu_1$ depend smoothly on $\iu$ and $\theta_{\rm w}$, and
such states are supersonic when $\theta_{\rm w}\in (0,\theta_{\rm w}^{\rm s})$
and subsonic when $\theta_{\rm w}\in [\theta_{\rm w}^{\rm s}, \theta_{\rm w}^{\rm d})$.

Once $\uu_1$ is determined, by  \eqref{1-i} and \eqref{1-m},
the pseudo-potential $\vphi_1$ below the weak oblique shock $\mathcal{S}_0$
is
\begin{equation}\label{1-n}
\vphi_1=-\frac 12|\xxi|^2+ \uu_1\cdot\xxi.
\end{equation}
Similarly, by  \eqref{1-h}--\eqref{1-i} and \eqref{1-m}--\eqref{1-k},
the pseudo-potential $\vphi_2$ below the normal shock $\mathcal{S}_1$
is of the form:
\begin{equation}\label{1-n-b}
\vphi_2=-\frac 12|\xxi|^2+ \uu_2\cdot \xxi+ k_2
\end{equation}
for constant state $\uu_2$ and constant $k_2$; see \eqref{constantStatesForm}.
Then it follows
from \eqref{1-o} and \eqref{1-n}--\eqref{1-n-b} that the corresponding densities $\rho_1$ and $\rho_2$
are constants, respectively.
In particular, we have
\begin{equation}\label{2-n1}
\rho_k^{\gam-1}=\rho_0^{\gamma-1}+\frac{\gam-1}{2}\big(\iu^2-|\uu_k|^2\big)\qquad\,\,\mbox{for $k=1,2$}.
\end{equation}

Denote $\Wedge :=\partial W\cap\partial\Lambda$.
Next, we define the sonic arcs $\Gamma_{\rm sonic}^1=P_1P_4$
on Fig. {\rm \ref{fig:global-structure-1}} and $\Gamma_{\rm sonic}^2=P_2P_3$ on
Figs. {\rm \ref{fig:global-structure-1}}--{\rm \ref{fig:global-structure-2}}.
The sonic circle $\partial B_{c_1}(\uu_1)$ of the uniform state $\varphi_1$ intersects line ${\mathcal S_0}$,
where $c_1=\rho_1^{\frac{\gamma-1}2}$ by \eqref{c-through-density-function}.
For the supersonic case $\theta_{\rm w}\in (0,\theta_{\rm w}^{\rm s})$,
there are two arcs of this sonic circle between ${\mathcal S_0}$ and $\Wedge$ in $\Lambda$.
We denote  by $\Gamma_{\rm sonic}^1$ the {\em lower} arc ({\it i.e.}, located to the left from another arc)
in the orientation on  Fig. {\rm \ref{fig:global-structure-1}}.
Note that $\Gamma_{\rm sonic}^1$ tends to point $O$ as $\theta_{\rm w} \nearrow \theta_{\rm w}^{\rm s}$
and is outside of $\Lambda$ for the subsonic case
$\theta_{\rm w}\in [\theta_{\rm w}^{\rm s}, \theta_{\rm w}^{\rm d})$.
Similarly, the sonic circle $\partial B_{c_2}(\uu_2)$ of the uniform state $\varphi_2$ intersects line ${\mathcal S_1}$,
where $c_2=\rho_2^{\frac{\gamma-1}2}$. There are two
arcs of this circle between ${\mathcal S_1}$ and the line containing $\Wedge$.
For all $\theta_{\rm w}\in (0, \theta_{\rm w}^{\rm d})$,
the {\em upper} arc ({\it i.e.}, located to the right of the other arc)
in the orientation on Figs. {\rm \ref{fig:global-structure-1}}--{\rm \ref{fig:global-structure-2}}
is within $\Lambda$, which is denoted as $\Gamma_{\rm sonic}^2$.

Then {Problem \ref{problem-5.2}} can be reformulated into the following
free boundary problem:

\smallskip
\begin{problem}[Free Boundary Problem]\label{fbp-a}
For $\theta_{\rm w}\in (0, \theta_{\rm w}^{\rm d})$,
find a free boundary {\rm (}curved shock{\rm )} $\shock$ and a function $\vphi$ defined in domain $\Omega$,
as shown in Figs. {\rm \ref{fig:global-structure-1}}--{\rm \ref{fig:global-structure-2}},
such that $\vphi$ satisfies
\begin{itemize}
\item[\rm (i)]
Equation \eqref{2-1} in $\Omega$,
\item[\rm (ii)]
$\vphi=\ivphi$ and $\rho D\vphi\cdot\nnu_{\rm s}=\rho_0 D\ivphi\cdot\nnu_{\rm s}$ {on} $\shock$,
\item[\rm (iii)]
$\vphi=\hat{\vphi}$ and $D\vphi=D\hat{\vphi}$ {on} $\Gamma_{\rm sonic}^1\cup\Gamma_{\rm sonic}^2\,\,$
when $\theta_{\rm w}\in (0, \theta_{\rm w}^{\rm s})$
and on $\Gamma_{\rm sonic}^2\cup \{O\}$ when $\theta_{\rm w}\in [\theta_{\rm w}^{\rm s}, \theta_{\rm w}^{\rm d})$
for $\hat{\vphi}:=\max(\vphi_1, \vphi_2)$,
\item[\rm (iv)]
$D\vphi\cdot \nnu_{\rm w}=0$ {on} $\Wedge$,
\end{itemize}
where $\nnu_{\rm s}$ and $\nnu_{\rm w}$ are
the unit normals to
$\shock$ and $\Wedge$ pointing to the interior
of $\Omega$, respectively.
\end{problem}

\begin{remark}
It can be shown that $\varphi_1>\varphi_2$ on $\Gamma_{\rm sonic}^1$ and the opposite inequality holds
on $\Gamma_{\rm sonic}^2$. This justifies the requirements in {\rm Problem \ref{fbp-a}(iii)} above.
\end{remark}

\begin{remark}
Similar to {\rm Problem  \ref{fbp-c}}, the conditions in
{\rm Problem \ref{fbp-a}}{\rm (ii)}--{\rm (iii)}
are the Rankine-Hugoniot conditions
\eqref{1-h}--\eqref{1-i} on $\Shock$ and
$\Gamma_{\rm sonic}^1\cup\Gamma_{\rm sonic}^2$ or $\Gamma_{\rm sonic}^2\cup \{O\}$,
respectively; see the discussions right after  {\rm Problem \ref{fbp-c}}.
\end{remark}

Let $\vphi$ be a solution of {Problem \ref{fbp-a}}
such that $\shock$ is a $C^1$--curve up to its endpoints and $\vphi\in C^1(\overline\Omega)$.
To obtain a solution of {Problem \ref{problem-5.2}} from $\vphi$,
we consider two cases:

\newcommand{\SzSeg}{{\mathcal{S}_{0,{\rm seg}}}}
\newcommand{\SoSeg}{{\mathcal{S}_{1,{\rm seg}}}}
For the supersonic case $\theta_{\rm w}\in (0, \theta_{\rm w}^{\rm s})$,
we divide region $\Lambda$ into four separate regions; see Fig. \ref{fig:global-structure-1}.
We denote by $\SzSeg$ the line segment $OP_1\subset \CS_0$, and
by $\SoSeg$ the portion (half-line) of $\CS_1$ with left endpoint $P_2$
so that $\SoSeg\subset\Lambda$.
Let $\Omega_{\mathcal{S}}$ be the unbounded domain below
curve $\ol{\SzSeg\cup\shock\cup \SoSeg}$
and above $\Wedge$
(see Fig. \ref{fig:global-structure-1}).
In $\Omega_{\mathcal{S}}$, let $\Omega_1$ be the bounded domain enclosed
by $\mathcal{S}_0, \Gamma^1_{\rm sonic}$,
and $\Wedge$.
Set $\Omega_{2}:=\Omega_\CS\setminus \ol{\Omega_1\cup\Omega}$.
Define a function $\vphi_*$ in $\Lambda$ by
\begin{equation}\label{extsol}
\vphi_*=
\begin{cases}
\ivphi& \qquad \text{in}\,  \Lambda\setminus \Omega_{\mathcal{S}},\\
\vphi_1& \qquad \text{in}\,\Omega_1,\\
\vphi& \qquad \text{in}\, \Gamma^1_{\rm sonic}\cup\Omega\cup\Gamma^2_{\rm sonic},\\
\vphi_2&\qquad \text{in}\,\Omega_2.
\end{cases}
\end{equation}
By {Problem \ref{fbp-a}}(ii)--(iii),
$\vphi_*$ is continuous in $\Lambda\setminus\Omega_\CS$
and $C^1$ in $\overline{\Omega_{\mathcal{S}}}$.
In particular, $\vphi_*$ is $C^1$ across $\Gamma^1_{\rm sonic}\cup\Gamma^2_{\rm sonic}$.
Moreover, using  {Problem \ref{fbp-a}}(i)--(iii), we obtain
that  $\vphi_*$ is a global entropy solution of equation \eqref{2-1}
in $\Lambda$.

For the subsonic case $\theta_{\rm w}\in [\theta_{\rm w}^{\rm s}, \theta_{\rm w}^{\rm d})$,
region $\Omega_1\cup \Gamma^1_{\rm sonic}$ in $\varphi_*$ reduces to one point $\{O\}$;
see Fig. \ref{fig:global-structure-2}.
The corresponding function $\varphi_*$ is a global entropy solution of equation \eqref{2-1}
in $\Lambda$.

\smallskip
The first unsteady analysis of the steady supersonic weak shock solution
as the long-time behavior of an unsteady flow
is due to Elling-Liu \cite{EllingLiu1},
in which they
succeeded in establishing a stability theorem for an important class of physical
parameters determined by certain assumptions for the wedge angle $\theta_{\rm w}$
less than the sonic angle $\theta^{\rm s}_{\rm w}\in (0, \theta^{\rm d}_{\rm w})$
for potential flow.

Recently, in Bae-Chen-Feldman \cite{BCF-14},
we have removed the assumptions
\cite{EllingLiu1}
and established the stability theorem for the steady (supersonic or transonic)
weak shock solutions as
the long-time asymptotics of
the global Prandtl-Meyer reflection configurations for unsteady potential flow
for all the admissible physical parameters
even up to the detachment angle $\theta^{\rm d}_{\rm w}$
(beyond the sonic angle $\theta^{\rm s}_{\rm w}<\theta^{\rm d}_{\rm w}$).

\medskip
To achieve this, we solve
the free boundary problem ({Problem  \ref{fbp-a}}),
involving transonic shocks,
for all the wedge angles
$\theta_{\rm w}\in (0, \theta_{\rm w}^{\rm d})$
by employing the
techniques developed in Chen-Feldman
\cite{CF-book2018}, described in \S \ref{RegReflProblSect} above.
Similar to Definition \ref{admisSolnDef}, we define admissible solutions in the present case:

\begin{definition}\label{admisSolnDef-Prandtl}
Let $\theta_{\rm w}\in (0, \theta_{\rm w}^{\rm d})$.
A function $\varphi\in C^{0,1}(\overline\Lambda)$ is an admissible solution of
{\rm Problem \ref{fbp-a}}
if $\varphi$ is a solution of {\rm Problem \ref{fbp-a}} extended to $\Lambda$
by \eqref{extsol} and satisfies the following properties{\rm :}
\begin{enumerate}[\rm (i)]
\item\label{RegReflSol-PropPrandtl0}
The structure of solutions is as follows{\rm :}
		
\smallskip
\begin{itemize}
\item
If $\theta_{\rm w}\in (0, \theta_{\rm w}^{\rm s})$,	then $\varphi$ has the configuration
shown on Fig. {\rm \ref{fig:global-structure-1}}			
such that
$\Gsh$ is $C^{2}$ in its relative interior,
$\varphi\in C^{0,1}(\Lambda)\cap C^1(\Lambda\setminus (\overline{\SzSeg}\cup\overline\Gsh\cup \overline{\SoSeg}))$,
and $\varphi\in  C^{1}(\overline{\Omega})\cap
C^2(\overline\Omega\setminus(\overline{\Gamma_{\rm sonic}^1}\cup\overline{\Gamma_{\rm sonic}^2}))
\cap C^3(\Omega)$.

\smallskip
\item If $\theta_{\rm w}\in [ \theta_{\rm w}^{\rm s}, \theta_{\rm w}^{\rm d})$,	
then $\varphi$ has the configuration  shown on Fig. {\rm \ref{fig:global-structure-2}}
such that
$\Gsh$ is $C^{2}$ in its relative interior,
$\varphi\in C^{0,1}(\Lambda)\cap C^1(\Lambda\setminus (\Gsh\cup \overline{\SoSeg}))$,
and $\varphi\in  C^{1}(\overline{\Omega})\cap
				C^2(\overline\Omega\setminus(\{O\}\cup \overline{\Gamma_{\rm sonic}^2}))\cap
				C^3(\Omega)$.
\end{itemize}

\smallskip		
\item\label{RegReflSol-PropPrandtl1}
Equation \eqref{2-1} is strictly elliptic in
$\overline\Omega\setminus (\overline{\Gamma_{\rm sonic}^1}\cup\overline{\Gamma_{\rm sonic}^2})$,
{\it i.e.},
$|D\varphi|<c(|D\varphi|^2, \varphi)$
in $\overline\Omega\setminus (\overline{\Gamma_{\rm sonic}^1}\cup\overline{\Gamma_{\rm sonic}^2})$.
	
\smallskip
\item\label{RegReflSol-PropPrandtl1-1}
$\partial_{\bn}\varphi_0>\partial_{\bn}\varphi>0$ on $\Gsh$, where $\bn$ is the normal
to $\Gsh$, pointing to the interior of $\Omega$.
		
\smallskip
\item \label{RegReflSol-PropPrandtl1-1-1} The inequalities hold{\rm :}
\begin{equation}\label{phi-between-in-omega-nonSt-Prandtl}
\max\{\varphi_{1}, \varphi_{2}\}\le\varphi\le \varphi_0  \qquad \mbox{in $\Omega$},
\end{equation}
		
\smallskip
\item\label{RegReflSol-PropPrandtl2}
The monotonicity properties hold:
\begin{equation}\label{MonotoneProperty}
D(\varphi_0-\varphi)\cdot \mathbf{e}_{\mathcal{S}_1}\ge 0,
\quad D(\varphi_0-\varphi)\cdot \mathbf{e}_{\mathcal{S}_0}\le 0 \qquad\,\,\,\, \mbox{in $\Omega$},
\end{equation}
where $\mathbf{e}_{\mathcal{S}_0}$ and $\mathbf{e}_{\mathcal{S}_{1}}$ are the unit
vectors along lines $\mathcal{S}_0$ and $\mathcal{S}_1$
pointing to the positive $\xi_1$--direction, respectively.
\end{enumerate}
\end{definition}

Similar to \eqref{coneOfMonotRegRefl-cone},
the monotonicity properties in \eqref{MonotoneProperty} imply that
\begin{equation}\label{coneOfMonotPrandtlRefl-cone}
  D(\varphi_1-\varphi)\cdot {\bf e}\le 0 \qquad\,\, \mbox{in $\overline\Omega\,\,$ for
  all ${\bf e}\in \overline{Cone(-\mathbf{e}_{\mathcal{S}_1}, \mathbf{e}_{\mathcal{S}_0})}$},
\end{equation}
where $Cone(-\mathbf{e}_{\mathcal{S}_1}, \mathbf{e}_{\mathcal{S}_0})
=\{-a\,\mathbf{e}_{\mathcal{S}_1}+b \,\mathbf{e}_{\mathcal{S}_0}\;: \; a, b>0\}$.
We note that $\mathbf{e}_{\mathcal{S}_0}$ and $\mathbf{e}_{\mathcal{S}_1}$ are not parallel
if $\theta_{\rm w}\ne 0$.
Then we establish the following theorem.

\smallskip
\begin{theorem}\label{mainPrandtlReflThm}
Let $\gamma>1$ and $u_{10}>c_0$.
For any $\theta_{\rm w}\in (0, \theta_{\rm w}^{\rm d})$, there exists a global entropy solution
$\varphi$ of {\rm Problem \ref{fbp-a}} such that the following regularity properties are satisfied for some $\alpha\in(0,1)${\rm :}
\begin{enumerate}[{\rm (i)}]
\item If $\theta_{\rm w}\in (0, \theta_{\rm w}^{\rm s})$,
the reflected shock $\overline{\SzSeg}\cup\Shock\cup\overline{\SoSeg}$ is $C^{2,\alpha}$--smooth,
and $\varphi\in C^{1,\alpha}(\overline\Omega)\cap C^\infty(\overline\Omega\setminus (\overline{\Gamma_{\rm sonic}^1}\cup
\overline{\Gamma_{\rm sonic}^2}))$;

\item If  $\theta_{\rm w}\in [\theta_{\rm w}^{\rm s}, \theta_{\rm w}^{\rm d})$,
the reflected shock $\overline{\Shock}\cup\overline{\SoSeg}$ is $C^{1,\alpha}$ near $O$
and $C^{2,\alpha}$ away from $O$, and
$\varphi\in C^{1,\alpha}(\overline\Omega)\cap C^\infty(\overline\Omega\setminus (\{O\}\cup\overline{\Gamma_{\rm sonic}^2}))$.
\end{enumerate}
Moreover, in both cases, $\varphi$ is $C^{1,1}$ across the sonic arcs, and
$\Shock$ is $C^\infty$ in its relative interior.

Furthermore, $\varphi$ is an admissible solution in the sense of Definition {\rm \ref{admisSolnDef-Prandtl}},
so $\varphi$ satisfies further properties listed in
Definition {\rm \ref{admisSolnDef-Prandtl}}.
\end{theorem}

We follow the argument described in \S \ref{RegReflProblSect} so that, for any small $\delta>0$,
we obtain the required uniform estimates of admissible solutions with wedge angles
$\theta_{\rm w} \in [0, \theta_{\rm w}^{\rm d}-\delta]$. Using these estimates,
we apply the Leray-Schauder degree theory
to obtain the existence for each $\theta_{\rm w} \in [0, \theta_{\rm w}^{\rm d}-\delta]$
in the class of admissible solutions,
starting from the unique normal solution for $\theta_{\rm w}=0$.
Since $\delta>0$ is arbitrary,
the existence of a global entropy solution for any $\theta_{\rm w}\in(0,\theta_{\rm w}^{\rm d})$
can be established.
More details can be found in Bae-Chen-Feldman \cite{BCF-14}; see also
Chen-Feldman \cite{CF-book2018}.

\smallskip
The existence results in Bae-Chen-Feldman \cite{BCF-14}
indicate that the steady weak supersonic/transonic shock solutions are the asymptotic limits
of the dynamic self-similar solutions, the Prandtl-Meyer reflection configurations,
in the sense of \eqref{1-k-c} in {Problem 4.5} for all $\theta_{\rm w}\in (0, \theta_{\rm w}^{\rm d})$
and all $\gamma>1$.

On the other hand, it is shown in Elling \cite{Elling} and  Bae-Chen-Feldman \cite{BCF-14} that,
for each $\gam >1$, there is no self-similar {\it strong} Prandtl-Meyer reflection configuration
for the unsteady potential flow in the class of admissible solutions.
This means that the situation for the dynamic stability of the strong steady oblique shocks
is more sensitive.

\section{Convexity of Self-Similar Transonic Shocks and Free Boundaries}

We now discuss some recent developments in the analysis of geometric properties of transonic shocks
as free boundaries in the 2-D self-similar coordinates for compressible fluid flows.
In Chen-Feldman-Xiang \cite{ChenFeldmanXiang}, we have developed
a general framework for the analysis of the convexity of
transonic shocks as free boundaries. For both applications discussed above,
the von Neumann problem for shock reflection-diffraction in \S \ref{RegReflProblSect}
and the Prandtl-Meyer problem for unsteady supersonic flow onto solid wedges
in \S\ref{PrandtlMeyerProblSect},
the admissible solutions satisfy the conditions
of this abstract framework, as shown in \cite{ChenFeldmanXiang}. For simplicity,
we present below the results on the convexity properties of transonic shocks for these two problems
(without discussion on the abstract framework).

For the regular shock reflection-diffraction
configurations, we recall that, for admissible solutions  in the sense of
Definition {\rm \ref{admisSolnDef}},
the inequality in \eqref{coneOfMonotRegRefl-cone} is shown to be strict
for any ${\bf e}\in Cone({\bf e}_{\xi_2}, {\bf e}_{\CS_1})$.
From this, it is proved that, for admissible solutions, the shock is a graph in the coordinate system $(S,T)$
with respect to basis $\{\ee, \ee^{\perp}\}$ for any unit vector
$\ee\in Cone({\bf e}_{\xi_2}, {\bf e}_{\CS_1})$,
where $\ee^{\perp}$ is the unit vector orthogonal to $\ee$ and oriented so that $T_{P_1}>T_{P_2}$,
and we have used notation $(S_P,T_P)$ for the coordinates of point $P$.
That is, there exists  $f_\ee\in C^\infty( (T_{P_2}, T_{P_1}))\cap C^1([T_{P_2}, T_{P_1}])$ such that
\begin{equation}\label{shock-graph-inMainThm}
\Gsh=\{(S,T):\; S=f_\ee(T), \; T_{P_2}<T<T_{P_1}\},\quad\;
\Omega\cap\{ T_{P_2}<T<T_{P_1}\} \subset\{S<f_\ee(T)\},
\end{equation}
where we have used the notational convention \eqref{P1-P4-P0-subs-Ch2}
in the subsonic/sonic case.

Since the convexity or concavity
of a shock as a graph depends on the orientation of the coordinate system
and $\Omega$ will be shown to be a convex domain (corresponding to
the concavity of $f_\ee$ in \eqref{shock-graph-inMainThm}),
we do not distinguish them
and instead use the term {\it convexity} for either case below.
Then we have

\begin{theorem}[Convexity of transonic shocks for the regular shock reflection-diffraction configurations]\label{thm:con}
If a solution of the von Neumann problem for shock reflection-diffraction is admissible in the sense of
Definition {\rm \ref{admisSolnDef}},
then its domain $\Omega$ is convex, and the shock curve $\Gsh$ is a strictly convex graph
in the following sense{\rm :}
For any
$\ee\in Cone({\bf e}_{\xi_2}, {\bf e}_{\CS_1})$,
the function $f_\ee$
in \eqref{shock-graph-inMainThm} satisfies
$$
f''_\ee<0 \qquad \mbox{ on $( T_{P_2}, T_{P_1})$}.
$$
That is, $\Gsh$ is uniformly convex on any closed subset of its relative interior.

Moreover, for the solution of {\rm Problem \ref{fbp-c}} extended to $\Lambda$
by \eqref{phi-states-0-1-2-MainThm}, with pseudo-potential
$\varphi\in C^{0,1}(\Lambda)$ satisfying
Definition {\rm \ref{admisSolnDef}}{\rm (\ref{RegReflSol-Prop0})}--\eqref{RegReflSol-Prop1-1-1},
the shock is strictly convex if and only if
Definition {\rm \ref{admisSolnDef}}\eqref{RegReflSol-Prop2} holds.
\end{theorem}

For the Prandtl-Meyer problem for unsteady supersonic flow onto solid wedges,
the results are similar.
We first note that, based on \eqref{coneOfMonotPrandtlRefl-cone}
(which is strict for ${\bf e}\in {Cone(-\mathbf{e}_{\mathcal{S}_1}, \mathbf{e}_{\mathcal{S}_0})}$)
and the maximum principle,
it is proved that, for admissible solutions in the sense of Definition {\rm \ref{admisSolnDef-Prandtl}},
the shock is a graph in the coordinate system $(S,T)$  with respect to basis $\{\ee, \ee^{\perp}\}$
for any unit vector $\ee\in Cone(-\mathbf{e}_{\mathcal{S}_1}, \mathbf{e}_{\mathcal{S}_0})$,
{\it i.e.}, \eqref{shock-graph-inMainThm} holds,
with $f_\ee\in C^\infty(( T_{P_2}, T_{P_1}))\cap C^1([T_{P_2}, T_{P_1}])$,
where we have used the notational convention $P_1=P_0$
for the subsonic/sonic case
$\theta_{\rm w}\in [\theta_{\rm w}^{\rm s}, \theta_{\rm w}^{\rm d})$.

\begin{theorem}[Convexity of transonic shocks for the Prandtl-Meyer reflection configurations]\label{thm:con-Prandtl}
If a solution of the Prandtl-Meyer problem is admissible in the sense of
Definition {\rm \ref{admisSolnDef-Prandtl}},
then its domain $\Omega$ is convex, and the shock curve $\Gsh$ is a strictly convex graph in the following sense{\rm :}
For any $\ee\in Cone(-\mathbf{e}_{\mathcal{S}_1}, \mathbf{e}_{\mathcal{S}_0})$,
the function $f_\ee$ in \eqref{shock-graph-inMainThm} satisfies
$$
f''_\ee<0 \qquad \mbox{ on $( T_{P_2}, T_{P_1})$}.
$$
That is, $\Gsh$ is uniformly convex on any closed subset of its relative interior.

Moreover, for the solution of {\rm Problem \ref{fbp-a}} extended to $\Lambda$
by \eqref{extsol} $($with the appropriate modification for the subsonic/sonic case$)$ with pseudo-potential
$\varphi\in C^{0,1}(\Lambda)$ satisfying
Definition {\rm \ref{admisSolnDef-Prandtl}}{\rm (\ref{RegReflSol-PropPrandtl0})}--\eqref{RegReflSol-PropPrandtl1-1-1},
the shock is strictly convex if and only if
Definition {\rm \ref{admisSolnDef-Prandtl}}\eqref{RegReflSol-PropPrandtl2} holds.
\end{theorem}

Theorems {\rm \ref{thm:con}}--{\rm \ref{thm:con-Prandtl}}
indicate that the  curvature of $\Gsh${\rm :}
$$
\kappa=-\frac{f''_\ee(T)}{\big(1+(f'_\ee(T))^2\big)^{3/2}}
$$
has a positive lower bound
on any closed subset of $(T_{P_2}, T_{P_1})$.

Now we discuss the techniques developed in \cite{ChenFeldmanXiang}
by giving the main steps in the proofs of Theorems \ref{thm:con}--\ref{thm:con-Prandtl}.
While the argument in \cite{ChenFeldmanXiang} is carried out in a more general setting,
we focus here on the specific cases of the regular shock reflection-diffraction and Prandtl-Meyer reflection configurations;
see  \cite{ChenFeldmanXiang} for the results in the more general setting and the detailed proofs.

For the von Neumann problem, define
$$
\phi:=\varphi-\varphi_1.
$$
For the Prandtl-Meyer problem, define
$$
\phi:=\varphi-\varphi_0.
$$
Then, in both cases, $\phi=0$ on $\Shock$.
From this, using Definition \ref{admisSolnDef}\eqref{RegReflSol-Prop1-1}
for the regular reflection-diffraction case
and Definition \ref{admisSolnDef-Prandtl}\eqref{RegReflSol-PropPrandtl1-1}
for the Prandtl-Meyer reflection case,
it follows that, in both problems, $\phi<0$ in $\Omega$ near $\Gsh$.
Since $\Gsh$ is the zero level set of $\phi$,
the conclusion of Theorems {\rm \ref{thm:con}}--{\rm \ref{thm:con-Prandtl}}
on the strict convexity of $\Gsh$ is equivalent
to the following{\rm :}
$\phi_{\bt\bt}>0$ along $\Gsh^0$, where $\Gsh^0$ is the
relative interior
of $\Gsh$.
Also, denote by $Con$ the cone from \eqref{coneOfMonotRegRefl-cone}
for the von Neuamnn problem and the cone from \eqref{coneOfMonotPrandtlRefl-cone}
for the Prandtl-Meyer problem.

First, we establish the relation between the strict convexity/concavity of a portion of the shock
and the possibility for $\partial_\ee \phi$, with $\ee\in Con$, to attain its
local  minimum or maximum with respect to
$\overline\Omega$ on that portion of the shock.
More precisely,
on a portion of ``wrong" convexity on which $f''_\ee \ge 0$ (equivalently,
$\phi_{\bt\bt}\le 0$), $\phi_\ee$  cannot attain its local minimum relative to $\overline\Omega$.
Then, assuming that a portion of the free boundary has a ``wrong" convexity $f''_\ee>0$,
we show that  $\phi_\ee$ for $\ee\in Con$ attains
its local minimum relative to $\Shock$ on the closure of that portion.
As we discussed above, it cannot be a local minimum with respect to $\overline\Omega$.
Starting from that, through a nonlocal argument, with the use of the
maximum principle for equation \eqref{equ:study},
considered as a linear elliptic PDE for $\phi$, in $\Omega$, and
the boundary
conditions on various parts of $\partial\Omega$, we reach a contradiction,
which implies that the shock is  convex, possibly non-strictly, {\it i.e.}, $f''_\ee\le 0$ on $( T_{P_2}, T_{P_1})$,
or equivalently, $\phi_{\bt\bt}\ge 0$ on $\Shock$.
Extending the previous argument
with use of the real analyticity of $\Shock^0$,
we improve the result to the locally uniform convexity as in Theorems
{\rm \ref{thm:con}}--{\rm \ref{thm:con-Prandtl}}.

\smallskip
Furthermore, with the convexity of reflected-diffracted transonic shocks,
the uniqueness and stability of global regular shock reflection-diffraction
configurations have also been established  in
the class of {\em admissible solutions}; see Chen-Feldman-Xiang \cite{CFX-Unique} for the details.

\medskip
The nonlinear method, ideas, techniques, and approaches that we have presented above for
solving M-D transonic shocks and free boundary problems should be useful to analyze
other longstanding and newly emerging problems.
Examples of such problems include the unsolved M-D steady transonic shock problems
for the full Euler equations (including steady detached shock problems),
the unsolved M-D self-similar transonic shock problems
(such as the 2-D Riemann problems and the conic body problems) for potential flow,
as well as the longstanding open transonic shock problems
for both the isentropic and the full Euler equations; also see Chen-Feldman \cite{CF-book2018}.
Certainly, further new ideas, techniques, and methods are still required to be developed in order
to solve these mathematically challenging and fundamentally important problems.

\bigskip
\medskip
\noindent
{\bf Acknowledgements}.
The research of Gui-Qiang G. Chen was supported in part by
the UK Engineering and Physical Sciences Research Council Awards
EP/L015811/1, EP/V008854, and EP/V051121/1, and the Royal Society--Wolfson Research Merit Award WM090014.
The research of Mikhail Feldman was
supported in part by the National Science Foundation under DMS-1764278 and
DMS-2054689.


\begin{thebibliography}{10}
\bibitem{AC}
H.~W. Alt and L.~A. Caffarelli,
\newblock Existence and regularity for a minimum problem with free boundary,
\newblock {\em J. Reine Angew. Math.} {\bf 325} (1981), 105--144.

\bibitem{AltCafFried_Compres-83}
H.~W. Alt and L.~A. Caffarelli, and A. Friedman,
\newblock
Axially symmetric jet flows,
{\em Arch. Ration. Mech. Anal.} {\bf 81} (1983), 97--149.


\bibitem{ACF}
H.~W. Alt, L.~A. Caffarelli, and A. Friedman,
\newblock A free boundary problem for quasilinear elliptic equations,
\newblock {\em Ann. Scuola Norm. Sup. Pisa Cl. Sci. (4)},
{\bf 11} (1984), 1--44.



\bibitem{AltCafFried_Compres}
H.~W. Alt, L.~A. Caffarelli, and A. Friedman,
\newblock Compressible flows of jets and cavities,
\newblock {\itshape J. Diff. Eqs.} \textbf{56} (1985), 82--141.



\bibitem{BaeChenFeldman}
M. Bae, G.-Q. Chen, and M. Feldman,
Regularity of solutions to regular shock reflection for potential flow,
{\itshape Invent. Math.}  \textbf{175} (2009), 505--543.

\bibitem{BCF-14}
M. Bae, G.-Q. Chen, and M. Feldman,\,
\textit{Prandtl-Meyer Reflection Configurations,
Transonic Shocks, and Free Boundary Problems},
Research Monograph, 118 pages,
Memoirs of Amer. Math. Soc.,
AMS: Providence, RI, 2022 (to appear).

\bibitem{BaeFeldman}
M. Bae and M. Feldman,
Transonic shocks in multidimensional divergent nozzles,
{\itshape Arch. Rational Mech. Anal.} \textbf{201} (2011), 777--840.


\bibitem{BD}
G. Ben-Dor,
\newblock{\itshape Shock Wave Reflection Phenomena},
\newblock{Springer-Verlag: New York}, 1991.

\bibitem{Busemann}
A. Busemann,
\newblock
Gasdynamik.
Handbuch der Experimentalphysik, Vol. IV, Akademische Verlagsgesellschaft,
Leipzig, 1931.

\bibitem{Ca1}
L.~A. Caffarelli,
\newblock A Harnack inequality approach to the regularity
of free boundaries, I: Lipschitz free boundaries are $C\sp
{1,\alpha}$, {\itshape Rev.  Mat. Iberoamericana}, \textbf{3} (1987) 139--162.

\bibitem{Ca2}
L.~A. Caffarelli,
\newblock A Harnack inequality approach to the regularity
of free boundaries, II: Flat free boundaries are Lipschitz,
{\itshape Comm. Pure
Appl. Math.} \textbf{42} (1989), 55--78.

\bibitem{Ca3}
L.~A. Caffarelli,
\newblock A Harnack inequality approach to the regularity
of free boundaries, III: Existence theory,
compactness, and dependence on $X$,
{\itshape Ann. Scuola Norm. Sup. Pisa
Cl. Sci.} (4), \textbf{15} (1989), 583--602.

\bibitem{CJK}
L.~A. Caffarelli, D. Jerison, C. Kenig,
{Some new
monotonicity theorems with applications to free boundary problems},
\textit{Ann. of Math.} \textbf{155} (2002), 369--404.

\bibitem{Caff-Salsa}
L.~A. Caffarelli and S. Salsa,
\textit{A Geometric Approach to Free Boundary Problems},
American Mathematical Society: Providence, RI, 2005.


\bibitem{CKK1b}
S. Cani\'{c}, B.~L. Keyfitz,  and E.~H. Kim,
\newblock A free boundary problem for a quasilinear degenerate elliptic
equation: regular reflection of weak shocks,
\newblock {\itshape  Comm. Pure Appl. Math.} \textbf{55} (2002), 71--92.

\bibitem{CanicKeyfitz}
S. Cani\'{c}, B.~L. Keyfitz, and G.~M. Lieberman,
\newblock A proof of existence of perturbed steady transonic
shocks via a free boundary problem,
\newblock {\itshape Comm. Pure Appl. Math.} \textbf{53} (2002), 484--511.

\bibitem{CCY2}
T. Chang, G.-Q. Chen, and S. Yang,
On the Riemann problem for two-dimensional Euler equations
I:     Interaction of shocks and rarefaction waves,
   \textit{Discrete Contin. Dynam. Systems}, \textbf{1} (1995), 555--584.

\bibitem{CCY3}
T. Chang, G.-Q. Chen, and S. Yang,
{On the Riemann problem for two-dimensional Euler equations
     II: Interaction of contact discontinuities},
   \textit{Discrete Contin. Dynam. Systems}, \textbf{6} (2000), 419--430.

\bibitem{CH}
T. Chang and L. Hsiao,
\textit{The Riemann Problem and Interaction of Waves in Gas
Dynamics}, Longman Scientific \& Technical: Harlow; and John Wiley
\& Sons, Inc.: New York, 1989.


\bibitem{Chapman} C.~J. Chapman,
\textit{High Speed Flow}, Cambridge University Press: Cambridge, 2000.


\bibitem{Chen}
G.-Q. Chen,
{\it Euler Equations and Related Hyperbolic Conservation Laws},
Chapter 1,
\textit{Handbook of Differential Equations,
Evolutionary Equations}, Vol. \textbf{2}, Eds. C. M. Dafermos and E. Feireisl,
Elsevier: Amsterdam, The Netherlands, 2005.

\bibitem{Chen2}
G.-Q. Chen,
{Supersonic flow onto solid wedges, multidimensional shock waves, and free boundary problems},
{\it Science China Mathematics}, {\bf 60 (8)} (2017), 1353--1370.


\bibitem{CCF}
G.-Q. Chen, J. Chen, and M. Feldman,
\newblock
Transonic shocks and free
boundary problems for the full Euler equations in infinite nozzles,
\newblock {\it J. Math. Pures Appl.} (9), {\bf 88} (2007), 191--218.


\bibitem{Chen-Chen-Feldman}
G.-Q. Chen, J. Chen, and M. Feldman,
Transonic flows with shocks past curved wedges for the full Euler equations,
\textit{Discrete Contin. Dyn. Syst.} \textbf{36} (2016), 4179--4211.


\bibitem{CCF3}
G.-Q. Chen, J. Chen, and M. Feldman,
\newblock
Stability and asymptotic behavior of transonic flows past wedges for the full Euler equations.
\textit{Interfaces and Free Boundaries}, \textbf{19} (2017), 591--626.

\bibitem{CCX}
G.-Q. Chen, J. Chen, and W. Xiang,
\newblock
Stability of attached transonic shocks in steady potential flow past three-dimensional wedges,
\textit{Commun. Math. Phys.} \textbf{387} (2021), 111--138.

\bibitem{C-Fang}
G.-Q. Chen and B.-X. Fang,
{Stability of transonic shock-fronts in steady potential flow
      past a perturbed cone},
      \textit{Discrete Conti. Dyn. Syst.} \textbf{23} (2009), 85--114.

\bibitem{C-Fang-2}
G.-Q. Chen and B.-X. Fang,
Stability of transonic shocks in steady supersonic flow past
multidimensional wedges, \textit{Adv. Math.} \textbf{314} (2017), 493--539.


\bibitem{CF-JAMS2003}
G.-Q. Chen and M. Feldman,
\newblock Multidimensional transonic shocks and free boundary problems for
  nonlinear equations of mixed type,
\newblock  {\itshape J. Amer. Math. Soc.} \textbf{16} (2003), 461--494.

\bibitem{ChenFeldman2}
G.-Q. Chen and M. Feldman,
Steady transonic shocks and free boundary problems in infinite
cylinders for the Euler equations,
\newblock {\itshape Comm. Pure Appl. Math.}
\textbf{\bf 57} (2004), 310--356.

\bibitem{ChenFeldman4}
G.-Q. Chen and M. Feldman,
Free boundary problems and transonic shocks for the Euler equations
in unbounded domains, {\itshape Ann. Scuola Norm. Sup. Pisa Cl. Sci.
{\rm (5)}}, \textbf{3}(2004), 827--869.


\bibitem{ChenFeldman3Arch}
G.-Q. Chen and M. Feldman,
Existence and stability of multidimensional transonic flows through
an infinite nozzle of arbitrary cross-sections,
{\itshape Arch. Ration. Mech. Anal.} \textbf{184} (2007), 185--242.


\bibitem{ChenFeldman}
G.-Q. Chen and M. Feldman,
Global solutions to shock reflection by large-angle wedges for potential flow,
{\itshape Ann. of Math.} \textbf{171} (2010), 1019--1134.

\bibitem{ChenFeldman4b}
G.-Q. Chen and M. Feldman,
Comparison principles for self-similar potential flow,
\textit{Proc. Amer. Math. Soc.} \textbf{140} (2012), 651--663.

\bibitem{CF-book2018}
G.-Q. Chen and M. Feldman,
\newblock {\em Mathematics of Shock Reflection-Diffraction and von Neumann's Conjecture}.
\newblock Research Monograph, Annals of Mathematics Studies, \textbf{197}, Princeton University Press, Princetion, 2018.

\bibitem{ChenFeldmanXiang}
G.-Q. Chen, M. Feldman, and W. Xiang,
Convexity of self-similar transonic shock waves for potential flow,
\textit{Arch. Ration. Mech. Anal.} \textbf{238} (2020), 47--124.

\bibitem{CFX-Unique}
G.-Q. Chen, M. Feldman, and W. Xiang,
\newblock Uniqueness of regular shock reflection/diffraction configurations
for potential flow,
\newblock
\newblock Preprint 2021.

\bibitem{CKZ21}
G.-Q. Chen, J. Kuang, and Y. Zhang,
\newblock Stability of conical shocks in the three-dimensional steady supersonic isothermal flows
past Lipschitz perturbed cones,
\textit{SIAM J. Math. Anal.} \textbf{53} (20210, 2811--2862.


\bibitem{ChenLi2008}
G.-Q. Chen and T.-H. Li,
Well-posedness for two-dimensional steady supersonic Euler flows past a Lipschitz wedge,
\textit{J. Diff. Eqs.} \textbf{244} (2008), 1521--1550.


\bibitem{CSV}
G.-Q. Chen, H. Shahgholian, and J.-V. V\'{a}zquez,
Free boundary problems: The forefront of current and future developments,
In: \textit{Free Boundary Problems and Related Topics}.
Theme Volume: Phil. Trans. R. Soc. \textbf{A373}: 20140285,
The Royal Society: London, 2015.


\bibitem{CY}
G.-Q. Chen and H. Yuan,
Uniqueness of transonic shock
solutions in a duct for steady potential flow,
\textit{J. Diff. Eqs.} \textbf{247} (2009), 564--573.

\bibitem{CY2}
G.-Q. Chen and H. Yuan,
Local uniqueness of steady spherical transonic shock-fronts for the three-dimensional full Euler equations,
\textit{Comm. Pure Appl. Anal.} \textbf{12} (2013), 2515--2542.


\bibitem{ChenZhangZhu}
G.-Q. Chen, Y. Zhang, and D. Zhu,
Existence and stability of supersonic Euler flows past
Lipschitz wedges,
\textit{Arch. Ration. Mech. Anal.} \textbf{181} (2006), 261--310.

\bibitem{CCJ2011}
J. Chen, C. Christoforou C, and K. Jegdi\'{c},
\newblock Existence and uniqueness analysis of a detached shock problem for the potential flow.
\textit{Nonlinear Anal.}  \textbf{74} (2011), 705--720.



\bibitem{ChS3}
S.-X. Chen,
Global existence of supersonic flow past a curved convex wedge,
\textit{J. Partial Diff. Eqs.}
\textbf{11} (1998), 73--82.

\bibitem{Sxchen3}
S.-X. Chen,
\newblock
Stability of transonic shock fronts in two-dimensional Euler
systems,
\textit{Trans. Amer. Math. Soc.} \textbf{357} (2005), 287--308.

\bibitem{Sxchen-book2020}
S.-X. Chen,
{\it Mathematical Analysis of Shock Wave Reflection},
Series in Contemporary Mathematics 4,
Shanghai Scientific and Technical Publishers, China;
Springer Nature Singapore Pte Ltd., Singapore,
2020.



\bibitem{ChenFang}
S.-X. Chen and B. Fang,
{Stability of transonic shocks in supersonic flow past a
wedge}, \textit{J. Diff. Eqs.} \textbf{233} (2007), 105--135.


\bibitem{CXY}
S.-X. Chen, Z. Xin, and H. Yin,
{Global shock waves for the supersonic flow past a perturbed
cone}, \textit{Commun. Math. Phys.} \textbf{228} (2002), 47--84.

\bibitem{CDK}
E. Chiodaroli, C. De Lellis, and O. Kreml,
Global ill-posedness of the isentropic system of gas dynamics,
\textit{Comm. Pure Appl. Math.} \textbf{68} (2015), 1157--1190.


\bibitem{CC}
J.~D. Cole and Cook, L.~P. Cook,
{\itshape Transonic Aerodynamics},
North-Holland: Amsterdam, 1986.

\bibitem{CF}
R. Courant  and  K.~O. Friedrichs,
{\itshape Supersonic Flow and Shock Waves},
Springer-Verlag: New York, 1948.

\bibitem{Da}
C.~M. Dafermos,
\textit{Hyperbolic Conservation Laws in Continuum Physics},
4th Ed., Springer-Verlag: Berlin, 2016.

\bibitem{DH}
P. Daskalopoulos and R. Hamilton,
The free boundary in the Gauss curvature flow with flat sides,
\newblock {\itshape J. Reine Angew. Math.}  {\bf 510} (1999), 187--227.

\bibitem{Elling}
V. Elling,
Non-existence of strong regular reflections in self-similar potential flow,
\textit{J. Diff. Eqs.} \textbf{252} (2012), 2085--2103.

\bibitem{EllingLiu1}
V. Elling and T.-P. Liu,
{Supersonic flow onto a solid wedge},
\textit{Comm. Pure Appl. Math.}
\textbf{61} (2008), 1347--1448.

\bibitem{Fang} B.-X. Fang,
{Stability of transonic shocks for the full Euler system in
supersonic flow past a wedge},
\textit{Math. Methods Appl. Sci.} \textbf{29} (2006), 1--26.

\bibitem{Fang-Liu-Yuan}
B.-X. Fang, L. Liu, and H. R. Yuan,
Global uniqueness of transonic shocks in two-dimensional steady compressible Euler flows,
\textit{Arch. Ration. Mech. Anal.} \textbf{207} (2013), 317--345.


\bibitem{FT} C. Ferrari and  F.~G. Tricomi,
\textit{Transonic Aerodynamics},
Academic Press: New York,  English Transl. of
\textit{Aerodinamica Transonica}, Cremonese (1962).


\bibitem{Friedman}
A. Friedman,
\textit{Variational Principles and Free-Boundary Problems},
2nd Ed., Robert E. Krieger Publishing Co., Inc.:
Malabar, FL, 1988
 [First edition, John Wiley \& Sons, Inc.: New York, 1982].

\bibitem{GilbargTrudinger}
D. Gilbarg  and N. Trudinger,
\newblock {\itshape Elliptic Partial Differential Equations of Second Order,}
\newblock 2nd Edition, Springer-Verlag: Berlin, 1983.

\bibitem{GlimmK}
J. Glimm, C.  Klingenberg, O.  McBryan, B. Plohr, D. Sharp, and S. Yaniv,
Front tracking and two-dimensional Riemann problems,
{\itshape Adv. Appl. Math.} \textbf{6} (1985), 259--290.


\bibitem{GlimmMajda}
J. Glimm and A. Majda,
{\itshape Multidimensional Hyperbolic Problems and Computations},
Springer-Verlag: New York, 1991.

\bibitem{Gu}  C.-H. Gu,
{A method for solving the supersonic flow past a curved wedge
(in Chinese)}, \textit{Fudan Univ. J.} \textbf{7} (1962), 11--14.

\bibitem{Guderley}
K.~G. Guderley,
{\itshape The Theory of Transonic Flow},
Translated from the German by J. R. Moszynski, Pergamon Press:
Oxford-London-Paris-Frankfurt; Addison-Wesley Publishing Co. Inc.:
Reading, Mass., 1962.

\bibitem{Hadamard}
J. Hadamard,
{\itshape Le\c{c}ons sur la Propagation des Ondes et
les \'{E}quations de l'Hydrodynamique}, Hermann: Paris, 1903
(Reprinted by Chelsea 1949).


\bibitem{Harabetian}
E. Harabetian,
Diffraction of a weak shock by a wedge,
{\itshape Comm. Pure Appl. Math.} \textbf{40} (1987), 849--863.

\bibitem{HK}
J.~K. Hunter and J.~B. Keller,
Weak shock diffraction, {\itshape Wave Motion}, \, \textbf{6} (1984), 79--89.


\bibitem{KB}
J.~B. Keller and A.~A. Blank,
Diffraction and reflection of pulses by wedges and corners, {\itshape
Comm. Pure Appl. Math.} \textbf{4} (1951), 75--94.

\bibitem{KinderlehrerNirenberg}
D. Kinderlehrer and L. Nirenberg,
Regularity in free boundary problems,
\newblock {\itshape Ann. Scuola Norm. Sup. Pisa Cl. Sci. {\rm (}4{\rm )}} \,
\textbf{4} (1977), 373--391.

\bibitem{KKMM}
C. Klingenberg, O. Kreml, V. M\'{a}cha, and S. Markfelder,
Shocks make the Riemann problem for the full Euler system
in multiple space dimensions ill-posed,
\newblock {\itshape Nonlinearity}, \textbf{33} (2020), 6517--6540.






\bibitem{KTa}
A. Kurganov and E. Tadmor,
{Solution of two-dimensional Riemann problems for gas dynamics
without Riemann problem solvers}, \textit{Numer. Methods Partial Diff. Eqs.}
\textbf{18} (2002), 584--608.

\bibitem{Lax}
P.~D. Lax,
{\itshape Hyperbolic Systems of Conservation Laws and the Mathematical
Theory of Shock Waves}, CBMS-RCSM, SIAM: Philiadelphia, 1973.

\bibitem{Lax2003}
P.~D. Lax,
{\itshape Max Shiffman} (1914-2000),
Notices of Amer. Math. Soc., pp. 1401, December 2003.

\bibitem{LaxLiu}
P.~D. Lax and X.-D. Liu,
Solution of two-dimensional Riemann problems of gas dynamics
by positive schemes,
\newblock {\itshape SIAM J. Sci. Comput.} \textbf{19} (1998), 319--340.



\bibitem{LXY}
J. Li, Z. Xin, and H. Yin,
Transonic shocks for the full compressible Euler system in a general two-dimensional
de Laval nozzle,
\textit{Arch. Ration. Mech. Anal.} \textbf{207} (2003), 533--581.

\bibitem{LXY2}
L. Li, G. Xu, and H.~C. Yin,
On the instability problem of a 3-D transonic oblique shock wave.
\textit{Adv. Math.} \textbf{282} (2015), 443--515.

\bibitem{LZY}
J. Li, T. Zhang, and S. Yang,
{\em The Two-Dimensional Riemann
Problem in Gas Dynamics}, Longman (Pitman Monographs 98): Essex, 1998.

\bibitem{Li}
T.-T. Li,
{On a free boundary problem}, \textit{Chinese Ann. Math.} \textbf{1} (1980), 351--358.


\bibitem{Lieberman86}
G.~M. Lieberman,
Regularity of solutions of nonlinear elliptic boundary value problems,
\newblock {\itshape J. Reine Angew. Math.}  {\bf 369} (1986), 1--13.

\bibitem{LiebermanTrudinger}
G.~M. Lieberman and N.~S. Trudinger,
Nonlinear oblique boundary value problems for nonlinear elliptic equations,
{\itshape Trans. Amer. Math. Soc.} \textbf{295} (1986), 509--546.

\bibitem{Lie}
W.-C. Lien and T.-P. Liu,
{Nonlinear stability of a self-similar 3-dimensional gas flow},
\textit{Commun. Math. Phys.} \textbf{204} (1999), 525--549.


\bibitem{Lighthill1}
M.~J. Lighthill,
The diffraction of a blast I,
{\itshape Proc. Roy. Soc. London} \, \textbf{198A} (1949), 454--470.

\bibitem{Lighthill2}
M.~J. Lighthill,
The diffraction of a blast II,
{\itshape Proc. Roy. Soc. London} \, \textbf{200A} (1950), 554--565.

\bibitem{LW}  F.~H. Lin and L. Wang,
A class of fully nonlinear elliptic equations with singularity at the boundary,
\newblock {\itshape J. Geom. Anal.} \textbf{8} (1998), 583--598.

\bibitem{Liu-Yuan}
L. Liu and H. R. Yuan,
\newblock Stability of cylindrical transonic shocks for the two-dimensional steady
compressible Euler system,
\textit{J. Hyperbolic Differ. Equ.} \textbf{5} (2008), 347--379.


\bibitem{Liu-Xu-Yuan}
L. Liu, G. Xu, and H. R. Yuan,
\newblock Stability of spherically symmetric subsonic flows and transonic
shocks under multidimensional perturbations,
\textit{Adv. Math.} \textbf{291} (2016), 696--757.



\bibitem{Liu}
T.-P. Liu,
\newblock Multi-dimensional gas flow: some historical perespectives.
\textit{Bull Inst Math Acad Sinica $($New Series$)$}, \textbf{6} (2011), 269--291



\bibitem{Mach}  E. Mach,
\"{U}ber den verlauf von funkenwellenin der ebene und im raume,
{\itshape Sitzungsber. Akad. Wiss. Wien} \, \textbf{\bf 78} (1878), 819--838.

\bibitem{Majda}
A. Majda,
{\em Compressible Fluid Flow and Systems of Conservation
 Laws in Several Space Variables},
Springer-Verlag: New York, 1984.


\bibitem{Meyer}
Th. Meyer,
\newblock
\"{U}ber zweidimensionale Bewegungsvorg\"{a}nge in einem Gas,
das mit \"{U}berschallgeschwindigkeit str\"{o}mt.
Dissertation, G\"{o}ttingen, 1908. {\em Forschungsheft des Vereins deutscher Ingenieure},
Vol. 62,  pp. 31--67, Berlin, 1908

\bibitem{Mises}
R.~V. Mises,
{\em Mathematical Theory of Compressible Fluid Flow},
Academic Press: New York, 1958.


\bibitem{Morawetz2}
C.~S. Morawetz,
Potential theory for regular and Mach reflection of a shock at a
wedge,
\newblock {\itshape Comm. Pure Appl. Math.} \textbf{47} (1994), 593--624.


\bibitem{Prandtl}
L. Prandtl,
\newblock Allgemeine \"{U}berlegungen \"{u}ber die Str\"{o}mung zusammendr\"{u}ckbarer
Flu\"{u}ssigkeiten. Zeitschrift f\"{u}r angewandte Mathematik und Mechanik,
\textbf{16} (1938), 129--142

\bibitem{Riemann}
B. Riemann,
\"{U}ber die Fortpflanzung ebener Luftvellen von endlicher
Schwingungsweite,
\textit{G\"{o}tt. Abh. Math. Cl.} \textbf{8} (1860), 43--65.

\bibitem{Schaeffer}
D.~G. Schaeffer,
{Supersonic flow past a nearly straight wedge}, \textit{Duke Math. J.}
\textbf{43} (1976), 637--670.

\bibitem{SCG}
C.~W. Schulz-Rinne, J.~P. Collins, and  H.~M. Glaz,
{Numerical solution of the Riemann problem for two-dimensional
gas dynamics}, \textit{SIAM J. Sci. Comput.} \textbf{14} (1993), 1394--1414.


\bibitem{Serre}
D. Serre,
Shock reflection in gas dynamics. In: {\itshape Handbook of
Mathematical Fluid Dynamics}, Vol. {\bf 4}, pp. 39--122,
Elsevier: North-Holland, 2007.

\bibitem{Serre2}
D. Serre,
von Neumann's comments about existence and uniqueness for the initial-boundary value problem in gas dynamics,
\textit{Bull. Amer. Math. Soc. $($N.S.$)$}, \textbf{47} (2010), 139--144.

\bibitem{Sh}
M. Shiffman,
On the existence of subsonic flows of a compressible fluid,
{\itshape J. Ration. Mech. Anal.} \textbf{1} (1952), 605--652.

\bibitem{Stokes}
G.~G. Stokes,
On a difficulty in the theory of sound,
\textit{Philos. Magazine, Ser. 3}, \textbf{33} (1845),
349--356.

\bibitem{Trudinger85}
N. Trudinger,
On an interpolation inequality and its applications
to nonlinear elliptic equations,
{\itshape Proc.  Amer. Math. Soc.} \textbf{95} (1985), 73--78.

\bibitem{VD}
M. Van Dyke,
{\itshape An Album of Fluid Motion},
The Parabolic Press: Stanford, 1982.

\bibitem{Neumann0}
J. von Neumann,
Theory of shock waves,
{\itshape Progress Report},
U.S. Dept. Comm. Off. Tech. Serv. No. \textbf{PB32719},
Washington, DC, 1943.


\bibitem{Neumann1}
J. von Neumann,
{Oblique reflection of shocks}, {\itshape Explo. Res. Rep.} \textbf{12}, Navy
Department, Bureau of Ordnance, Washington, DC, 1943.

\bibitem{Neumann2}
J. von Neumann,
Refraction, intersection, and reflection of shock waves,
{\em NAVORD Rep.} \textbf{203-45}, Navy Department, Bureau of Ordnance,
Washington, DC, 1945.



\bibitem{Neumann}
J. von Neumann,
\textit{Collected Works}, Vol. \textbf{6}, Pergamon: New York, 1963.


\bibitem{Neumann10}
J. von Neumann,
Discussion on the existence and uniqueness or multiplicity of solutions of
the aerodynamical equation [Reprinted from MR0044302 (1949)],
\textit{Bull. Amer. Math. Soc. {\rm (}N.S.{\rm )}} \textbf{47} (2010), 145--154.

\bibitem{Whitham}
G. B. Whitham,
\textit{Linear and Nonlinear Waves},
John Wiley \& Sons, Inc.: New York,  1974.

\bibitem{WC}
P. Woodward and P. Colella,
{The numerical simulation of two-dimensional fluid flow with
strong shocks}, \textit{J. Comp. Phys.} \textbf{54} (1984), 115--173.

\bibitem{XY}
Z. Xin and H. Yin,
{Transonic shock in a nozzle: two-dimensional case}, \textit{Comm.
Pure Appl. Math.} \textbf{58} (2005), 999--1050.

\bibitem{YinZhou}
H. Yin and C. Zhou,
\newblock On global transonic shocks for the steady supersonic Euler flows past sharp 2-D wedges.
\newblock \textit{J. Diff. Eqs.} \textbf{246} (2009), 4466--4496


\bibitem{Yuan1}
H. Yuan,
{On transonic shocks in two-dimensional variable-area ducts
for steady Euler system}, \textit{SIAM J. Math. Anal.}
\textbf{38} (2006), 1343--1370.

\bibitem{ZZ}
T. Zhang and Y. Zheng,
{Conjecture on the structure of solutions of the Riemann problem
for two-dimensional gas dynamics}, \textit{SIAM J. Math. Anal.}
\textbf{21} (1990), 593--630.


\bibitem{Zh2}
Y.-Q. Zhang,
{Steady supersonic flow past an almost straight wedge with
large vertex angle}, \textit{J. Diff.  Eqs.} \textbf{192} (2003), 1--46.


\bibitem{Zhe}
Y. Zheng,
{\em Systems of Conservation Laws{\rm :} Two-Dimensional Riemann
Problems}, Birkh\"{a}user: Boston, 2001.

\end{thebibliography}
\end{document}